\documentclass[12pt]{amsart}

\usepackage{amsmath, amsthm, amssymb}
\input xy
\xyoption{all}
\usepackage{mathrsfs}
\usepackage{enumerate}
\usepackage{color}
\usepackage{tikz}
\usetikzlibrary{matrix}
\usepackage{caption}
\usepackage{subcaption}

\newtheorem{dummy}{}[section]
\newtheorem{theorem}[dummy]{Theorem}

\newtheorem{lemma}[dummy]{Lemma}

\newtheorem{conjecture}[dummy]{Conjecture}
\theoremstyle{definition}
\newtheorem{definition}[dummy]{Definition}

\newtheorem{remark}[dummy]{Remark}

\newcommand{\Z}{\ensuremath{\mathbb{Z}}}
\newcommand{\Q}{\ensuremath{\mathbb{Q}}}
\newcommand{\C}{\ensuremath{\mathbb{C}}}
\renewcommand{\P}{\ensuremath{\mathbb{P}}}

\newcommand{\M}{\ensuremath{\overline{\mathcal{M}}}}
\renewcommand{\L}{\ensuremath{\mathcal{L}}}
\renewcommand{\O}{\ensuremath{\mathcal{O}}}
\newcommand{\cC}{\ensuremath{\mathcal{C}}}
\newcommand{\cP}{\ensuremath{\mathcal{P}}}
\newcommand{\ev}{\mathrm{ev}}
\newcommand{\vir}{\mathrm{vir}}
\newcommand{\vdim}{\ensuremath{\textrm{vdim}}}
\renewcommand{\t}{\ensuremath{\mathbf{t}}}
\renewcommand{\d}{\ensuremath{\partial}}
\newcommand{\m}{\ensuremath{\vec{m}}}

\newcommand{\one}{\ensuremath{\mathbf{1}}}

\newcommand{\bbeta}{\ensuremath{\vec{\beta}}}
\newcommand{\tw}{\mathrm{tw}}
\newcommand{\GIT}{\ensuremath{\mathbin{/\mkern-6mu/}}}
\newcommand{\Hom}{\ensuremath{\textrm{Hom}}}
\renewcommand{\H}{\ensuremath{\mathcal{H}}}
\newcommand{\G}{\ensuremath{\mathcal{G}}}
\newcommand{\Ze}{\ensuremath{Z^{\epsilon}}}
\newcommand{\Xe}{\ensuremath{X^{\epsilon}}}

\renewcommand{\i}{\ensuremath{\textrm{i}}}
\newcommand{\n}{\ensuremath{n'}}
\newcommand{\ri}{\overline{\mathcal I}}
\newcommand{\WC}{\mathrm{WC}}
\newcommand{\ct}{\mathrm{ct}}
\newcommand{\Xect}{X^{\epsilon,\ct}}
\newcommand{\Xecto}{X^{\epsilon,\ct-1}}
\newcommand{\CR}{\mathrm{CR}}
\newcommand{\fM}{\mathfrak{M}}

\usepackage{geometry}
\usepackage{color}

\begin{document}

\title[Higher-genus wall-crossing in the GLSM]{Higher-genus wall-crossing in the gauged linear sigma model}
\author{Emily Clader, Felix Janda, and Yongbin Ruan \\ With an appendix by Yang Zhou}
\address{Department of Mathematics, San Francisco State University, 1600 Holloway Avenue, San Francisco, CA 94132, USA}
\email{eclader@sfsu.edu}
\address{Department of Mathematics, University of Michigan, 2074 East Hall, 530 Church Street, Ann Arbor, MI 48109, USA}
\email{janda@umich.edu}
\address{Department of Mathematics, University of Michigan, 2074 East Hall, 530 Church Street, Ann Arbor, MI 48109, USA}
\email{ruan@umich.edu}
\address{Department of Mathematics, Stanford University, Building 380, Stanford, California 94305, USA}
\email{yangzhou@stanford.edu}


\begin{abstract}
We introduce a technique for proving all-genus wall-crossing formulas in the gauged linear sigma model as the stability parameter varies, without assuming factorization properties of the virtual class.  Implementing this technique to the gauged linear sigma model associated to a complete intersection in weighted projective space, we obtain a uniform proof of the wall-crossing formula in both the geometric and the Landau--Ginzburg phase.
\end{abstract}

\maketitle

\section{Introduction}

The gauged linear sigma model (GLSM) was introduced in the physics literature by Witten \cite{Witten} in the early 1990's as an example of a two-dimensional quantum field theory whose target is a complete intersection in a GIT quotient.  Since then, it has been extensively studied in physics (see \cite{HIV, HV} for references) and has become a powerful tool in understanding many aspects of both physics and geometry.

The mathematical development of the GLSM as a curve-counting theory was carried out by Fan, Jarvis, and the third author in \cite{FJRGLSM}, where the foundations were laid for the theory in the ``compact-type" subspace.  The resulting theory is divided into chambers (or ``phases") by variation of GIT on the target geometry.  In the geometric phase, the GLSM recovers the quasimap theory developed by Ciocan-Fontanine, Kim, and Maulik \cite{CFKM}, while in the Landau--Ginzburg phase, it recovers Fan--Jarvis--Ruan--Witten (FJRW) theory \cite{FJR1} when the target is a hypersurface in weighted projective space.

Our interest in the GLSM is motivated by the celebrated Landau--Ginzburg/Calabi--Yau (LG/CY) correspondence.  Originally framed mathematically as a connection between the Gromov--Witten theory of a hypersurface and the FJRW theory of its defining polynomial, the correspondence can be re-cast in this new language as a wall-crossing (or ``phase transition") between different chambers of the GLSM.  Two types of walls are relevant here: between phases there are chamber walls, and within each phase, there are walls as one varies the choice of a stability parameter $\epsilon$.  So far, very little is understood about the transition across the chamber walls, though a general conjecture posits that it should involve analytic continuation of generating functions.  The $\epsilon$-wall-crossing, on the other hand, has been made remarkably explicit by Ciocan-Fontanine and Kim, as we explain below; this is the subject of the current work.  Our ultimate goal is to combine the two types of wall-crossings to yield a proof of the LG/CY correspondence in all genus.  In genus zero, this program has been carried out when the target is a hypersurface (through the combined work of Ciocan-Fontanine--Kim \cite{CFKZero}, Ross and the third author \cite{RR}, and Chiodo and the third author \cite{CR}), and in genus one, it has been carried out for the quintic hypersurface without marked points (through the work of Kim--Lho \cite{KL} and Guo--Ross \cite{GR1, GR2}).

\subsection{Statement of results}

For a complete intersection $Y$ in projective space, the definition of quasimaps introduced by Ciocan-Fontanine, Kim, and Maulik generalizes the notion of stable maps to $Y$ via the additional datum of a positive rational number $\epsilon$.  When $\epsilon \rightarrow \infty$, quasimaps coincide with the usual stable maps, so quasimap theory recovers Gromov--Witten theory.  When $\epsilon \rightarrow 0$, on the other hand, quasimaps become stable quotients \cite{MOP}, and the resulting theory is thought to correspond to the mirror B-model of $Y$ \cite{CZ, CFKZero, CFKBigI}.

Quasimap theory changes only at certain discrete values of $\epsilon$, so there is a wall-and-chamber structure on the space of stability parameters.  Ciocan-Fontanine and Kim proved a wall-crossing formula in \cite{CFKZero} exhibiting how the genus-zero theory varies with $\epsilon$, and in \cite{CFKHigher, CFKMirror}, they formulated the following conjecture in all genus:

\begin{conjecture}[See \cite{CFKMirror}]
  \label{conj:main}
  Let $Y$ be a complete intersection in projective space, and fix
  $g, n \geq 0$.
  Then
\[\sum_{\beta} q^\beta [\M^{\epsilon}_{g, n}(Y,\beta)]^{\vir} = \sum_{\beta_0, \beta_1, \ldots, \beta_k} \frac{q^{\beta_0}}{k!} b_{\bbeta*}c_* \left(\prod_{i = 1}^k q^{\beta_i} \ev_{n + i}^*(\mu_{\beta_i}^\epsilon(-\psi_{n + i})) \cap [\M^\infty_{g, n + k}(Y,\beta_0)]^{\vir} \right),\]
  where $\mu^{\epsilon}_{\beta}(z)$ are certain coefficients of the
  $I$-function of $Y$, $b_{\bbeta}$ is a morphism that converts marked
  points to basepoints, and $c$ is the natural contraction morphism
  from $\infty$-stable to $\epsilon$-stable quasimaps.
\end{conjecture}

Ciocan-Fontanine and Kim proved this conjecture in \cite{CFKMirror}, using virtual push-forward techniques and MacPherson's graph construction.  A different proof was given in \cite{CJR} under the assumption that $n \geq 1$, using localization on a ``twisted graph space", which has the advantage that it can be adapted to the more general context of the GLSM.

To explain this generalized setting, we recall that the GLSM depends on the choice of (1) a GIT quotient
\[X_{\theta} = [V\!\GIT_{\theta}\!G],\] in which $V$ is a complex
vector space, $G \subset \text{GL}(V)$, and $\theta$ is a character of
$G$; (2) a polynomial function $W\colon X_{\theta} \rightarrow \C$
known as the ``superpotential''; and (3) an action of $\C^*$ on $V$
known as the ``$R$-charge''.
Associated to any such choice of input data and any stability
parameter $\epsilon$, there is a state space $\H_{\theta}$ and a moduli space
$Z^{\epsilon,\theta}_{g,n,\beta}$ parameterizing $\epsilon$-stable Landau--Ginzburg
quasimaps to the critical locus of $W$ (see Section~\ref{sec:defs}
below).
Under the technical requirement that the theory admits a ``good
lift'', $Z^{\epsilon,\theta}_{g,n,\beta}$ is compact and admits a virtual fundamental
class, which can be paired with elements of the compact-type subspace
$\H^{\ct}_{\theta} \subset \H_{\theta}$ to yield invariants.

If $Y$ is a nonsingular complete intersection in weighted projective space $\P(w_1, \ldots, w_M)$ defined by the vanishing of polynomials $F_1, \ldots, F_N$ of degrees $d_1, \ldots, d_N$, then one obtains a GIT quotient by letting $G = \C^*$ act on $V= \C^{M+N}$ with weights 
\[(w_1, \ldots, w_M, -d_1, \ldots, -d_N).\]
Choosing $\theta \in \Hom(G,\C^*) \cong \Z$ to be any positive character, the resulting GIT quotient is
\[X_+ := \bigoplus_{j=1}^N \O_{\P(w_1, \ldots, w_M)}(-d_j).\]
Denote the coordinates of $X_+$ by $(x_1, \ldots, x_M, p_1, \ldots, p_N)$, and define a superpotential by
\begin{equation*}
  W(x_1, \ldots, x_M, p_1, \ldots, p_N) := \sum_{j=1}^N p_jF_j(x_1, \ldots, x_M).
\end{equation*}
Finally, define an $R$-charge by letting $\C^*$ act with weight $1$ on
the $p$-coordinates.
Then the critical locus of $W$ is simply the complete intersection
$Y$, the existence of a good lift is automatic, and the moduli space of the GLSM reduces to the moduli space of ordinary quasimaps to $Y$:
\[Z^{\epsilon,+}_{g,n,\beta} = \M_{g,n}^{\epsilon}(Y,\beta).\]
The state space $\H_+$ is isomorphic to the cohomology of $Y$, and the
compact-type subspace consists of classes pulled back from the ambient
$\P(w_1, \ldots, w_M)$.
Thus, one obtains a theory that is thought to coincide with the quasimap theory of $Y$ with ambient insertions.\footnote{In order to know that the two theories indeed coincide, one would need to match the construction of the virtual cycle in the GLSM, defined via Kiem--Li's cosection technique, with the construction of the virtual cycle in quasimap theory by Ciocan-Fontanine--Kim--Maulik.  This has been carried out (on the level of correlators) for the quintic hypersurface by Chang--Li \cite{CL}, and their proof can be generalized to complete intersections in $\P^{N-1}$.  See Remark~\ref{rem:geomvc} below for further discussion.}

On the other hand, when the character $\theta$ is negative, the resulting GIT quotient is
\[X_- := \bigoplus_{i=1}^M \O_{\P(d_1, \ldots, d_N)}(-w_i)\]
and the critical locus of the above $W$ is the zero section.  With the $R$-charge as above, the moduli space $Z^{\epsilon,-}_{g,n,\beta}$ parameterizes marked curves equipped with a line bundle $L$ and a section
\[\vec{p} \in \Gamma\left( \bigoplus_{j=1}^N (L^{\otimes - d_j} \otimes \omega_{\log})\right)\]
satisfying certain stability conditions depending on the parameter $\epsilon$.  In particular, when $N=1$, one has a section $p \in \Gamma(L^{\otimes -d} \otimes \omega_{\log})$, and for $\epsilon \gg 0$, the stability condition enforces that $p$ is nowhere vanishing; thus, $p$ trivializes $L^{\otimes -d} \otimes \omega_{\log}$, and the resulting moduli space of $d$-spin structures is precisely the moduli space of FJRW theory.

The assumption that the theory admits a good lift, explained in detail
in \cite{FJRGLSM}, amounts in this setting to the requirement that
$d_1 = \cdots = d_N$, in which case the theory is known in the physics
literature as a ``hybrid model'' and was developed mathematically in
the narrow sectors (a subset of the state space properly contained in
$\H^{\ct}_-$) by the first author \cite{Clader}.

Conjecture~\ref{conj:main} can be directly generalized to this setting, and our main theorem is a verification that the conjecture holds:

\begin{theorem}[See Theorem~\ref{thm:main}]
\label{thm:1}
Let $Y \subset \P(w_1, \ldots, w_M)$ be a nonsingular complete intersection defined by the vanishing of a collection of polynomials of degrees $d_1, \ldots, d_N$, where $w_i|d_j$ for all $i$ and $j$.  Fix $\theta \in \{+,-\}$ and suppose that the associated GLSM admits a good lift.

For any $g$ and $n$ and a tuple of insertions 
$\alpha_1, \ldots, \alpha_n \in \H^{\ct}_{\theta}$, one has
\begin{align*}
&\sum_{\beta} q^\beta  \left(\prod_{a=1}^n \ev_a^*(\alpha_a) \cap [Z^{\epsilon,\theta}_{g,n,\beta}]^{\vir}\right) =\\
&\sum_{\beta_0, \beta_1, \ldots, \beta_k} \frac{q^{\beta_0}}{k!} b_{\bbeta*}c_* \left(\prod_{a=1}^n \ev_a^*(\alpha_a)\prod_{i = 1}^k q^{\beta_i} \ev_{n + i}^*(\mu_{\beta_i}^\epsilon(-\psi_{n + i})) \cap [Z^{\infty,\theta}_{g,n+k,\beta_0}]^{\vir} \right).
\end{align*}
\end{theorem}

\begin{remark}
In particular, taking $\theta = +$ and $w_1 = \cdots = w_M = 1$, the proof of Theorem~\ref{thm:1} reproduces the results of \cite{CJR} for ambient insertions, assuming the equivalence of quasimap theory with the positive phase of the GLSM mentioned above.  For this reason, we focus in what follows on the case where $\theta = -$, remarking where necessary on the appropriate modifications for the positive (that is, the ``geometric") phase.
\end{remark}

The basic idea of the proof of Theorem~\ref{thm:1} is the same as the authors' technique in the geometric phase, as presented in \cite{CJR}.  Namely, we construct a larger moduli space with a $\C^*$-action in which the theories at $\epsilon=\infty$ and arbitrary $\epsilon$ arise as fixed loci.  (This larger moduli space is closely related to the space of mixed-spin $p$-fields considered by Chang--Li--Li--Liu \cite{CLLL, CLLL2}.)  A direct application of the proof in \cite{CJR}, however, would require a factorization property of the virtual class along strata corresponding to nodal curves, and the verification of such a property is an unsolved issue dating back to the early days of FJRW theory; the problem is that, while the theory is only defined for compact-type insertions, decomposition at nodes may involve insertions that are not of compact type.  In ongoing research, Ciocan-Fontanine--Favero--Gu\'{e}r\'{e}--Kim--Shoemaker \cite{CFFGKS} are working toward a resolution of this problem via the theory of matrix factorizations, and an analytic approach has also been proposed by Tian--Xu \cite{TX1, TX2} and Fan--Jarvis--Ruan \cite{FJRAnalytic}.

Rather than awaiting the completion of the general theory in the non-compact-type case, though (since the compact-type theory is sufficient for the applications in which we are interested), we present in this work a technique that avoids the necessity of factorization along nodes.  The idea is to introduce a new moduli space for each boundary stratum in $Z^{\epsilon,\theta}_{g,n,\beta}$ and to carefully choose the discrete data in the twisted graph space so that only very particular degenerations occur in the contributing fixed loci.  We hope that this strategy may be useful in other situations where the factorization of the GLSM virtual class has been an obstacle.

As in the geometric phase, the proof of Theorem~\ref{thm:1} initially requires that $n \geq 1$, since it relies crucially on the fact that a certain localization expression changes in a nontrivial way when an insertion is varied.  However, in the appendix (written by Yang Zhou), the result for $n=0$ is deduced from the result for $n \geq 1$ by leveraging a wall-crossing between ``heavy points" and ``light points."  These ideas were developed in work of Yang Zhou \cite{Z}, in which an alternative proof of Theorem~\ref{thm:1} is given in the case where $Y$ is a hypersurface.

\subsection{Plan of the paper}

In Section~\ref{sec:defs}, we review the necessary theory of the GLSM in the setting of interest here.  We state our main theorem in precise form in Section~\ref{sec:statements}, and we state the progressive refinements of the theorem that are necessary for the proof.  We introduce the twisted graph space in Section~\ref{sec:twistedgs}, where we explicitly calculate the fixed-locus contributions to the localization formula.  We also review the setting where the target is a point, in which case the twisted graph space specializes to $\M_{g,n}(\P^1,\delta)$ and was studied in detail in \cite{CJR}.  Finally, in Section~\ref{sec:proof}, we present the proof of Theorem~\ref{thm:1}.  The structure of the proof is to compute, via localization, the pushforward from the twisted graph space to $Z^{\epsilon,\theta}_{g,n,\beta}$ of a certain difference of cohomology classes.  The contributions from degree-zero components are closely related to the calculations for a point, and using these calculations, we show that the class in question changes by an irrational function of the equivariant parameter when an insertion is varied.  This implies that the class vanishes, and from here, we deduce the wall-crossing theorem when $n \geq 1$.  In the appendix, the notion of ``light" marked points is introduced, and the $n=0$ case of the wall-crossing is deduced from the $n \geq 1$ case by proving an auxiliary wall-crossing that relates heavy to light markings and then applying the dilaton and divisor equations in the light-marking setting.

\subsection{Acknowledgments}

The authors are grateful to Ionu\c{t} Ciocan-Fontanine, Bumsig Kim, Dustin Ross, and Yang Zhou for many useful conversations and comments.  The third author was partially supported by NSF grant DMS 1405245 and NSF FRG grant DMS 1159265.

\section{Definitions and setup}
\label{sec:defs}

We review the definition of the gauged linear sigma model (GLSM) for the Landau--Ginzburg phase in the setting of Theorem~\ref{thm:1}.  This is a very special case of the general construction of the GLSM by Fan, Jarvis, and the third author.  When $\theta$ is negative and $\epsilon$ is sufficiently large, it recovers the hybrid model defined in \cite{Clader}.

\subsection{Gauged linear sigma model}

Throughout what follows, we fix a nonsingular complete intersection
\[ Y := \{F_1= \cdots = F_N = 0\} \subset \P(w_1, \ldots, w_M),\]
where $F_1, \ldots, F_N$ are polynomials of degrees $d_1, \ldots, d_N$ and $w_i|d_j$ for all $i$ and $j$.

The GLSM, in general, depends on three pieces of input data: a GIT quotient $X_{\theta}=[V \GIT_{\theta} G]$, a polynomial function $W\colon X \rightarrow \C$ known as the superpotential, and an action of $\C^*$ on $V$ known as the $R$-charge.  In our case,
\[X_{\theta}:= [(\C^M \times \C^N) \GIT_{\theta} G],\]
where
\[G:= \{(g^{w_1}, \ldots, g^{w_M}, g^{-d_1}, \ldots, g^{-d_N}) \; | \; g \in \C^*\}\cong \C^*\]
acts diagonally on $V:= \C^M \times \C^N$, and there are two phases depending on whether the character $\theta \in \Hom_{\Z}(\C^*, \C^*) \cong \Z$ is positive or negative.  Denoting the coordinates on $V$ by $(x_1, \ldots, x_M, p_1, \ldots, p_N)$, the unstable locus when $\theta$ is positive is $\{x_1 = \cdots = x_M = 0\}$, so we have
\[X_+ = \frac{(\C^M \setminus \{0\}) \times \C^N}{\C^*} \cong \bigoplus_{j=1}^N \O_{\P(w_1, \ldots, w_M)}(-d_j).\]
Similarly, when $\theta$ is negative, the unstable locus of $\theta$ is $\{p_1= \cdots = p_N = 0\}$, so
\[X_-  = \frac{\C^M \times (\C^N \setminus \{0\})}{\C^*}\cong \bigoplus_{i=1}^M \O_{\P(d_1, \ldots, d_N)}(-w_i).\]
The superpotential in either phase is
\[W(x_1, \ldots, x_M, p_1, \ldots, p_N) := \sum_{i=1}^N p_j F_j(x_1, \ldots, x_M), \]
and the $R$-charge acts by multiplication on the $p$-coordinates.  Using that $Y$ is nonsingular, one can check that the critical locus $Z_+$ of $W$ when $\theta$ is positive is the complete intersection $Y$ inside the zero section of $X_+$, and the critical locus $Z_-$ when $\theta$ is negative is the entire zero section in $X_-$.

From here forward, we restrict to the negative phase of the GLSM unless otherwise specified, and we assume that all of the degrees $d_j$ are equal.  To ease notation, we make use of the abbreviations
\[d:= d_1 = \cdots = d_N,\]
\[X:= X_- = \bigoplus_{i=1}^M \O_{\P(d, \ldots, d)}(-w_i),\]
and
\[Z:=Z_- = \P(d,\ldots, d) \subset X.\]

\subsection{State space}
\label{sec:ss}

The state space of the GLSM is defined as the relative Chen--Ruan cohomology group
\begin{equation}
  \label{ss}
  \H := H^*_{\CR}(X, W^{+\infty};\C),
\end{equation}
in which $W^{+\infty}$ is a Milnor fiber of $W$---that is, $W^{+\infty} = W^{-1}(A)$ for a sufficiently large real number $A$.  Thus, the state space decomposes into summands indexed by the components of the inertia stack $\mathcal{I}X$, which are labeled by elements $g \in G$ with nontrivial fixed locus $\text{Fix}(g) \subset X$.

More concretely, the only $(g^{w_1}, \ldots, g^{w_M}, g^{-d}, \ldots, g^{-d})\in G$ with nontrivial fixed locus are those for which $g^d = 1$.  An element $(\vec{x}, \vec{p}) \in \text{Fix}(g)$ must have $x_i=0$ whenever $g^{w_i}\neq 1$, but there is no constraint on the $x_i$ for which $g^{w_i}=1$.  Thus, denoting $\P(\vec{d}) :=\P(d,\ldots, d)$ and
\begin{equation*}
\label{eq:F(g)}
F(g):= \{ i \; | \; g^{w_i} = 1\} \subset \{1, \ldots, M\},
\end{equation*}
we set
\[X_g:= \bigoplus_{i \in F(g)} \O_{\P(\vec{d})}(-w_i) \subset X\]
and $W_g:= W|_{X_g}$.  Then
\[\H = \bigoplus_{g \in \Z_d} H^*(X_g, W_g^{+\infty};\C).\]
An element $g \in \Z_d$ (or its corresponding component of $\H$) is referred to as {\it narrow} if $\text{Fix}(g)$ is compact, which amounts to requiring that $F(g) = \emptyset$.  In other words, if we set
\begin{equation}
\label{eq:nar}
\text{nar}:= \left\{m \in \left.\left\{0,\frac{1}{d}, \cdots, \frac{d-1}{d}\right\}\right| \not\exists \; i \text{ such that } mw_i \in \Z\right\},
\end{equation}
then the narrow sectors are indexed by $g = e^{2\pi i m}$ with $m \in \text{nar}$.  Narrow sectors have $W|_{\text{Fix}(g)} \equiv 0$, so their contribution to $\H$ is simply $H^*(\P(\vec{d})) \cong H^*(\P^{N-1})$.

For any $g \in \Z_d$, there is a natural map
\[\eta_g\colon H^{k-2|F(g)|}(X_g) \rightarrow H^k(X_g, W_g^{+\infty}).\]
To define $\eta_g$, we set $X_g^{\text{ct}} \cong \P(\vec{d})$ to be the zero section inside $X_g$.  The map $\eta_g$ is defined as the composition
\begin{equation}
\label{eq:etag}
H^{k-2|F(g)|}(X_g) \xrightarrow{\sim} H^{k-2|F(g)|}(X_g^{\text{ct}}) \xrightarrow{\sim} H^k(X_g, X_g \setminus X_g^{\text{ct}}) \rightarrow H^k(X_g, W_g^{+\infty}),
\end{equation}
where the second map is the Thom isomorphism and the remaining maps are induced by the inclusions.  Together, the $\eta_g$ define a homomorphism
\begin{equation}
\label{eq:eta}
\eta\colon H^*_{\CR}(X) \rightarrow \H,
\end{equation}
and we define the \emph{compact-type state space} as
\begin{equation*}
  \H^{\ct} := \text{image}(\eta) \subset \H.
\end{equation*}
Note that $\eta_g$ is an isomorphism when $g$ is narrow, so the
compact-type state space contains all of the narrow sectors.

Instead of working with $\H^{\ct}$, in what follows we take insertions from 
\begin{equation*}
  \widetilde\H := H^*_{\CR}(X),
\end{equation*}
which surjectively maps to $\H^{\ct}$ via $\eta$.  If one wishes to take insertions from $\H^{\ct}$ instead, then one must first lift them to $H^*_{\CR}(X)$.  We conjecture that the resulting correlators are independent of the choice of lift; see Lemma~\ref{lem:vanishing} below.

\begin{remark}
\label{rem:geomss}
In the geometric chamber, there is an isomorphism
\[\H_+ := H^*_{\CR}(X_+, W^{+\infty}; \C) \cong H^*_{\CR}(Y)\]
(see \cite[Proposition 3.4]{CN}).  This has a twisted sector
\begin{equation*}
\label{eq:Xg+}
X_{g,+} = \bigoplus_{j=1}^N \O_{\P_g}(-d)
\end{equation*}
for any $g$ such that $g^{w_i} = 1$ for some $i$, in which $\P_g \subset \P(w_1, \ldots, w_M)$ is the sub-projective space spanned by the coordinates $x_i$ for which $g^{w_i} = 1$.  The zero section is $X_{g,+}^{\text{ct}} \cong \P_g$ in this case, and if $Y_g:= \P_g \cap Y$, then $\eta_g$ can be identified with the restriction $H^*(\P_g) \rightarrow H^*(Y_g)$.  Thus, the compact-type state space $\H_{+}$ consists of the ambient classes in each sector of $H^*_{\CR}(Y)$ on the geometric side, and inserting classes
from $\widetilde\H_+:= H^*_{\CR}(X_+)$ simply means lifting an ambient class
$\iota^*\alpha$ (in which $\iota\colon Y_g \rightarrow \P_g$ is the
inclusion) to an insertion of $\alpha$.
\end{remark}

\subsection{Moduli space}

The definition of the GLSM moduli space, which depends on the choice of a stability parameter $\epsilon$, is based on the notion of quasimaps introduced by Ciocan-Fontanine, Kim, and Maulik \cite{CFKM} and studied extensively by Ciocan-Fontanine and Kim \cite{CFK3, CFKZero, CFKBigI, CFKMirror, CFKHigher}.

Fix a genus $g$, a degree $\beta \in \Z$, a nonnegative integer $n$, and a positive rational number $\epsilon$.

\begin{definition}
\label{def:quasimap}
A \emph{prestable Landau--Ginzburg quasimap} to $Z$ consists of an $n$-pointed prestable orbifold curve $(C; q_1, \ldots, q_n)$ of genus $g$ with nontrivial isotropy only at special points, an orbifold line bundle $L$ on $C$, and a section
\[\vec{p} = (p_1, \ldots, p_N) \in \Gamma((L^{\otimes -d} \otimes \omega_{\log})^{\oplus N})\]
(where
\[\omega_{\log}:= \omega_C \otimes \O_C([q_1] + \cdots + [q_n]),\]
is the logarithmic dualizing sheaf) such that the zero set of
$\vec{p}$ (that is, the set of points $q \in C$ such that
$p_1(q) = \cdots = p_N(q) = 0$) is finite.
An \emph{$\epsilon$-stable Landau--Ginzburg quasimap} is a prestable
Landau--Ginzburg quasimap, satisfying the following further
conditions:
\begin{itemize}
\item {\it Representability}:  For every $q \in C$ with isotropy group $G_q$, the homomorphism $G_q \rightarrow \C^*$ giving the action of the isotropy group on the bundle $L$ is injective.
\item {\it Nondegeneracy}: The zero set of $\vec{p}$ is disjoint from
  the marked points and nodes of $C$, and for each zero $q$ of
  $\vec{p}$, the order of the zero (that is, the common order of
  vanishing of $p_1, \ldots, p_N)$ satisfies
\begin{equation}
\label{nondegen}
\text{ord}_q(\vec{p}) \leq \frac{1}{\epsilon}.
\end{equation}
(Zeroes of $\vec{p}$ are referred to as {\it basepoints} of the quasimap.)
\item {\it Stability}:  The $\Q$-line bundle
\begin{equation}
\label{stab}
(L^{\otimes -d} \otimes \omega_{\log})^{\otimes \epsilon} \otimes \omega_{\log}
\end{equation}
is ample.
\end{itemize}
The {\it degree} of a quasimap is defined as
\begin{equation}
\label{eq:beta}
\beta:=\deg(L^{\otimes -d} \otimes \omega_{\log}).
\end{equation}
\end{definition}

Notice that when $N=1$ and $\epsilon > 2$, condition \eqref{nondegen} implies that the section $p_1$ is nowhere-vanishing, so it gives a trivialization $L^{\otimes -d} \otimes \omega_{\log} \cong \O_C$.  In this case, the degree must be zero, and condition \eqref{stab} amounts to the requirement that $(C; q_1, \ldots, q_n)$ be a stable orbifold curve.  The definition of an $\epsilon$-stable quasimap, then, recovers the notion of a $d$-spin curve.  On the other hand, when $\epsilon \leq \frac{1}{\beta}$, condition \eqref{nondegen} puts no restriction on the orders of the basepoints, and \eqref{stab} is equivalent to imposing the analogous requirement for all $\epsilon > 0$.  The resulting moduli space is analogous to the moduli space of stable quotients \cite{MOP}.

\begin{remark}
\label{rem:rational}
An alternative way to view the choice of stability parameter $\epsilon$ is to replace, in the definition of the GLSM, the character $\theta \in \Hom_{\Z}(\C^*, \C^*)$ by the rational character $\theta \cdot (d\epsilon) \in \Hom_{\Z}(\C^*, \C^*) \otimes_{\Z} \Q \cong \Q$; see \cite[Section 2]{CFKBigI}.   We return to this perspective later when defining the twisted graph space.
\end{remark}

The properness of the moduli space of the GLSM in general is a subtle question that requires the existence of a ``good lift" of the character $\theta$.  However, in the hybrid model, the trivial lift is good (see \cite[Example 7.2.2]{FJRGLSM}), so we have the following result:

\begin{theorem}[Fan--Jarvis--Ruan \cite{FJRGLSM}]
There is a proper Deligne--Mumford stack $Z^{\epsilon}_{g,n,\beta}$ parameterizing genus-$g$, $n$-pointed, $\epsilon$-stable Landau--Ginzburg quasimaps of degree $\beta$ to $Z$ up to isomorphism.
\end{theorem}

\begin{remark}
To put things another way, one can view an $\epsilon$-stable Landau--Ginzburg quasimap to $Z$ as as a section
\[(\vec{x}, \vec{p}) \in \Gamma\left(\bigoplus_{i=1}^M L^{\otimes w_i} \oplus \bigoplus_{j=1}^N (L^{\otimes -d_j} \otimes \omega_{\log})\right)\]
whose image lies in the affine cone over the critical locus $Z$ and whose order of contact with the unstable locus is bounded.  This perspective is important for defining the virtual cycle, and it also points to how one must modify Definition~\ref{def:quasimap} for the positive phase.  Namely, Landau--Ginzburg quasimaps to $Z_+$ consist of sections $(\vec{x},\vec{p})$ as above whose image lies in the affine cone over the critical locus of $Z_+$, which precisely recovers the definition of quasimaps to the complete intersection $Y$ as in \cite{CFKM}.  Replacing $\vec{p}$ by $\vec{x}$ in \eqref{nondegen} and $L^{\otimes -d} \otimes \omega_{\log}$ by $L$ in both \eqref{stab} and \eqref{eq:beta} yields the definition of $\epsilon$-stability for quasimaps.  Thus, $Z^{\epsilon,+}_{g,n,\beta}$ is the moduli space $\M^{\epsilon}_{g,n}(Y,\beta)$ of (ordinary) quasimaps defined in \cite{CFKM}.
\end{remark}

\subsection{Multiplicities and evaluation maps}

Recall that the {\it multiplicity} of an orbifold line bundle $L$ at a point $q \in C$ with isotropy group $\Z_r$ is defined as the number $m \in \Q/\Z$ such that the canonical generator of $\Z_r$ acts on the total space of $L$ in local coordinates around $q$ by
\[(x,v) \mapsto \left(e^{2\pi \i \frac{1}{r}}x, e^{2\pi \i m}v\right).\]

In our case, we can take the multiplicities to lie in the set $\{0, \frac{1}{d}, \ldots, \frac{d-1}{d}\}$, and for a tuple $\m = (m_1, \ldots, m_n)$ with $m_i \in \{0, \frac{1}{d}, \ldots, \frac{d-1}{d}\}$, we denote by
\[Z^{\epsilon}_{g,\m,\beta} \subset Z^{\epsilon}_{g,n,\beta}\]
the open and closed substack consisting of quasimaps for which the multiplicity of $L$ at $q_i$ is $m_i$.  We occasionally wish to leave some multiplicities undetermined, so we denote by
\[Z^{\epsilon}_{g,\m+k,\beta} \subset Z^{\epsilon}_{g,n+k,\beta}\]
the substack on which the multiplicity of $L$ at $q_i$ is $m_i$ for $1 \leq i \leq n$, while the last $k$ marked points are allowed any multiplicity.

A crucial feature of the multiplicities is that they determine the relationship between $L$ and its pushforward $|L|$ to the coarse underlying curve.  Specifically, suppose that $C' \subset C$ is an irreducible component of $C$ with special points $\{q_k\}$ at which the multiplicities of $L$ are $\{m_k\}$.  Then, if $\rho\colon C' \rightarrow |C'|$ is the natural map to the coarse underlying curve, we have 
\begin{equation*}
\label{coarse}
L = \rho^*|L| \otimes \O_{C'} \left(\sum_{k} m_k [\overline{q_k}]\right),
\end{equation*}
for the degree-$1$ divisors $\overline{q_k}$ pulled back from the coarse underlying curve $|C|$.
Applying this equation in the case where $C$ is smooth yields a compatibility condition on the multiplicities, since the degree of $\rho^*|L|$ is an integer.  Namely, we have:
\begin{equation}
\label{eq:compatibility}
\frac{-\beta+2g-2+n}{d} - \sum_{i=1}^n m_i \in \Z. 
\end{equation}

Note that equation \eqref{eq:compatibility} is independent of $m_i$ if and only if $m_i = \frac{1}{d}$, so this is the only case in which there is a forgetful map on $\Ze_{g,\m,\beta}$ forgetting $q_i$ and its orbifold structure.  Thus, the role of the unit in the $\epsilon \gg 0$ hybrid theory---in particular, in the string and dilaton equation---is played by the fundamental class in the narrow sector $H^*(\P^{N-1}) \subset \widetilde{\H}$ indexed by $e^{2\pi i \frac{1}{d}} \in \Z_d$.  For this reason, we denote this element of $\widetilde{\H}$ by $\one$ in what follows.

\begin{remark}
The same string and dilaton equations do {\it not} hold in the theory of $\epsilon$-stable Landau--Ginzburg quasimaps for more general $\epsilon$.  One can deduce the appropriate modification from Theorem \ref{thm:1}; see the discussion in \cite[Section 3.4]{CFKHigher}.
\end{remark}

\begin{remark}
\label{rem:geommults}
In the geometric chamber, our convention is that $\beta = \deg(L)$, so the analogue of condition \eqref{eq:compatibility} is $\beta - \sum_{i=1}^n m_i \in \Z$.  In particular, there is a forgetful map on $Z^{\epsilon,+}_{g,\vec{m},\beta} = \M^{\epsilon}_{g,\vec{m}}(Y,\beta)$ forgetting $q_i$ only if $m_i = 0$, and the role of the unit in this chamber is played by the usual fundamental class $\one \in H^*_{\CR}(Y) \cong \H_+$.
\end{remark}

For each $i=1, \ldots, n$, there is an evaluation map
\begin{equation*}
  \ev_i\colon Z^{\epsilon}_{g,n,\beta} \rightarrow \ri X.
\end{equation*}
Here, $\ri X$ is the rigidified inertia stack of $X$, which in this
case is simply
\begin{equation}
\label{eq:riX}
  \ri X := \bigsqcup_{a \in \Z_d} \bigoplus_{i \in F(a)} \O_{\P(\vec d/\gcd(a,d))}(-w_i/\gcd(a,d)).
\end{equation}
(In fact, the evaluation maps land in
$\ri Z= \bigsqcup_{a \in \Z_d} \P(\vec d/\gcd(a,d))$, but given that
our insertions are elements of $\widetilde{\H}:= H^*_{\CR}(X)$, it is
more natural and more consistent with the geometric phase if the
target of evaluation is understood as $\ri X$.)
To define the evaluation maps, let
$\pi\colon \cC \rightarrow \Ze_{g,n,\beta}$ be the universal curve and $\mathcal{L}$ the universal line bundle,
and let
\[\sigma \in \Gamma\left( \bigoplus_{i=1}^M \L^{\otimes w_i} \oplus \bigoplus_{j=1}^N (\L^{\otimes -d} \otimes
\omega_{\pi,\log})\right)\]
be the universal section, whose first $M$ coordinates are zero.  If $\Delta_i \subset \cC$ denotes the stacky divisor corresponding to the $i$th marked point, then
\[\sigma\big|_{\Delta_i} \in \Gamma\left( \bigoplus_{i=1}^M \mathcal{L}^{\otimes w_i} \oplus  \bigoplus_{j=1}^N \L^{\otimes -d} \bigg|_{\Delta_i}\right),\]
using the fact that $\omega_{\pi, \log}|_{\Delta_i}$ is trivial.
Thus, evaluating $\sigma\big|_{\Delta_i}$ at the fiber over a point
$(C;q_1, \ldots, q_n; L; \vec{p})$ in the moduli space yields an
element of $\P(\vec d/\gcd(a,d))$, and by definition, $\ev_i$ sends
$\Ze_{g,\m,\beta}$ to the zero section $\P(\vec d/\gcd(a,d))$ in the
component of \eqref{eq:riX} indexed by $a:= e^{2\pi\i m_i } \in \Z_d$.

\begin{remark}
\label{rem:geomev}
In the geometric chamber, these recover the composition of the usual evaluation map with the inclusion,
\[\M^{\epsilon}_{g,n}(Y,\beta) \rightarrow \ri Y \rightarrow \ri X_+,\]
so in particular, they factor through the inclusion $\ri Y \rightarrow \ri \P(w_1, \ldots, w_M)$.
\end{remark}

\subsection{Virtual cycle and correlators}
\label{subsec:vcycle}

By the results of \cite[Section 5.1]{FJRGLSM}, there exists a virtual cycle
\[[Z_{g,n,\beta}^{\epsilon}]^{\vir} \in A_*(Z_{g,n,\beta}^{\epsilon}).\]
It is constructed using the cosection technique of Kiem--Li \cite{KL}, following the application of the technique by Chang--Li \cite{CL} to the Gromov--Witten theory of the quintic threefold and Chang--Li--Li \cite{CLL} to spin theory.  We outline the basic idea, referring the reader to \cite{FJRGLSM} for details.

The key point is that $Z_{g,n,\beta}^{\epsilon}$ sits inside of the noncompact moduli space $X^{\epsilon}_{g,n,\beta}$ of Landau--Ginzburg quasimaps to $X$, which parameterizes tuples $(C;q_1, \ldots, q_n; L; (\vec{x}, \vec{p}))$ with
\[(\vec{x}, \vec{p}) \in \Gamma\left(\bigoplus_{i=1}^M L^{\otimes w_i} \oplus \bigoplus_{j=1}^N (L^{\otimes -d} \otimes \omega_{\log})\right),\]
satisfying the same conditions as in Definition \ref{def:quasimap}.  By \cite[Proposition 2.5]{CL}, $X^{\epsilon}_{g,n,\beta}$ admits a relative perfect obstruction theory
\begin{equation}
\label{eq:E}
\mathbb{E}^{\bullet} := \left(R\pi_*\left(\bigoplus_{i=1}^M \L^{\otimes w_i} \oplus \bigoplus_{j=1}^N (\L^{\otimes -d} \otimes \omega_{\pi, \log})\right)\right)^{\vee} \rightarrow \mathbb{L}^{\bullet}_{X^{\epsilon}_{g,n,\beta}/D_{g,n,\beta}}
\end{equation}
relative to the moduli space $D_{g,n,\beta}$ parameterizing only
$(C; q_1, \ldots, q_n;L)$, in which
$\pi\colon \cC \rightarrow X^{\epsilon}_{g,n,\beta}$ is the
universal curve and $\L$ the universal line bundle on
$\cC$.
One hopes to define a homomorphism
\begin{equation*}
  \O b_{X^{\epsilon}_{g,n,\beta}/D_{g,n,\beta}} \rightarrow \O_{X^{\epsilon}_{g,n,\beta}}
\end{equation*}
via the derivatives of the superpotential that descends to a cosection
$\sigma\colon \O b_{X^{\epsilon}_{g,n,\beta}} \rightarrow
\O_{X^{\epsilon}_{g,n,\beta}}$ whose fiber is zero exactly over
$\Ze_{g,n,\beta} \subset \Xe_{g,n,\beta}$.
On the components $\Xe_{g,\m,\beta}$ for which each $m_i$ is narrow,
this procedure works, and the cosection technique outputs a virtual
cycle supported on $\Ze_{g,\m,\beta}$.

On the components of $\Xe_{g,n,\beta}$ where not all multiplicities are narrow, on the other hand, an additional step is necessary.  Let $\Xect_{g,n,\beta} \subset \Xe_{g,n,\beta}$ denote the subspace where $x_1(q_k) = \cdots = x_M(q_k) = 0$ for each marked point $q_k$.  (We have $\Xect_{g,\m,\beta} = \Xe_{g,\m,\beta}$ if all multiplicities are narrow, since this condition implies that $L^{\otimes w_i}$ has nonzero multiplicity at $q_k$ for each $i$ and $k$, and hence all sections of $L^{\otimes w_i}$ must vanish at $q_k$.) 

The relative perfect obstruction theory $\mathbb{E}^{\bullet}$ can be modified to yield a relative perfect obstruction theory for $\Xect_{g,n,\beta}$.  Namely, we set
\begin{equation}
\label{eq:pot}
\mathbb{E}^{\bullet}_{\ct}:= \left(R\pi_*\left(\bigoplus_{i=1}^M \left(\L^{\otimes w_i}\otimes\O\left(-\textstyle\sum_{k=1}^n\Delta_k\right)\right) \oplus \bigoplus_{i=1}^N(\L^{\otimes -d} \otimes \omega_{\pi,\log})\right)\right)^{\vee} \rightarrow \mathbb{L}^{\bullet}_{X^{\epsilon, \ct}_{g,n,\beta}/D_{g,n,\beta}},
\end{equation}
where $\Delta_k$ again denotes the stacky divisor in the universal curve
corresponding to the $k$th marked point with its orbifold structure.  Again, using \cite[Proposition
2.5]{CL}, one checks that \eqref{eq:pot} is indeed a relative perfect obstruction theory.
(In the narrow case, the bundles $L^{\otimes w_i}$ and
$L^{\otimes w_i} \otimes \O(-\sum [q_k])$ have the same coarse
underlying bundle and hence the same cohomology, so
$\mathbb{E}^{\bullet}_{\ct} = \mathbb{E}^\bullet$.) 
The derivatives of $W$ now define a homomorphism
\[\O b_{\Xect_{g,n,\beta}/D_{g,n,\beta}} \rightarrow \O_{\Xect_{g,n,\beta}}.\]
Namely, the fiber over a point  $(C;q_1, \ldots, q_n; (\vec{x}, \vec{p})) \in X^{\epsilon, \ct}_{g,n,\beta}$ is the homomorphism
\[\bigoplus_{i=1}^M H^1\left(L^{\otimes w_i} \otimes \O\left(-\sum_{k=1}^n [q_k]\right)\right) \oplus \bigoplus_{j=1}^N H^1(L^{\otimes -d} \otimes \omega_{\log}) \rightarrow \C\]
\[(\dot{x}_1, \ldots, \dot{x}_M, \dot{p}_1, \ldots, \dot{p}_N) \mapsto \sum_{i=1}^M \frac{\d W}{\d x_i}(\vec{x}, \vec{p})\cdot \dot{x}_i + \sum_{j=1}^N \frac{\d W}{\d p_j}(\vec{x}, \vec{p})\cdot \dot{p}_j.\]
A straightforward application of Serre duality checks that this indeed lands in $\C$.  Furthermore, one can check that it descends to a cosection---that is, a homomorphism
\[\sigma\colon \O b_{X^{\epsilon, \ct}_{g,n,\beta}} \rightarrow \O_{X^{\epsilon, \ct}_{g,n,\beta}}\]
out of the absolute obstruction sheaf.  The degeneracy locus of $\sigma$ (the locus of points in $\Xect_{g,n\beta}$ over which the fiber of $\sigma$ is the zero homomorphism) is precisely $Z^{\epsilon}_{g,n,\beta} \subset \Xect_{g,n,\beta}$, and the cosection technique outputs a virtual cycle supported on this locus.

From this discussion, it is straightforward to compute the virtual dimension of each component $Z^{\epsilon}_{g,\m,\beta}$.  In particular, we have
\begin{align}
\label{eq:vdim}
\vdim(Z^{\epsilon}_{g,\m,\beta}) &= \vdim(\Xect_{g,\m,\beta})\\
\nonumber&= \dim(D_{g,\m,\beta}) + (h^0-h^1)\left(\bigoplus_{i=1}^M \left(L^{\otimes w_i}\otimes \O\left(-\sum_{k=1}^n [q_k]\right)\right) \oplus \bigoplus_{j=1}^N (L^{\otimes -d} \otimes \omega_{\log})\right)\\
\nonumber&=\vdim(\M^{\epsilon}_{g,\m}(Z,\beta)) + \sum_{i=1}^M \chi\left(L^{\otimes w_i}\otimes \O\left(-\sum_{k=1}^n [q_k]\right)\right).
\end{align}

\begin{remark}
\label{rem:genuszero}
In genus zero, the definition of the virtual cycle simplifies substantially.  Indeed, the condition that $w_i|d$ implies that $\deg(|L^{\otimes w_i} \otimes \O(-\sum [q_k])|) < 0$ for each $i$ (see \cite[Section 4.2.9]{Clader}).  Thus, the cosection is identically zero, so the cosection-localized virtual class is the usual virtual class of $Z^{\epsilon}_{0,n,\beta} = X^{\epsilon}_{0,n,\beta}$ defined by way of the perfect obstruction theory $\mathbb{E}^{\bullet}$.  Furthermore, $\mathbb{E}^{\bullet}$ is quasi-isomorphic to a vector bundle over a smooth space, so we have
\[[Z^{\epsilon}_{0,n,\beta}]^{\vir} = e\left( \bigoplus_{i=1}^M R^1\pi_*(\L^{\otimes w_i} \otimes \O\left(-\textstyle\sum_{k=1}^n\Delta_k\right))\right) \cap [Z^{\epsilon}_{0,n,\beta}].\]
\end{remark}

Recall that the psi classes are defined by
\[\psi_i = c_1(\mathbb{L}_i) \in A^*(Z^{\epsilon}_{g,n,\beta}), \; \; i=1, \ldots, n,\]
where $\mathbb{L}_i$ is the line bundle whose fiber over a moduli point is the cotangent line to the coarse curve $|C|$ at $q_i$. 

We are now ready to define correlators in the GLSM, following \cite{FJRGLSM}:

\begin{definition}
\label{def:correlators}
Given
\[\phi_1, \ldots, \phi_n \in \widetilde\H\]
and nonnegative integers $a_1, \ldots, a_n$, the associated {\it genus-$g$, degree-$\beta$, $\epsilon$-stable GLSM correlator} is defined as
\begin{equation*}
  \langle \phi_1\psi^{a_1} \cdots \phi_n \psi^{a_n} \rangle^{\epsilon}_{g,n,\beta} := \int_{[\Ze_{g,n,\beta}]^{\vir}} \ev_1^*(\phi_1) \psi_1^{a_1} \cdots \ev_n^*(\phi_n)\psi_n^{a_n}.
\end{equation*}
\end{definition}

\begin{remark}
\label{rem:geomvc}
The cosection construction outlined above can also be used to define a virtual cycle $[\M^{\epsilon}_{g,n}(Y,\beta)]^{\vir}$ in the geometric chamber.  To do so, we set $X^{\epsilon,+}_{g,n,\beta}$ to be the moduli space parameterizing tuples $(\vec{x},\vec{p})$ as above, satisfying representability and the geometric-phase nondegeneracy and stability conditions that
\[\text{ord}_q(\vec{x}) \leq \frac{1}{\epsilon}\]
and $L^{\epsilon} \otimes \omega_{\log}$ is ample.  This admits a relative perfect obstruction theory $\mathbb{E}^{\bullet}_+$ exactly as in \eqref{eq:E}.  Inside $X^{\epsilon,+}_{g,n,\beta}$, we define $X^{\epsilon,\text{ct},+}_{g,n,\beta}$ as the subspace in which $p_1(q_k) = \cdots = p_N(q_k) = 0$ for each marked point $q_k$, whose relative perfect obstruction theory $\mathbb{E}^{\bullet}_{\text{ct},+}$ is obtained from $\mathbb{E}^{\bullet}_+$ by replacing $\omega_{\pi, \log}$ by $\omega_{\pi}$.  The derivatives of $W$ define a homomorphism out of the obstruction sheaf of $\mathbb{E}^{\bullet}_{\text{ct},+}$ that descends to a cosection, yielding a virtual cycle supported on $Z^{\epsilon,+}_{g,n,\beta} = \M^{\epsilon}_{g,n}(Y,\beta)$.

It is a highly nontrivial statement that this construction agrees with the virtual cycle defined by Ciocan-Fontanine--Kim--Maulik in \cite{CFKM}, which was what appeared in our previous work \cite{CJR}.  In fact, in the case where $X_+$ is the quintic threefold and $\epsilon \gg 0$, Chang--Li proved in \cite{CL} that the two definitions of the virtual cycle yield the same correlators up to a sign, and their proof can be adapted to show that the two definitions agree (up to a sign) for any complete intersection in $\P^{N-1}$.  Thus, while we work with the cosection-localized virtual class in what follows, the results of \cite{CJR} for ambient insertions can be deduced as a special case.

We note, furthermore, that the analogue of Remark~\ref{rem:genuszero} in the geometric chamber (known as the orbifold quantum Lefschetz hyperplane principle \cite{Tseng}) requires that $w_i|d_j$ for all $i$ and $j$, since this condition is equivalent to the requirement that the bundle $\bigoplus_{j=1}^N \O_{\P(w_1, \ldots, w_M)}(-d_j)$ be pulled back from the coarse underlying space of $\P(w_1, \ldots, w_M)$.  
\end{remark}

Recall that in replacing the compact-type state space $\H^{\ct}$ by
$\widetilde\H$, it was necessary to choose lifts of elements of
$\H^{\ct}$ under the map
$\eta\colon \widetilde{\H} \rightarrow \H^{\ct}$.
Now that we have defined correlators, we can more precisely address the issue of their independence of the choice of lift.  The statement relies on the following conjecture:

\begin{conjecture}[Broad vanishing]
\label{conj:vanishing}
Let $X^{\text{ct}} \subset X$ denote the zero section,\footnote{Although $X^{\text{ct}}=Z$, we use different notation to clarify the parallel in the geometric phase.} and let $\phi \in H^*_{\CR}(X)$ be such that
\[e(T_{\ri X/\ri X^{\text{ct}}})\phi=0.\]
Then, for any $i \in \{1,\ldots, n\}$, we have
\[\ev_i^*(\phi) \cap [Z^\epsilon_{g, n, \beta}]^\vir = 0.\]
\end{conjecture}

For example, Conjecture~\ref{conj:vanishing} holds in $r$-spin theory, and the analogue in the geometric phase also holds (see Remark~\ref{rem:geomvanishing} below).

\begin{remark}
The name ``broad vanishing" refers to the fact that, if the theory satisfies the further condition that $|F(g)| \geq N$ for any $g \in \Z_d$ such that $F(g) \neq \emptyset$ (c.f. \cite[Condition {\bf (A2)}]{ClaRo}), then Conjecture~\ref{conj:vanishing} implies $\ev_i^*(\phi) \cap [\Ze_{g,n,\beta}]^{\vir} = 0$ whenever $\phi$ is not narrow.  Indeed, under this assumption, $e(T_{\ri X/\ri X^{\text{ct}}})$ vanishes on all noncompact twisted sectors, so any non-narrow element of $H^*_{\CR}(X)$ satisfies $e(T_{\ri X/\ri X^{\text{ct}}})\phi=0$.
\end{remark}

From here, the independence of the correlators of the choice of lift under $\eta$ is nearly immediate:
\begin{lemma}
  \label{lem:vanishing}
  Assume Conjecture~\ref{conj:vanishing}, and let $\phi \in \ker(\eta)$, where $\eta$ is as in \eqref{eq:eta}.
  Then
  \begin{equation*}
    \ev_i^*(\phi) \cap [Z^\epsilon_{g, n, \beta}]^\vir = 0
  \end{equation*}
  for any $i \in \{1, \dotsc, n\}$.
\end{lemma}
\begin{proof}
  By the definition of the Thom isomorphism, the composition of $\eta$
  with the natural restriction map
  $H^*_{\CR}(X, W^{+\infty}) \rightarrow H^*_{\CR}(X)$ is given by
  \begin{equation*}
    \rho\colon H^*_{\CR}(X) \rightarrow H^*_{\CR}(X), \qquad \rho(\phi) = e(T_{\ri X/\ri X^{\text{ct}}}) \phi.
  \end{equation*}
Since $\ker(\eta) \subset \ker(\rho)$, the lemma follows from Conjecture~\ref{conj:vanishing}.
\end{proof}

\begin{remark}
\label{rem:geomvanishing}
In the geometric phase, Conjecture 2.12 does hold, so the proof of Lemma~\ref{lem:vanishing} goes through in that setting with no additional assumptions.  Indeed, by Remark \ref{rem:geomev}, the evaluation maps factor through the inclusion
\[\iota\colon \ri Y \rightarrow \ri \P(w_1, \ldots, w_M),\]
and by Remark~\ref{rem:geomss}, we have $\eta = \iota^*$.  Thus, in the geometric phase one actually has the stronger statement  that $\ev_i^*(\phi) = 0$ whenever $\phi \in \ker(\eta)$.
\end{remark}

\subsection{The $J$-function}
\label{J}

The small $J$-function for $\epsilon$-stable quasimap theory was defined by Ciocan-Fontanine and Kim in \cite{CFKZero} using localization on a graph space, generalizing the original definition in Gromov--Witten theory due to Givental \cite{GiventalEquivariant, GiventalElliptic, GiventalMirror} and the interpretation in terms of contraction maps due to Bertram \cite{Bertram}.  The appropriate modifications for Landau--Ginzburg theory, which we recall below, were carried out by Ross and the third author in \cite{RR}.

Let $\G \Ze_{0,1,\beta}$ denote a ``graph space" version of the moduli space $\Ze_{0,1,\beta}$, which parameterizes the same objects as $\Ze_{0,1,\beta}$ together with the additional datum of a degree-$1$ map $C \rightarrow \P^1$, or equivalently, a parameterization of one component of $C$, and the ampleness condition \eqref{stab} is not required on the parameterized component.  Note that the compatibility condition \eqref{eq:compatibility} implies that the multiplicity at the single marked point must be
\begin{equation}
\label{eq:m1}
m_1 = \left\langle\frac{-\beta-1}{d}\right\rangle,
\end{equation}
in which the symbol $\langle a\rangle$ for a rational number $a$ is defined by the requirement that $0 \leq a <1$ and $\langle a\rangle \equiv a \mod \Z$.

There is an action of $\C^*$ on $\G \Ze_{0,1,\beta}$ given by multiplication on the parameterized component $C_0 \cong \P^1$.  The fixed loci consist of quasimaps for which the marked point and all of the degree $\beta$ lies either over $0$ or over $\infty$ on $C_0$.  We denote by $F^{\epsilon}_{\beta} \subset \G \Ze_{0,1,\beta}$ the fixed locus on which the marked point lies at $\infty$ and all of the degree lies over $0$.  More precisely, when $\beta > 1/\epsilon$, an element of $F^{\epsilon}_{\beta}$ consists of an $\epsilon$-stable Landau--Ginzburg quasimap to $Z$ attached at a single marked point to $C_0$, so
\[F^{\epsilon}_{\beta} \cong \Ze_{0,1,\beta}.\]
When $\beta \leq 1/\epsilon$, on the other hand, such a quasimap would not be stable; instead, $C_0$ is the entire source curve, and the quasimap has a single basepoint of order $\beta$ at $0$.  In either case, there is an evaluation map
\[\ev_{\bullet}\colon F^{\epsilon}_{\beta} \rightarrow \ri X,\]
defined by evaluation at the single marked point $\infty \in C_0$.

The evaluation map we require must be slightly modified from the above, in order to account for the orbifold structure of $X$ (see \cite[Section 3.1]{CCFK}).  Let
\[\mathbf{r}_{\beta}: H^*(\ri X) \rightarrow H^*(\ri X)\]
be multiplication by the number $d_{m_1}$, in which $d_{m_1} := d/\gcd(d \cdot m_1, d)$, and $m_1$ is as defined by \eqref{eq:m1}, and let
\[\iota: \ri X \rightarrow \ri X\]
be the involution coming from inversion on the indexing set $\Z_d$ of the twisted sectors, as well as inverting the banding.  Then we define
\[(\widetilde{\ev_{\bullet}})_* = \iota_* \circ \mathbf{r}_{\beta} \circ (\ev_{\bullet})_* \colon H^*_{\C^*}(F^{\epsilon}_{\beta}) \rightarrow H^*_{\C^*}(\ri X).\]

Analogously to $Z^{\epsilon}_{0,1,\beta}$, the graph space admits a virtual cycle defined by cosection localization.  From here, the $J$-function is defined as follows:

\begin{definition}
Let $z$ denote the equivariant parameter for the action of $\C^*$ on $\G \Ze_{0,1,\beta}$, and let $q$ be a formal Novikov variable.  The {\it small $\epsilon$-stable $J$-function} is
\begin{equation*}
  J^{\epsilon}(q,z) := -z^2\sum_{\beta \geq 0} q^{\beta} (\widetilde{\ev_{\bullet}})_* \left( \frac{[F_{\beta}^{\epsilon}]^{\vir}}{e_{\C^*}(N^{\vir}_{F_{\beta}^{\epsilon}/\G \Ze_{0,1,\beta}})} \right) \in \widetilde\H[[q,z,z^{-1}]].
\end{equation*}
\end{definition}

\begin{remark}
\label{rem:prefactor}
In the geometric phase, this coincides up to the prefactor of $-z^2$ with the image under the pushforward
\[H^*_{\CR}(Y) \rightarrow H^*_{\CR}(X_+)\]
of the small $J$-function of \cite{CFKZero}.  The discrepancy in the prefactor is due to two differences between our set-up and theirs.  First, we use a one-pointed graph space instead of a zero-pointed graph space (a necessary modification, since the marked point in general carries orbifold structure), which changes the localization contributions by a factor of $-z$ to cancel the contribution of automorphisms moving the unmarked point at $\infty$.  Second, our conventions differ by an overall factor of $z$; for example, in the case where $Y$ is semi-positive, the $J$-function of \cite{CFKZero} is of the form $I_0(q) + O(z^{-1})$, whereas ours is of the form $I_0(q) z + O(z^{0})$.
\end{remark}

Using Remark \ref{rem:genuszero} and the discussion above, we can make the $J$-function more explicit.  In particular, when $\beta > 1/\epsilon$, we have
\[[F_{\beta}^{\epsilon}]^{\vir} = [\Ze_{0,1,\beta}]^{\vir}\]
and
\[e_{\C^*}(N^{\vir}_{F_{\beta}^{\epsilon}/\G \Ze_{0,1,\beta}}) = z(z-\psi_1),\]
so
\[(\widetilde{\ev_{\bullet}})_* \left( \frac{[F_{\beta}^{\epsilon}]^{\vir}}{e_{\C^*}(N^{\vir}_{F_{\beta}^{\epsilon}/\G \Ze_{0,1,\beta}})} \right) = \left\langle \frac{\one_{\langle(-\beta-1)/d\rangle}}{z(z-\psi_1)}\right\rangle\one_{\langle(\beta+1)/d\rangle},\]
where the denominator should be understood as a geometric series in $\psi_1$.  The terms of $J^{\epsilon}$ with $\beta \leq 1/\epsilon$ are referred to as {\it unstable terms}, and can be computed explicitly as described in Step 5 of the proof of \cite[Lemma 2.1]{RR}.

Taking $\epsilon \rightarrow 0+$ (that is, requiring the stability condition \eqref{stab} for all $\epsilon >0$),  every term of the $J$-function becomes unstable, so one obtains a generating function that can be computed exactly.  The result is known as the $I$-function $I(q,z) := J^{0+}(q,z)$, and can be calculated explicitly, as was carried out for hypersurfaces in \cite{RR}.  Truncating $I(q,z)$ to powers of $q$ less than or equal to $1/\beta$, more generally, yields an explicit expression for the unstable part of $J^{\epsilon}$ for any $\epsilon$. 

We denote by
\begin{equation*}
  [J^{\epsilon}]_+(q,z) \in \widetilde\H[[q,z]]
\end{equation*}
the part of the $J$-function with non-negative powers of $z$, which has contributions only from the unstable terms, and we let $\mu_{\beta}^{\epsilon}(z)$ denote the $q^{\beta}$-coefficient in $-z\one + [J]^{\epsilon}_+(q,z)$:
\[\sum_{\beta} q^{\beta} \mu_{\beta}^{\epsilon}(z) = -z \one + [J^{\epsilon}]_+(q,z).\]
This series, which is sometimes called the ``mirror transformation", plays a particularly important role in the wall-crossing formula.

We occasionally require a generalization of $J^{\epsilon}(q,z)$ in
which descendent insertions are allowed.
This is defined, for $\t =\t(z) \in \widetilde\H[\![z]\!]$, by
\begin{equation}
\label{bigJ}
J^{\epsilon}(q,\t,z):= -z^2 \sum_{ k \geq 0} \frac{q^{\beta}}{k!}  (\widetilde{\ev_{\bullet}})_* \left( \prod_{i=1}^k\ev_i^*(\t(\psi_i)) \cap \frac{[F_{k,\beta}^{\epsilon}]^{\vir}}{e_{\C^*}(N^{\vir}_{F_{k,\beta}^{\epsilon}/\G \Ze_{0,k+1,\beta}})} \right),
\end{equation}
where $F^{\epsilon}_{k,\beta} \subset \mathcal{G}\Ze_{0,k+1,\beta}$ is defined as the fixed locus in which all but the last marked point and all of the degree are concentrated over $0 \in C_0$, while the last point marked lies at $\infty \in C_0$.  The small $J$-function is recovered by setting $\t = 0$.

\section{Statement of results}
\label{sec:statements}

Ciocan-Fontanine and Kim conjectured wall-crossing formulas for stable quasimap invariants in \cite{CFKZero, CFKHigher}, and they have proven these conjectures in many cases; in \cite{CFKMirror}, they prove a wall-crossing theorem in all genus for any projective complete intersection.  In this section, we review Ciocan-Fontanine--Kim's conjecture in the Landau--Ginzburg context, and we state our main theorem in precise form.

\subsection{Wall-crossing for $J$-functions}

In genus zero, the wall-crossing conjecture states that the function $J^{\epsilon}(q,z)$ lies on the Lagrangian cone defined by the $\infty$-stable hybrid theory---that is, there exists $\t \in \H[[z]]$ such that $J^{\epsilon}(q,z) = J^{\infty}(q,\t,z)$.  The $\t$ in this equation can be determined explicitly from the fact that $J^{\infty}(q,\t,z) = z\one  + \t(-z) + O(z^{-1})$, so we must have
\begin{equation}
\label{eq:JWC}
J^{\epsilon}(q,z) = J^{\infty}(q, z\one + [J^{\epsilon}]_+(q,-z), z) = J^{\infty} \left(q, \sum_{\beta} q^{\beta} \mu^{\epsilon}_{\beta}(-z),z\right).
\end{equation}
This statement was proved in the hypersurface case by Ross and the third author \cite{RR}, and their proof can be adapted to the hybrid setting; see \cite{ClaRo2}.  In the geometric phase, \eqref{eq:JWC} was proved by Ciocan-Fontanine--Kim in \cite{CFKZero}.

\subsection{Wall-crossing for virtual cycles}
\label{subsec:conj}

More generally, Ciocan-Fontanine and Kim lift their conjecture in any genus to the level of virtual cycles.  To state the analogue in Landau--Ginzburg theory, which is the statement of our main theorem in the negative phase, we require some further notation.  The relevant ideas are based on \cite[Section 3.2]{CFKZero}.

Let $\bbeta = (\beta_1, \ldots, \beta_k)$ be a tuple of nonnegative
integers.
Let $m_i := \left\langle \frac{\beta_i+1}{d} \right\rangle$ for all
$i$, and let $\m = (m_1, \ldots, m_k)$.
Then there is a morphism
\begin{equation*}
  \label{eq:beta}
  b_{\bbeta}\colon \Ze_{g,n+\m, \beta - \sum_{i=1}^k \beta_i}  \rightarrow \Ze_{g,n,\beta},
\end{equation*}
defined as follows.  For $(C; q_1, \ldots, q_{n+k}; L;\vec{p}) \in \Ze_{g,n+\m, \beta - \sum_{i=1}^k \beta_i}$, let $\widetilde{C}$ be the partial coarsening of $C$ obtained by forgetting the last $k$ marked points and their orbifold structure.  The bundle
\[L\otimes \O\left( \sum_{i=1}^k \frac{-\beta_i-1}{d}[\overline{q}_{n+i}]\right)\]
has multiplicity zero at each of the last $k$ marked points, and hence it is pulled back from a bundle $\widetilde{L}$ on $\widetilde{C}$.  For $j=1, \ldots, N$, let
\begin{align*}
  \widetilde{p}_j \in &\Gamma\left(\widetilde{L}^{\otimes -d} \otimes \omega\left(\sum_{i=1}^n [q_i]\right)\right)\\
  =  &\Gamma\left(\left(L^{\otimes -d} \otimes \omega\left(\sum_{i=1}^{n+k} [q_i]\right)\right) \otimes \O\left(\beta_i[q_{n+i}] \right)  \right)
\end{align*}
be the bundle on $\widetilde{C}$ obtained from $p_j \in \Gamma(L^{\otimes -d} \otimes \omega_{\log})$ by the natural map.  Then, setting $\widetilde{\vec{p}} = (\widetilde{p}_1, \ldots, \widetilde{p}_N)$, we define
\[b_{\bbeta}(C;q_1, \ldots, q_{n+k}; L;\vec{p}) := (\widetilde{C}; q_1, \ldots, q_{n};\widetilde{L};\widetilde{\vec{p}}),\]
assuming the latter is stable.  It is possible, however, that $\widetilde{C}$ contains rational tails---that is, genus-zero components with a single special point---on which
\begin{equation}
\label{eq:deg}
\deg\left(\widetilde{L}^{\otimes -d} \otimes \omega\left(\sum_{i=1}^n [q_i]\right)\right) \leq \frac{1}{\epsilon},
\end{equation}
and which thus violate the ampleness condition \eqref{stab}.  To define $b_{\bbeta}$, we contract such a component and replace it with a basepoint of order equal to the left-hand side of \eqref{eq:deg} at the point where the component was attached, as formalized below.  This may create a new rational tail, so we repeat the process inductively until stability is achieved.

The definition of the morphism
\begin{equation*}
  c\colon Z^{\infty}_{g,n,\beta} \rightarrow Z^{\epsilon}_{g,n,\beta},
\end{equation*}
which contracts unstable rational tails and replaces them with
basepoints, is similar.
Namely, suppose that
$(C;q_1, \ldots, q_n; L;\vec{p}) \in Z^{\infty}_{g,n,\beta}$, where
$C$ has a rational tail $C_0$ and
$\beta_0 := \deg(L^{\otimes -d} \otimes \omega_{\log}|_{C_0}) \leq
1/\epsilon$.
Let $\widetilde{C}$ be obtained from $\overline{C \setminus C_0}$ by
forgetting the orbifold structure at the point $q$ where $C_0$ meets
the rest of $C$.
Then the line bundle
\[L\big|_{\overline{C \setminus C_0}} \otimes \O\left(\frac{-\beta_0-1}{d}[\overline{q}]\right)\]
has multiplicity zero at $q$, so as above, it is pulled back from a bundle $\widetilde{L}$ on $\widetilde{C}$.  For $j=1, \ldots, N$, we let
\begin{align*}
  \widetilde{p}_j \in &\Gamma\left(\widetilde{L}^{\otimes -d} \otimes \omega\left(\sum_{i=1}^n [q_i]\right)\right)\\
  =  &\Gamma\left(\left(L^{\otimes -d} \otimes \omega\left(\sum_{i=1}^{n} [q_i] + [q]\right) \right) \otimes \O\left(\sum_{i=1}^{k} \beta_0[q] \right)  \right)
\end{align*}
be obtained from the restriction of $p_j$ to $\overline{C \setminus C_0}$.  We then set $\widetilde{\vec{p}} = (\widetilde{p}_1, \ldots, \widetilde{p}_N)$ and let
 \[c(C;q_1, \ldots, q_{n+k}; L;\vec{p}) := (\widetilde{C}; q_1, \ldots, q_{n};\widetilde{L};\widetilde{\vec{p}}),\]
assuming the latter is stable.  If not, we iterate the procedure until we reach an $\epsilon$-stable Landau--Ginzburg quasimap.

By a slight abuse of notation, we extend $\ev_i^*$ to $\widetilde\H[\psi_i]$ by linearity in $\psi_i$, and similarly, we allow $b_{\bbeta *}$ and $c_*$ to operate linearly in $q$.  Equipped with these definitions, we can give a precise statement of Theorem~\ref{thm:1} in the Landau--Ginzburg phase:

\begin{theorem}
	\label{thm:main}
	Fix $g \geq 0$, $n \geq 1$, and a tuple of insertions
        $\alpha_1, \ldots, \alpha_n \in \widetilde\H$.
        Then
\begin{align*}
&\sum_{\beta} q^\beta  \left(\prod_{a=1}^n \ev_a^*(\alpha_a) \cap [\Ze_{g,n,\beta}]^{\vir}\right) =\\
&\sum_{\beta_0, \beta_1, \ldots, \beta_k} \frac{q^{\beta_0}}{k!} b_{\bbeta*}c_* \left(\prod_{a=1}^n \ev_a^*(\alpha_a)\prod_{i = 1}^k q^{\beta_i} \ev_{n + i}^*(\mu_{\beta_i}^\epsilon(-\psi_{n + i})) \cap [Z^{\infty}_{g,n+k,\beta_0}]^{\vir} \right).
\end{align*}
	(Note that, by the definition of $\mu^{\epsilon}_{\beta}(z)$, the expression inside $c_*(\cdots)$ is supported on the substack of $Z^{\infty}_{g,n+k,\beta_0}$ on which the multiplicity at the marked point $q_{n+i}$ is $\langle \frac{\beta_i+1}{d}\rangle$, so the morphism $b_{\vec{\beta}}$ is well-defined.)
\end{theorem}

\begin{remark}
Note that, in contrast to the statement of Theorem~\ref{thm:1}, we now assume that $n \geq 1$.  The $n=0$ case is handled separately in the appendix.
\end{remark}

\begin{remark}
The analogue of Theorem~\ref{thm:main} in the geometric phase is Theorem 2.6 of \cite{CJR}, except that in that case (1) the insertions are not required to be ambient, and (2) the weights $w_i$ are all equal to one.  Thus, assuming the equivalence of Ciocan-Fontanine--Kim--Maulik's virtual cycle and the cosection-localized virtual cycle (Remark~\ref{rem:geomvc}), the proof discussed below reproduces, in a fundamentally similar but more complicated way, the ambient case of \cite{CJR}.
\end{remark}

\subsection{Twisted theory}
\label{subsec:twisted}

In order to prove Theorem~\ref{thm:main}, we must introduce a twist of the virtual class by a certain equivariant Euler class.

Denoting the universal curve over $\Ze_{g,n,\beta}$ again by
$\pi\colon \cC \rightarrow \Ze_{g,n,\beta}$ and the universal line
bundle on $\cC$ by $\L$, we let $\C^*$ act trivially on
$\Ze_{g,n,\beta}$.
We denote the equivariant parameter by $\lambda$, and to distinguish
this action from the one on the graph space considered above, we write
$\C^*$ as $\C^*_{\lambda}$.
Furthermore, we denote by $\C_{(\lambda)}$ a nonequivariantly trivial
line bundle with a $\C^*$-action of weight $1$.

Define
\begin{equation*}
  [\Ze_{g,n,\beta}]^{\vir}_\tw := \frac{[\Ze_{g,n,\beta}]^{\vir}}{e_{\C^*_{\lambda}}\left(R\pi_*\left(\mathcal{P^\vee} \otimes \C_{(\lambda)}\right)\right)},
\end{equation*}
where
\[\cP:= \L^{\otimes -d} \otimes \omega_{\pi,\log}.\]
Using this, we can define a twisted $J$-function by
\begin{equation*}
  J^{\epsilon}_\tw(q,z) := -z^2 e_{\C^*_{\lambda}}(\O_{Z}(-1) \otimes \C_{(\lambda)}) \sum_{\beta \geq 0} q^{\beta}  (\widetilde{\ev_{\bullet}})_* \left( \frac{[\Ze_{0,1,\beta}]^{\vir} \cap e_{\C^*_{\lambda} \times \C^*_z}(-R\pi_*\left(\cP^\vee \otimes \C_{(\lambda)})\right)}{ e_{\C^*_z}(N_{F^{\epsilon}_{\beta}})} \right),
\end{equation*}
where $\C^*_z$ denotes the $\C^*$-action on the graph space.  From here, we define a twisted mirror transformation by
\[\sum_{\beta} q^{\beta} \mu_{\beta}^{\epsilon, \tw}(z) = -\one z + [J^{\epsilon}_{\tw}(z)]_+.\]

The twisted version of Theorem~\ref{thm:main} is the following:

\begin{theorem}
	\label{thm:twistedWC}
        Fix $g \geq 0$, $n \geq 1$, and a tuple of insertions
        $\alpha_1, \ldots, \alpha_n \in \widetilde\H$.
        Then
\begin{align*}
&\sum_{\beta} q^\beta  \left(\prod_{a=1}^n \ev_a^*(\alpha_a) \cap [\Ze_{g,n,\beta}]^{\vir}_{\tw}\right) =\\
&\sum_{\beta_0, \beta_1, \ldots, \beta_k} \frac{q^{\beta_0}}{k!} b_{\bbeta*}c_* \left(\prod_{a=1}^n \ev_a^*(\alpha_a)\prod_{i = 1}^k q^{\beta_i} \ev_{n + i}^*(\mu_{\beta_i}^{\epsilon,\tw}(-\psi_{n + i})) \cap [Z^{\infty}_{g,n+k,\beta_0}]^{\vir}_{\tw} \right).
\end{align*}
\end{theorem}

In fact, the twisted wall-crossing theorem implies the untwisted one:

\begin{lemma}
\label{lem:twuntw}
Theorem~\ref{thm:twistedWC} implies Theorem~\ref{thm:main}.
\begin{proof}
The proof is identical to Lemma 2.10 of \cite{CJR}: one recovers the untwisted theorem from the twisted one by taking the top power of $\lambda$ on both sides.
\end{proof}
\end{lemma}

\begin{remark}
In the geometric phase, the definitions of $[\M^{\epsilon}_{g,n}(Y,\beta)]^{\vir}_{\tw}$, $J^{\epsilon}_{\tw}$, and $\mu^{\epsilon,\tw}_{\beta}$ are exactly parallel to the above, except that $\mathcal{P}$ is replaced by $\mathcal{L}$; see \cite[Section 2.4]{CJR}.
\end{remark}

\subsection{Wall-crossing for dual graphs}
\label{subsec:graphWC}

In contrast to the quasimap wall-crossing of \cite{CJR}, the proof of Theorem~\ref{thm:twistedWC} requires an induction on topological types, and for this, we must introduce one further refinement of the wall-crossing statement.

Let $\Gamma$ be a prestable dual graph---that is, a set of vertices, edges, and numbered legs---for which each vertex $v$ is decorated with a genus $g(v) \geq 0$ and a degree $\beta(v) \geq 0$, and each half-edge $h$ (including the legs) is decorated with a multiplicity $m(h) \in \{0, \frac{1}{d}, \ldots, \frac{d-1}{d}\}$.  We assume that for each edge $e$, the multiplicities at the two half-edges $h$ and $h'$ satisfy
\begin{equation}
\label{eq:comp1}
m(h) + m(h') \in \Z
\end{equation}
and for each vertex $v$, the multiplicities at the set $H(v)$ of incident half-edges satisfy
\begin{equation}
\label{eq:comp2}
\frac{-\beta(v) + 2g(v) -2 + |H(v)|}{d} - \sum_{h \in H(v)} m(h) \in \Z.
\end{equation}
Denote by $\vec{m}(v)$ the tuple $\{m(h)\}_{h \in H(v)}$.

In addition, we equip $\Gamma$ with a number of further decorations.  Namely, let
\[v_{\bullet} \in V(\Gamma)\]
be any be any vertex such that $\beta(v_{\bullet}) > 0$, and let
\[\n\colon V(\Gamma) \rightarrow \mathbb{N}\]
be any function with $\n(v_{\bullet}) = 0$.
We denote by $\Gamma_{\n}$ the graph obtained from $\Gamma$ by adding
$\n(v)$ additional legs to each vertex $v$ and assigning multiplicity
$\frac{1}{d}$ to each of them.

Similarly to the recursive structure of the boundary of the moduli space of curves,
for any stable $(\Gamma, v_{\bullet}, \n)$ and any stability parameter
$\epsilon$, we construct an explicit fiber product $\Ze_{\Gamma}$ and a morphism
\begin{equation*}
  \iota_\Gamma\colon \Ze_{\Gamma} \to \Ze_{g,n+n',\beta}
\end{equation*}
that is a finite cover of the closure
of the locus of $\epsilon$-stable Landau--Ginzburg quasimaps with
decorated dual graph $\Gamma_{\n}$.  Here, denoting by $V(\Gamma)$, $E(\Gamma)$, and $L(\Gamma)$ the vertex, edge and leg sets, respectively, we have
\begin{equation}
  \label{eq:dd}
  \begin{aligned}
    g &= h^1(\Gamma) + \sum_{v \in V(\Gamma)} g(v),\\
\beta &= \sum_{v \in V(\Gamma)} \beta(v),\\
n &= |L(\Gamma)|,\\
n' &= \sum_{v \in V(\Gamma)} \n(v).
  \end{aligned}
\end{equation}
Explicitly, for each vertex $v \in V(\Gamma)$, denote by $\m'(v)$ the tuple of multiplicities at half-edges incident to $v$ together with $\n(v)$ additional multiplicites of $\frac{1}{d}$.  Then $\Ze_\Gamma$ is the fiber
product of the moduli spaces $\Ze_{g(v), \m'(v), \beta(v)}$ for
each $v \in V(\Gamma)$, glued over $(\ri Z)^{|E(V)|}$
via the evaluation maps at the half-edges.  The degree of $\iota_\Gamma$ onto its image is given by
\begin{equation}
\label{eq:factor}
  \frac{|\operatorname{Aut}(\Gamma_{\n})|}{\prod_{h} d_{m(h)}},
\end{equation}
where the product is over all half-edges.  Note that $\Ze_{\Gamma}$ depends on the choice of
$\n$, although we suppress it from the notation.

There is also another description of $\Ze_{\Gamma}$ useful for the
definition its virtual class.
Let $D_\Gamma$ be the fiber product over $B\C^*$ of the stacks
$D_{g(v), n(v)+n'(v), \beta(v)}$, which is defined analogously to
$\Ze_\Gamma$.
This stack comes with a local complete intersection morphism
$j_\Gamma\colon D_\Gamma \to D_{g, n + n', \beta}$ (see
\cite[Lemma~4.4]{CLL}).
Using the cartesian diagram
\begin{equation}
\label{eq:vpullback}
  \begin{tikzpicture}
    \matrix (m) [matrix of math nodes,row sep=2em,column sep=4em,minimum width=2em] {
    \Ze_\Gamma & \Ze_{g, n + n', \beta} \\
    D_\Gamma & D_{g, n + n', \beta}, \\};
  \path[-stealth]
    (m-1-1) edge node[above] {$\iota_\Gamma$} (m-1-2) edge (m-2-1)
    (m-1-2) edge (m-2-2)
    (m-2-1) edge node[above] {$j_\Gamma$} (m-2-2);
  \end{tikzpicture}
\end{equation}
we can easily define a virtual cycle for $\Ze_{\Gamma}$:
namely, defining $\Xect_{\Gamma}$ in analogy to
Section~\ref{subsec:vcycle}, we set
$[\Ze_{\Gamma}]^\vir := [\Xect_{\Gamma}]^\vir_\sigma$ by
virtual pullback along $\iota_\Gamma$.
Alternatively, we can define a relative perfect obstruction theory
similarly to \eqref{eq:pot}, except that we now utilize the
universal curve over $\Xect_{\Gamma}$ and the base of
the obstruction theory is $D_\Gamma$.

Similarly to Section~\ref{subsec:twisted}, there is also a twisted virtual cycle for $Z^{\epsilon}_{\Gamma}$, but this requires the additional choice of a coloring of the vertices
\[V(\Gamma) = V_0 \cup V_{\infty}\]
such that $v_{\bullet} \in V_0$.  For any such coloring, we define
\begin{equation}
\label{eq:Gammatw}
  [Z^{\epsilon}_{\Gamma}]^{\vir}_\tw = \frac{[Z^{\epsilon}_{\Gamma}]^{\vir}}{\prod_{v \in V_0} e_{\C^*_{\lambda}}(R\pi_{v,*}(\mathcal{P^\vee}) \otimes \C_{(\lambda)}) \prod_{v \in V_\infty} e_{\C^*_{\lambda}}(R\pi_{v,*}(\cP) \otimes \C_{(-\lambda)})},
\end{equation}
in which $\pi_v$ denotes the projection from the component of the universal curve corresponding to the vertex $v$.

Our goal, now, is to state a wall-crossing theorem for each choice of
decorated dual graph.
We require a bit of notation. 
First, for each $l \geq 0$, let $\Gamma + l$ be the graph obtained from
$\Gamma$ by adding $l$ additional legs of multiplicity $\frac{1}{d}$
to the vertex $v_{\bullet}$, and let
\begin{equation*}
  \pi_l\colon Z^{\epsilon}_{\Gamma + l} \rightarrow Z^{\epsilon}_{\Gamma}
\end{equation*}
be the morphism that forgets these additional legs.

The insertions at these $l$ additional legs in the wall-crossing
theorem are expressed in terms of a universal series
\begin{equation*}
  \epsilon^{\lambda_0}(z) \in R[[z]],
\end{equation*}
where the ground ring is $R := \widetilde\H(\lambda)[[y]]$ for a
formal parameter $y$ that is given geometric meaning in what follows.  The precise definition of $\epsilon^{\lambda_0}(z)$ is given in
Sections~\ref{subsec:eqP1} and \ref{sec:proof} below, in terms of the contribution to a
virtual localization from unmarked trees of rational curves; for now,
we require just one property:
\begin{equation}
  \label{eq:Ty}
  \epsilon^{\lambda_0}(z) \in y \cdot R[[z]],
\end{equation}
which is immediate from the definition in Section \ref{subsec:eqP1}.

For each $i \in \{1, \ldots, l\}$, let
\begin{equation*}
  \ev_{v_{\bullet},i}\colon Z^{\epsilon}_{\Gamma + l} \rightarrow X_{(1/d)} \subset \ri X
\end{equation*}
be the evaluation map at the $i$th of the $l$ additional markings on
$v_{\bullet}$, which by construction maps to the
multiplicity-$\frac{1}{d}$ sector $X_{(1/d)}$ in $\ri X$.

\begin{theorem}
  \label{thm:graphWC}
  Fix $g \geq 0$, $\beta \geq 0$, and $n \geq 1$, and let $\Gamma$ be
  a prestable dual graph with these discrete data.
  Fix decorations
  $v_{\bullet}$ and
  $\n\colon V(\Gamma) \rightarrow \mathbb{Z}^{\geq 0}$ such that
  $(\Gamma, v_{\bullet}, \n)$ is stable.
  Then, for any coloring $V(\Gamma) = V_0 \cup V_{\infty}$ as above and any tuple of insertions
  $\alpha_1, \ldots, \alpha_n \in \widetilde\H$,
  \begin{multline}
    \sum_{l=0}^{\infty} \frac{1}{l!} \pi_{l*} \left( \prod_{a \in L'} \ev_{v_{\bullet},a}^*(\epsilon^{\lambda_0}(\psi)) \prod_{\substack{ v \in V(\Gamma)\\ a \in L(v)}} \ev_{v,a}^*(\alpha_a)\cap [Z^{\epsilon}_{\Gamma + l}]^{\vir}_{\tw}\right)=\sum_{l=0}^{\infty} \frac{1}{l!} \pi_{l*} \left( \sum_{K = \{k_v\}_{v \in V(\Gamma)}} \frac{1}{\prod_v k_v!}\right.\\
    \left.\sum_{\substack{\vec{\beta} = \{\beta_i^{v}\}_{v \in V(\Gamma), i \in [k_v]}\\ \beta_0^v + \beta_1^v + \cdots + \beta_{k_v}^v = \beta(v)}} \hspace{-0.75cm} b_{\vec{\beta} *} c_{*} \left(\prod_{a \in L'} \ev_{v_{\bullet},a}^*(\epsilon^{\lambda_0}(\psi) ) \prod_{\substack{ v \in V(\Gamma)\\ a \in L(v)}} \ev_{v,a}^*(\alpha_a) \prod_{\substack{ v \in V(\Gamma)\\ j \in [k_v]}} \ev_{v,j}^*(\mu_{\beta_j^v}^{\epsilon, \tw}(-\psi)) \cap [Z^{\infty}_{\Gamma_{K,\bbeta} + l}]^{\vir}_\tw\right) \right),
  \end{multline}
  where $L(v)$ denotes the leg set of $v$ (excluding the extra legs specified by $\n$), $L'$
  denotes the set of extra legs, and $\ev_{v,a}$ denote the evaluation
  maps at the corresponding marked points.  The morphisms $b_{\vec{\beta}}$ and $c$
  are defined analogously to Section~\ref{subsec:conj}, and
  $\Gamma_{K, \bbeta}$ is the graph obtained from $\Gamma$ by
  modifying each vertex $v$ to have degree $\beta_0^v$ and $k_v$
  additional legs.
\end{theorem}

\begin{remark}
Although we leave the multiplicities at the $k_v$ additional legs of each vertex unspecified, the expression on the right-hand side is in fact supported on the locus where the multiplicities are $\left\langle \frac{\beta_1^v+1}{d}\right\rangle, \ldots, \left\langle \frac{\beta_{k_v}^v+1}{d}\right\rangle$.
\end{remark}

This statement implies the twisted, and hence the untwisted, wall-crossing of the previous subsection:

\begin{lemma}
\label{lem:graphWC}
Theorem~\ref{thm:graphWC} implies Theorem~\ref{thm:twistedWC}.
\begin{proof}
  Take $\Gamma$ to be a graph with a single vertex $v$, and the
  additional decorations to be $v = v_{\bullet} \in V_0$ and
  $\n(v) = 0$.
  Then $Z^{\epsilon}_{\Gamma} = \Ze_{g,n,\beta}$ and
  $Z^{\infty}_{\Gamma_{K,\bbeta}} = Z^{\infty}_{g,n+k,\beta_0}$.
  The resulting equality depends on the parameter $y$, but setting
  $y=0$ and applying \eqref{eq:Ty}, we see that only the $l=0$ term
  contributes.
  This term is precisely the statement of Theorem~\ref{thm:twistedWC}.
\end{proof}
\end{lemma}

\begin{remark}
The definition of $Z^{\epsilon}_{\Gamma}$ and its virtual cycle via \eqref{eq:vpullback} work equally well in the geometric phase, after replacing the condition \eqref{eq:comp2} by
\[\beta - \sum_{h \in H(v)} m(h) \in \Z.\]
Replacing $\mathcal{P}$ by $\mathcal{L}$ in \eqref{eq:Gammatw} gives the definition of the twisted virtual cycle.  The additional legs in $\Gamma + l$ should have multiplicity zero in the geometric phase (see Remark~\ref{rem:geommults}), but after this modification, Theorem~\ref{thm:graphWC} and Lemma~\ref{lem:graphWC} generalize immediately.
\end{remark}

\section{Twisted graph space}
\label{sec:twistedgs}

The proof of Theorem~\ref{thm:graphWC} is by $\C^*$-localization on a ``twisted graph space".  This space, which is closely related to the space of mixed-spin $p$-fields considered by Chang--Li--Li--Liu \cite{CLLL, CLLL2}, was introduced in our previous work \cite{CJR} in the geometric phase.  We adapt the definition to the Landau--Ginzburg phase in this section and derive the properties we require for the proof of Theorem~\ref{thm:graphWC}.

\subsection{Definition of the twisted graph space}
\label{subsec:tgs}

In the language of the GLSM, the twisted graph space $\P X^{\epsilon}_{g,n,\beta,\delta}$ consists of Landau--Ginzburg quasimaps to the GIT quotient
\[\P X:=((\C^M \times \C^N) \times \C^2)\GIT_{\widetilde{\theta}} (\C^*)^2,\]
in which the action of $(\C^*)^2$ is
\begin{equation*}
  (g,t) \cdot (\vec{x}, \vec{p} , z_1, z_2) = (g^{w_1}x_1, \ldots, g^{w_M}x_M, g^{-d}p_1, \ldots, g^{-d}p_N, g^dtz_1, tz_2)
\end{equation*}
and the rational character
$\widetilde{\theta}\colon (\C^*)^2 \rightarrow \C^*$
(c.f. Remark~\ref{rem:rational}) is
\[\widetilde{\theta}(g,t)=\theta(g)^{d\epsilon} t^3\]
for the negative character $\theta$ used to define $X$. The superpotential is extended trivially to the new factors, and the $R$-charge acts with weight $1$ on the $p$-coordinates, weight $-1$ on $z_1$, and weight $0$ on all other coordinates.

More concretely,
\[\P X^{\epsilon}_{g,n,\beta,\delta} = \{(C;q_1, \ldots, q_n; L_1, L_2; x_1, \ldots, x_M, p_1, \ldots, p_N, z_1, z_2)\},\]
where $(C; q_1, \ldots, q_n)$ is an $n$-pointed prestable orbifold curve of genus $g$ with nontrivial isotropy only at special points, $L_1$ and $L_2$ are orbifold line bundles with $P_1:=L_1^{\otimes -d} \otimes \omega_{\log}$,
\begin{equation}
\label{eq:TGSdegrees}
\beta = \deg(P_1), \; \; \; \delta = \deg(L_2),
\end{equation}
and
\[(\vec{x}, \vec{p}, \vec{z}) := (x_1, \ldots, x_M, p_1, \ldots, p_N, z_1, z_2) \in \Gamma\left(\bigoplus_{i=1}^M L_1^{\otimes w_i} \oplus \bigoplus_{j=1}^N P_1 \oplus (P_1^{\vee} \otimes L_2) \oplus L_2\right).\]
We require this data to satisfy the following conditions:
\begin{itemize}
\item {\it Representability}:  For every $q \in C$ with isotropy group $G_q$, the homomorphism $G_q \rightarrow \C^*$ giving the action of the isotropy group on the bundle $L_1$ is injective.
\item {\it Nondegeneracy}:  The sections $z_1$ and $z_2$ never simultaneously vanish.  Furthermore, for each point $q$ of $C$ at which $z_2(q) \neq 0$, we have
\begin{equation}
\label{eq:stab1}
\text{ord}_q(\vec{p}) \leq \frac{1}{\epsilon},
\end{equation}
and for each point $q$ of $C$ at which $z_2(q) = 0$, we have
\begin{equation}
\label{eq:stab2}
\text{ord}_q(\vec{p}) = 0.
\end{equation}
(We note that this can be phrased, as in \cite[Section 2.1]{CFKBigI}, as a length condition bounding the order of contact of $(\vec{x}, \vec{p}, \vec{z})$ with the unstable locus of $\P X$.)
\item {\it Stability}: The $\Q$-line bundle
\begin{equation}
\label{eq:TGSstability}
(P_1^{\epsilon} \otimes L_2^{\otimes 3}) \otimes \omega_{\log}
\end{equation}
is ample.
\end{itemize}

We decompose $\P\Xe_{g,n,\beta,\delta}$ by multiplicities, denoting by
\[\P\Xe_{g,\m,\beta,\delta} \subset \P\Xe_{g,n,\beta,\delta}\]
the open and closed substack on which the multiplicity of $L_1$ at the marked point $q_i$ is $m_i$.   (The multiplicity of $L_2$ is always zero, since otherwise the sections $z_1$ and $z_2$ would both vanish.) 

Note that the nondegeneracy condition implies that $z_1$ and
$z_2$ define a map
\[f\colon C \to \P(P_1^\vee \oplus \O_C).\]
Analogously to what was stated in Section~\ref{subsec:vcycle}, we let
\[\P \Xect_{g,n,\beta,\delta} \subset \P\Xe_{g,n,\beta,\delta}\]
be the locus on which $x_1(q_k) = \cdots = x_M(q_k) = 0$ for each
marked point $q_k$.
Then there is a relative perfect obstruction theory\footnote{This is
  not exactly the standard definition of the perfect obstruction
  theory of a GLSM moduli space, which in this case would be relative
  to a moduli space parameterizing $(C;q_1, \ldots, q_n; L_1, L_2)$.
  However, by the Euler sequence, the resulting (cosection-localized)
  virtual cycles agree.}
\begin{equation*}
  \widetilde{\mathbb{E}}^{\bullet} := \left(R\pi_*\left(\bigoplus_{i=1}^M \left(\L_1^{\otimes w_i}\left(-\textstyle \sum_{k=1}^n\Delta_k\right)\right) \oplus \bigoplus_{j=1}^N\cP_1\right) \oplus R\pi_* f^* T_{\P(\cP_1^\vee \oplus \O)/\cC}\right)^\vee \to \mathbb L_{\P \Xect_{g,n,\beta,\delta}/D_{g,n,\beta}}^\bullet.
\end{equation*}

The derivatives of the superpotential once again define a homomorphism
out of the relative obstruction sheaf, which descends to a cosection
\begin{equation*}
  \widetilde{\sigma}\colon \O b_{\P \Xect_{g,n,\beta,\delta}} \rightarrow \O_{\P \Xect_{g,n,\beta,\delta}}.
\end{equation*}
The degeneracy locus of $\widetilde{\sigma}$ is the substack
\begin{equation*}
  \P \Ze_{g,n,\beta,\delta} \subset \P \Xect_{g,n,\beta,\delta}
\end{equation*}
parameterizing Landau--Ginzburg quasimaps that land in the critical
locus $\P Z \subset \P X$ of the extended superpotential
$\widetilde{W}:= \sum_{j=1}^N p_j F_j(\vec{x})$.
Explicitly, $\P Z$ is a $\P^1$-bundle over $Z= \P(\vec{d})$ and $\P \Ze_{g,\m,\beta,\delta}$ is the locus in
$\P \Xect_{g,n,\beta,\delta}$ on which the sections
$x_1, \ldots, x_M$ are all identically zero.  This locus is proper, because the trivial lift of $\widetilde{\theta}$ is again a good lift.  The cosection construction thus yields a virtual cycle
\begin{equation*}
  [\P\Ze_{g,n,\beta,\delta}]^{\vir} \in A_*(\P\Ze_{g,n,\beta,\delta}),
\end{equation*}
and a similar computation to \eqref{eq:vdim} shows that
\[\vdim(\P\Ze_{g,n,\beta,\delta}) = \vdim(\Ze_{g,n,\beta}) + 2\delta+ \beta+1-g.\]

To define evaluation maps, let
\[\varsigma \in \Gamma(\cP_1^{\oplus N} \oplus (\cP_1^{\vee} \otimes \L_2) \oplus \L_2)\]
denote the universal section over the universal curve
$\pi\colon \cC \rightarrow \P \Ze_{g,n,\beta,\delta}$. 
If $\Delta_i \subset \cC$ denotes the divisor corresponding to
the $i$th marked point, then
\begin{equation*}
  \varsigma|_{\Delta_i} \in \Gamma\left(\left. \bigoplus_{j=1}^N \left(\L_1^{\otimes -d} \otimes \O\left(-\textstyle\sum_{k=1}^n\Delta_k\right)\right) \oplus (\L_1^{\otimes d} \otimes \L_2) \oplus \L_2\right|_{\Delta_i}\right),
\end{equation*}
since $\omega_{\pi, \log}|_{\Delta_i}$ is trivial.
Thus, restricting $\varsigma|_{\Delta_i}$ to a fiber gives an element
of $\ri \P Z \subset \ri \P X$, and we obtain an evaluation map
\begin{equation*}
  \ev_i\colon \P\Ze_{g,n,\beta,\delta} \rightarrow \ri \P X
\end{equation*}
that sends $\P\Ze_{g,\m,\beta,\delta}$ to the twisted sector indexed by
$e^{2\pi \i m_i} \in \Z_d$.  Analogously to the definition of the compact-type state space in Section~\ref{sec:ss}, the insertions for the twisted graph space are drawn from $H^*_{\CR}(\P X)$, which, as a vector space, is isomorphic to $\widetilde{\H} \otimes H^*(\P^1)$.

As in Section~\ref{subsec:graphWC}, one can define substacks of
$\P \Ze_{g,n,\beta,\delta}$ and finite covers thereof associated to
decorated dual graphs.
Namely, let $(\Gamma, v_{\bullet}, \n)$ be as in
Section~\ref{subsec:graphWC}, and define $\widetilde{\Gamma}$ by adding an
additional decoration
\begin{equation*}
  D\colon V(\Gamma) \rightarrow \mathbb{N}
\end{equation*}
to each vertex of $\Gamma$.  We say that $(\widetilde\Gamma, v_{\bullet}, \n)$ is {\it stable} if every vertex $v \in V(\widetilde\Gamma)$ for which $D(v) = 0$ becomes stable after adding $\n(v)$ additional legs to each vertex $v$.  In particular, each element of $\P \Ze_{g,n,\beta,\delta}$ defines a stable decorated dual graph, where $D(v)$ is the degree of the restriction of $L_2$ to the component indexed by $v$.

For any stable decorated graph $(\widetilde\Gamma, v_{\bullet}, \n)$
as above, and any stability parameter $\epsilon$, we denote by
$\P \Ze_{\widetilde\Gamma}$ the fiber product of the moduli spaces
$\P Z^{\epsilon}_{g(v), \m'(v), \beta(v), D(v)}$ over
$\prod_{e \in E(\widetilde\Gamma)} \ri \P Z$.  (Here, as above, $\m'(v)$ records the multiplicities at half-edges incident to $v$ as well as $\n(v)$ additional multiplicities of $\frac{1}{d}$.)  This fiber product imposes the agreement, for each edge
$e \in E(\widetilde\Gamma)$ with half-edges $h_1$ and $h_2$ incident
to vertices $v_1$ and $v_2$, of the evaluation maps corresponding to
$h_1$ and $h_2$.
As before, there is a virtual local complete intersection map
\begin{equation*}
  \iota_{\widetilde\Gamma}\colon \P \Ze_{\widetilde\Gamma} \to \P\Ze_{g,n+n',\beta,\delta},
\end{equation*}
which is a finite cover of the closure of the locus in the twisted
graph space whose elements have decorated dual graph
$\widetilde\Gamma_{\n}$, where $\widetilde\Gamma_{\n}$ is obtained
from $\widetilde{\Gamma}$ by adding $\n(v)$ additional legs of
multiplicity $\frac{1}{d}$ to each vertex $v$.

By an analogous construction, we obtain a moduli space
$\P \Xect_{\widetilde\Gamma}$ equipped with an obstruction theory and
cosection.
By cosection localization, we obtain a virtual cycle
$[\P Z^{\epsilon}_{\widetilde\Gamma}]^{\vir}$ supported on the
degeneracy locus $\P Z^{\epsilon}_{\widetilde\Gamma}$.

\begin{remark}
The definition of $\P X^{\epsilon,+}_{g,n,\beta,\delta}$ in the geometric phase is given by replacing $\vec{p}$ by $\vec{x}$ in \eqref{eq:stab1} and \eqref{eq:stab2} and replacing $P_1$ by $L_1$ in \eqref{eq:TGSdegrees} and \eqref{eq:TGSstability}.  Setting $\P X^{\epsilon,\text{ct},+}_{g,n,\beta,\delta}$ to be the locus in which $p_1(q_k) = \cdots = p_N(q_k) = 0$ for each marked point $q_k$, the above yields a virtual cycle supported on the $\P Z^{\epsilon,+}_{g,n,\beta,\delta}$.  Analogously, one obtains moduli spaces spaces $\P Z^{\epsilon,+}_{\widetilde\Gamma}$ with virtual cycles for any decorated graph.
\end{remark}

\subsection{$\C^*$-action and fixed loci}
\label{fixedloci}

There is an action of $\C^*$ on $\P\Xe_{g,n,\beta,\delta}$, given by
acting on the $z_1$-coordinate with weight $-1$.
Let $\lambda$ be the equivariant parameter, which is defined as the
dual of the weight-one representation of $\C^*$.

The fixed loci of this action are indexed by decorated graphs.  We denote such a graph by $\Lambda$ (to avoid confusion with dual graphs $\Gamma$) and decorate it as follows:
\begin{itemize}
\item Each vertex $v$ has an index $j(v) \in \{0, \infty\}$, a genus $g(v)$, and a degree $\beta(v) \in \mathbb{N}$.
\item Each edge $e$ has a degree $\delta(e) \in \mathbb{N}$.
\item Each half-edge $h$ (including the legs) has a multiplicity $m(h) \in \left\{0,\frac{1}{d}, \ldots, \frac{d-1}{d}\right\}$.
\item The legs are labeled with the numbers $\{1, \ldots, n\}$.
\end{itemize}
By the ``valence" of a vertex $v$, denoted $\text{val}(v)$, we mean the total number of incident half-edges, including legs.

The fixed locus in $\P\Xe_{g,n,\beta,\delta}$ indexed by $\Lambda$ parameterizes Landau--Ginzburg quasimaps for which:
\begin{itemize}
\item Each edge $e$ corresponds to a genus-zero component $C_e$ on which $\deg(L_2) = \delta(e)$, where $z_1$ and $z_2$ each vanish at exactly one point (the ``ramification points"), and all of the marked points, all of the nodes, and all of the degree $\beta$ is concentrated at the ramification points.  That is,
\[\deg(P_1|_{C_e}) = \sum_{q \in C_e} \text{ord}_{q}(\vec{p}),\]
so if both ramification points are special points, it follows that $\deg(P_1|_{C_e}) = 0$.
\item Each vertex $v$ for which $j(v) = 0$ (with unstable exceptional cases noted below) corresponds to a maximal sub-curve $C_v$ of $C$ over which $z_1 \equiv 0$, and each vertex $v$ for which $j(v) = \infty$ (again with unstable exceptions) corresponds to a maximal sub-curve over which $z_2 \equiv 0$.  The labels $g(v)$ and $\beta(v)$ denote the genus of $C_v$ and the degree of $P_1|_{C_v}$, respectively, and the legs incident to $v$ indicate the marked points on $C_v$.
\item A vertex $v$ is {\it unstable} if stable sub-curves of the type described above do not exist (where, as always, we interpret legs as marked points and half-edges as half-nodes).  In this case, $v$ corresponds to a single point of the component $C_e$ for each adjacent edge $e$, which may be a node at which $C_e$ meets $C_{e'}$, a marked point of $C_e$, or a basepoint on $C_e$ of order $\beta(v)$.
\item The index $m(l)$ on a leg $l$ indicates the multiplicity of $L_1$ at the corresponding marked point.
\item A half-edge $h$ incident to a stable vertex $v$ corresponds to a node at which components $C_e$ and $C_v$ meet, and $m(h)$ indicates the multiplicity of $L_1$ at the branch of the node on $C_v$.  If $v$ is unstable and hence $h$ corresponds to a single point on a component $C_e$, then $m(h)$ is the {\it negative} in $\Q/\Z$ of the multiplicity of $L_1$ at this point.
\end{itemize}
In particular, we note that the decorations at each stable vertex $v$ yield a tuple
\[\m(v) \in \left\{0, \frac{1}{d}, \ldots, \frac{d-1}{d}\right\}^{\text{val}(v)}\]
recording the multiplicities of $L_1$ at every special point of $C_v$.

The crucial observation, now, is the following.  For a stable vertex $v$ such that $j(v) = 0$, we have $z_1|_{C_v} \equiv 0$, so the stability condition \eqref{eq:stab1} implies that $\text{ord}_q(\vec{p}) \geq 1/\epsilon$ for each $q \in C_v$.  That is, the restriction of $(C; q_1, \ldots, q_n; L_1;\vec{p})$ to $C_v$ defines an element of $\Xe_{g(v), \m(v), \beta(v)}$.  On the other hand, for a stable vertex $v$ such that $j(v) = \infty$, we have $z_2|_{C_v} \equiv 0$, so the stability condition \eqref{eq:stab2} implies that $\text{ord}_q(\vec{p}) = 0$ for each $q \in C_v$.  Thus, the restriction of $(C; q_1, \ldots, q_n; L_1;\vec{p})$ to $C_v$ defines an element of $X^{\infty}_{g(v), \m(v), \beta(v)}$.  Finally, for each edge $e$, the restriction of $(\vec{x}, \vec{p})$ to $C_e$ defines a constant map to $X$ (possibly with an additional basepoint at the ramification point where $z_1=0$).

\begin{remark}
\label{rem:deltabeta}
It is important in what follows to observe that, if $C_e$ is an edge component containing a basepoint of order $\beta(e)$, then one must have $\delta(e) > \beta(e)$.  Indeed, if this is not the case, then $z_1 \equiv 0$.  Given that $z_2$ must vanish somewhere, this is impossible without violating the nondegeneracy condition in the definition of $\P\Xe_{g,n,\beta,\delta}$.
\end{remark}

Denote
\begin{equation}
  \label{eq:FLambda}
  FX_{\Lambda}:= \prod_{\substack{v \text{ stable}\\ j(v) = 0}} \Xe_{g(v), \m(v), \beta(v)} \times_{\ri X} \prod_{\text{edges }e} \ri X^{\frac 1{\delta(e)}} \times_{\ri X} \prod_{\substack{v \text{ stable}\\ j(v) = \infty}} X^{\infty}_{g(v), \m(v), \beta(v)},
\end{equation}
where $\ri X^{\frac 1{\delta(e)}}$ is the $\delta(e)$th root stack for
the line bundle $\O(-1) \otimes \C_{(\lambda)}$ on $\ri X$; the reason
for the appearance of the root stack is explained in
Section~\ref{subsubsec:edges} below.  Here, the fiber products are taken over the evaluation map at the half-node on the vertex moduli space and the composition of the evaluation map at the half-node on the edge moduli space with the involution $\iota: \ri X \rightarrow \ri X$ coming from inversion (to ensure that all nodes are balanced).

The preceding discussion implies that there is a canonical family of
$\C^*$-fixed elements of $\P\Xe_{g,n,\beta,\delta}$ over
$FX_{\Lambda}$, yielding a morphism
\begin{equation*}
  \iota_{\Lambda}\colon FX_{\Lambda} \rightarrow \P\Xe_{g,n,\beta,\delta}.
\end{equation*}
This is not exactly the inclusion of the associated fixed locus,
because elements of $\P\Xe_{g,n,\beta,\delta}$ have additional
automorphisms from permuting the components via an automorphism of
$\Lambda$, and scaling the fibers of $L_1|_{C_e}$ by $d$th roots of
unity.\footnote{The scalings of $L_1$ on adjacent components do not
  yield independent automorphisms, however, when the components meet
  at a node with nontrivial isotropy; see \cite[Proposition 1.18]{JKV}.}
In particular, $\iota_{\Gamma}$ decomposes into a finite (\'etale) map
of degree
\begin{equation}
  \label{eq:degree}
   \frac{|\text{Aut}(\Lambda)| }{\prod_{h \in \widetilde H(\Lambda)} d_{m(h)}}
\end{equation}
and a closed embedding.
Here, $\widetilde H(\Lambda)$ is a set of half-edges containing one
half-edge for each node of a generic curve in $FX_\Lambda$.
This means that $\widetilde H(\Lambda)$ contains all non-leg
half-edges at stable vertices and exactly one half-edge at each
non-legged valence-two vertex with genus zero and $\beta$-degree zero.
Furthermore, we define $d_{m(h)} = d/\gcd(d \cdot m(h),d)$, which does
not depend on which half-edge at an unstable vertex is chosen.

The $\C^*$-action preserves the locus $\P \Xect_{g,n,\beta,\delta} \subset \P \Xe_{g,n,\beta,\delta}$, and we can similarly define $FX^{\ct}_{\Lambda} \subset FX_{\Lambda}$ by the requirement that $x_1(q_k) = \cdots = x_M(q_k) = 0$ whenever $q_k$ corresponds to a leg of $\Lambda$.  The perfect obstruction theory on $\P \Xect_{g,n,\beta,\delta}$ and the cosection $\widetilde{\sigma}$ are both $\C^*$-equivariant, so by \cite{CKL}, each fixed locus admits a cosection-localized virtual class, and the virtual localization formula expresses $[\P\Ze_{g,n,\beta,\delta}]^{\vir}$ in terms of these contributions.  We set
\begin{equation}
  \label{eq:FZLambda}
  FZ_{\Lambda}:= \prod_{\substack{v \text{ stable}\\ j(v) = 0}} \Ze_{g(v), \m(v), \beta(v)} \times_{\ri X^\ct} \prod_{\text{edges }e} \ri X^{\ct, \frac 1{\delta(e)}} \times_{\ri X^{\ct}} \prod_{\substack{v \text{ stable}\\ j(v) = \infty}} Z^{\infty}_{g(v), \m(v), \beta(v)},
\end{equation}
where, again, the fiber product is taken over evaluation on the vertex side and evaluation followed by inversion on the edge side.  We denote by $[FZ_{\Lambda}]^{\vir}$ the pullback of the
cosection-localized virtual class on the fixed locus associated to
$\Lambda$ under the \'etale map.

More generally, one can do all of this on each of the stacks
$\P \Ze_{\widetilde\Gamma}$ associated to a decorated dual graph
$(\widetilde\Gamma, v_{\bullet}, \n)$.  That is, let $\widetilde\Gamma$ be a stable decorated dual graph as in Section~\ref{subsec:tgs}, and let $\Lambda$ be any graph obtained from $\widetilde\Gamma$ by replacing each vertex $v$ by a localization graph for $\P X^{\epsilon}_{g(v),\m'(v),\beta(v),D(v)}$.  For each such $\Lambda$, there is a fixed locus in $\P X^{\epsilon}_{\widetilde\Gamma}$, and for the same fiber product $FX_{\Lambda}^{\ct}$ defined above, we have
a morphism
\begin{equation*}
  \iota_{\widetilde\Gamma, \Lambda}\colon FX^{\ct}_\Lambda \rightarrow \P \Xe_{\widetilde\Gamma},
\end{equation*}
whose image is the fixed locus associated to $\Lambda$.

\begin{remark}
The fixed loci have an analogous structure in the geometric phase; see \cite[Section 3.2]{CJR}.
\end{remark}

\subsection{Localization contributions}
\label{subsec:contributions}

The virtual localization formula, first proved by
Graber--Pandharipande \cite{GP} and adapted to the setting of
cosection-localized virtual classes by Chang--Kiem--Li \cite{CKL},
expresses $[\P\Ze_{g,n,\beta,\delta}]^{\vir}$ in terms of contributions from
each localization graph:
\begin{equation}
  \label{eq:virloc}
  [\P\Ze_{g,n,\beta,\delta}]^{\vir}=  \sum_{\Lambda} \frac{\prod_{h \in \widetilde H(\Lambda)} d_{m(h)}}{|\text{Aut}(\Lambda)|} \iota_{\Lambda, *} \left(\frac{[FZ_{\Lambda}]^{\vir}}{e(N_{\Lambda}^{\vir})}\right).
\end{equation}
There is a technical assumption necessary in order to be able to apply
this formula: the restriction of the perfect obstruction theory to
each fixed locus has to admit a global resolution by a perfect
complex.  In our case, this is clear from the analysis of the perfect
obstruction theory in the following sections.
  
More generally, all of this can be done for each dual graph $\widetilde\Gamma$.  Note that if $\Lambda$ is a localization graph contributing to $[\P\Ze_{\widetilde\Gamma}]^{\vir}$, then the space $FZ_{\Lambda}$ appears in the localization both for $[\P\Ze_{g,n,\beta,\delta}]^{\vir}$ and for
$[\P\Ze_{\widetilde\Gamma}]^{\vir}$.  However, its contribution is {\it a priori} different: in the localization for $[\P\Ze_{g,n,\beta,\delta}]^{\vir}$, we have a virtual cycle $[FZ_{\Lambda}]^{\vir}$ as defined above via the obstruction theory on $\P \Xect_{g,n,\beta,\delta}$, whereas in the localization for $[\P\Ze_{\widetilde\Gamma}]^{\vir}$, we have a virtual cycle that we denote $[FZ_{\Lambda}]^{\vir}_{\widetilde\Gamma}$, defined via the obstruction theory on $\P \Xect_{\widetilde\Gamma}$.  Similarly, we have virtual normal bundles $N^{\vir}_{\Lambda}$ for the fixed locus associated to $\Lambda$ inside $\P \Xect_{g,n,\beta,\delta}$, and
$N^{\vir}_{\Lambda,\widetilde\Gamma}$ for the fixed locus inside $\P \Xect_{\widetilde\Gamma}$.

In fact, however, these two contributions are closely related:
\begin{lemma}
  \label{lem:comparison}
  We have
  \begin{equation*}
    [FZ_{\Lambda}]^{\vir}
    = [FZ_{\Lambda}]_{\widetilde\Gamma}^{\vir}
  \end{equation*}
  and
  \begin{equation*}
    e(N_{\Lambda}^{\vir})
    = e(N_{\Lambda, \widetilde\Gamma}^{\vir}) \prod_{e \in E(\widetilde\Gamma)} c_1(T_{e, 1} \otimes T_{e, 2}),
  \end{equation*}
  where, for each edge $e \in E(\widetilde\Gamma)$, $T_{e, 1}$ and $T_{e,2}$ denote the relative tangent line bundles
  at the corresponding node of the \emph{orbifold} curve.
\end{lemma}
\begin{proof}
  We have the following commutative diagram, whose right square is
  cartesian:
  \begin{center}
    \begin{tikzpicture}
      \matrix (m) [matrix of math nodes,row sep=2em,column sep=4em,minimum width=2em] {
      FX^{\ct}_\Lambda & \P \Xect_{\widetilde\Gamma} & \P \Xect_{g,n,\beta,\delta} \\
      & D_{\widetilde\Gamma} & D_{g, n, \beta}. \\};
    \path[-stealth]
      (m-1-1) edge node[above] {$\iota_{\widetilde\Gamma, \Lambda}$} (m-1-2)
      (m-1-2) edge node[left] {$\pi$} (m-2-2) edge node[above] {$\iota_{\widetilde\Gamma}$} (m-1-3)
      (m-2-2) edge node[above] {$\iota$} (m-2-3)
      (m-1-3) edge (m-2-3);
    \end{tikzpicture}
  \end{center}
  By definition, the perfect obstruction theory on
  $\P \Xect_{\widetilde\Gamma}$ relative to
  $D_{\widetilde\Gamma}$ is
  $\iota^* \widetilde{\mathbb E}^\bullet$.
  By \cite[Lemma~4.4]{CLL}, the map $\iota$ is a finite local complete
  intersection morphism, so its cotangent complex
  $\mathbb L_\iota^\bullet$ is concentrated in degree $-1$, and it is
  quasi-isomorphic to the complex 
  \begin{equation*}
    \left(\bigoplus_{e \in E(\widetilde\Gamma)} T_{e, 1} \otimes T_{e, 2}\right)^\vee \to 0
  \end{equation*}
  of vector bundles.
  Thus, the absolute perfect obstruction theories
  $\mathbb E_{\widetilde\Gamma, \mathrm{abs}}^\bullet$ and
  $\widetilde{\mathbb E}_{\mathrm{abs}}^\bullet$ fit into a
  distinguished triangle
  \begin{equation*}
    \pi^* (\mathbb{L}^{\bullet}_{\iota})[-1] \xrightarrow{+1} \widetilde{\mathbb E}_{\mathrm{abs}}^\bullet \to \mathbb E_{\widetilde\Gamma, \mathrm{abs}}^\bullet \to \pi^* (\mathbb{L}^{\bullet}_{\iota}).
  \end{equation*}
  The virtual classes $[FZ_{\Lambda}]^{\vir}$ and
  $[FZ_{\Lambda}]_{\widetilde\Gamma}^{\vir}$ are defined by the $\C^*$-fixed part of
  $\iota_{\widetilde\Gamma, \Lambda}^* \left(\widetilde{\mathbb
    E}_{\mathrm{abs}}^\bullet\right)$ and
  $\iota_{\widetilde\Gamma, \Lambda}^* \left(\mathbb E_{\widetilde\Gamma,
    \mathrm{abs}}^\bullet\right)$, respectively, and the virtual normal bundles are defined by the
  corresponding moving parts.
  The lemma therefore follows from the fact that, as in the calculation in Section~\ref{subsubsec:nodes} below, the line bundles $T_{e, 1} \otimes T_{e, 2}$ are moving.
\end{proof}

The goal of the remainder of this subsection is to compute the
contributions to \eqref{eq:virloc} of each graph $\Lambda$ explicitly. 
(All of this can also be done for each
$[\P \Ze_{\widetilde\Gamma}]^{\vir}$, but by Lemma
\ref{lem:comparison}, the contributions of each localization graph
$\Lambda$ in that case are straightforward to deduce from the
corresponding contributions to $[\P\Ze_{g,n,\beta,\delta}]^{\vir}$.)  We first concentrate on the virtual normal bundle, identifying its equivariant Euler class with an explicit cohomology class pulled back from $FZ_{\Lambda}$.    To do so, we apply the normalization exact sequence to express the moving part of the perfect obstruction theory in terms of vertex, edge, and node factors.  Then, in Section~\ref{subsubsec:virtual}, we study the fixed part of the perfect obstruction theory, decomposing the cosection-localized virtual class along certain nodes.  Finally, in Section~\ref{subsubsec:comparison}, we push forward $[\P Z^{\epsilon}_{g,n,\beta,\delta}]^{\vir}$ to $Z^{\epsilon}_{g,n,\beta}$ by forgetting $L_2$ and $\vec{z}$, and we use the above calculations to express the contribution of each graph $\Lambda$ to the resulting pushforward in terms of a genus-$g$ contribution and contributions from certain trees of rational components.  Crucially, we observe that the genus-$0$ contributions depend in only a very controlled way on the particular GIT quotient $X$ at hand, which sets the stage for a more explicit study of those contributions in the following subsection.
  
\subsubsection{Vertex contributions}
\label{subsubsec:vertices}

For a vertex $v$ of $\Lambda$, the sections $\vec x$ and $\vec p$ are
$\C^*$-fixed, so the deformations of these sections, together with the
deformations of the curve $C_v$ and the line bundle $L_1|_{C_v}$, are
part of the $\C^*$-fixed part of the perfect obstruction theory.
The map $f$ is moving, and its contribution to the inverse Euler class of the virtual normal bundle is
\begin{equation*}
  \frac{1}{e^{\C^*}(R\pi_* f^* T_{\P_v/\cC_v})} = \frac{1}{e^{\C^*}(R\pi_* N_{j(v)})},
\end{equation*}
where 
\begin{equation*}
  N_0 = \cP_1^\vee \otimes \C_{(\lambda)}, \qquad N_\infty = \cP_1 \otimes \C_{(-\lambda)}
\end{equation*}
are the normal bundles to the zero and infinity sections,
respectively, in the projectivized bundle
$\P_v := \P((\cP_1^{\vee} \otimes \C_{(\lambda)}) \oplus \O_{\cC_v})$
on the universal curve $\cC_v$.

\subsubsection{Edge contributions}
\label{subsubsec:edges}

For an edge $e$, the corresponding factor of $FZ_\Lambda$ is
isomorphic to the $\delta(e)$th root stack
$\widetilde{X}^{\text{ct}}_e := \mathcal IX^{\ct, \frac 1{\delta(e)}}$
of a component $X_e^\ct$ of the \emph{ordinary} inertia stack
$\mathcal IX^\ct$ with respect to the line bundle
$\O_{X_e^\ct}(-1) \otimes \C_{(\lambda)}$.
By virtue of its universal property, on $\widetilde{X}^{\text{ct}}_e$,
there is a universal $\delta(e)$th root $\mathcal R$ of
$\O_{\widetilde{X}^{\text{ct}}_e}(-1) \otimes \C_{(\lambda)}$.
With this root, the universal family over $X_e^\ct$ takes the form
\begin{equation*}
  \xymatrix{
    \cC_e := \P(\mathcal R \oplus \O_{\widetilde{X}^{\text{ct}}_e}) \ar^-{f}[r]\ar_{\pi}[d] & \P((\cP_1^\vee \otimes \C_{(\lambda)})\oplus \O_{\widetilde{X}^{\text{ct}}_e}) =: \P_e \\
    \widetilde{X}^{\text{ct}}_e,
  }
\end{equation*}
at least in the case that the multiplicity of the edge is zero.
In general, $\cC_e$ has additional orbifold structure at the
zero and infinity section; however, since the computations of this section are of cohomological nature, we can ignore this orbifold structure and view $\cC_e$ as above.

Let us make the universal sections and universal map explicit.
By the universal property of the projectivized bundle $\cC_e$ over
$\widetilde{X}^{\text{ct}}_e$, there is a line bundle $\O_{\cC_e}(1)$ together with
universal sections
\begin{equation*}
  (x, y) \in H^0((\O_{\cC_e}(1) \otimes \mathcal R) \oplus \O_{\cC_e}(1))
\end{equation*}
such that $x$ (respectively, $y$) vanishes precisely at the zero
section (respectively, at the infinity section) of $\cC_e$ with
order one.
The universal line bundle $\cP_1$ is given by
\begin{equation*}
  \cP_1 = \pi^*\O_{\widetilde{X}^{\text{ct}}_e}(1) \otimes \O_{\cC_e}(\beta(e)[0]) = \mathcal R^{\delta(e)} \otimes \C_{(\lambda)},
\end{equation*}
where $[0]$ is the zero section in the projectivized bundle and
\begin{equation*}
  \beta(e) :=
  \begin{cases}
    \beta(v) & \text{if there is a genus zero, valence one vertex }v\text{ at }e, \\
    0 & \text{otherwise.}
  \end{cases}
\end{equation*}
Note that the universal section $x$ gives an isomorphism
\begin{equation*}
  \O_{\cC_e}([0]) \cong \O_{\cC_e}(1) \otimes \mathcal R.
\end{equation*}

The universal map $f$ is given by the sections
\begin{equation*}
  (\zeta_1, \zeta_2) = (x^{\delta(e) - \beta(e)}, y^{\delta(e)}) \in H^0((\L_2 \otimes \cP_1^\vee \otimes \C_{(\lambda)}) \oplus \L_2),
\end{equation*}
where
\begin{equation*}
  \L_2 = \O_{\cC_e}(\delta(e)).
\end{equation*}
Note that
\begin{equation*}
  f^* T_{\P_e/\cC_e} = \O_{\cC_e}((\delta(e) - \beta(e))[0] + \delta(e)[\infty]),
\end{equation*}
where $[\infty]$ denotes the infinity section of $\cC_e$.

It is useful for what follows to introduce the notation
\begin{equation*}
  \lambda_{j} = \begin{cases} \lambda -H & \text{ if } j=0\\ -\lambda + H & \text{ if } j=\infty\end{cases}
\end{equation*}
for the relative tangent line classes at the zero and infinity section
of $\P_e$.
We interpret $\lambda_j$ as an equivariant cohomology class on
$\widetilde{X}^{\text{ct}}_e$, where $H$ denotes the hyperplane class.

Equipped with this notation, we turn to the contributions of the edge.  First of all, there are contributions from the sections $\vec x$ and $\vec p$.  For these, it is useful to consider both the moving and fixed parts of the obstruction theory, since the total contribution to the localization from these sections is essentially identical to the corresponding contribution to the graph-space localization used to define the $J$-function.  There are two differences between this setting and the $J$-function, however.    First, we must divide the $J$-function by a factor of $z$ to cancel
the contribution $-z^{-1}$ in the localization on the ordinary graph
space and the prefactor of $-z^2$; see Remark~\ref{rem:prefactor}.
Second, the tangent space at the zero section of the universal curve
of the ordinary graph space has first Chern class $z$, while in our
situation the corresponding Chern class is $\lambda_0/\delta(e)$.
Thus, one sees that $z$ must be replaced in the $J$-function by
$\lambda_0/\delta(e)$. 
In all, then, we thus far have a contribution of
\begin{equation}
  \label{eq:xcont}
  \frac 1{d_{m(h)}} \bigg[z^{-1} J^{\epsilon}(q,z)\big|_{z = \lambda_0/\delta(e)}\bigg]_{q^{\beta(e)}}.
\end{equation}

Furthermore, there are contributions from deformations of the map $f$.
The vector bundle
\begin{equation*}
  R^0\pi_*\big(f^* T_{\P_e/\cC_e}\big)
  = R^0\pi_*\bigg(\O_{\cC_e}\big((\delta(e) - \beta(e))[0] + \delta(e)[\infty]\big)\bigg)
\end{equation*}
has a trivial, $\C^*$-fixed factor $R^0\pi_* \O$.
The complement of this factor is moving, and the inverse of its Euler
class can be written as
\begin{equation*}
  \frac 1{\prod_{b = 1}^{\delta(e) - \beta(e)} \frac{b\lambda_0}{\delta(e)} \prod_{b = 1}^{\delta(e)} \frac{b\lambda_\infty}{\delta(e)}}
  = \frac{\prod_{b = \delta(e) - \beta(e) + 1}^{\delta(e)} \frac{b\lambda_0}{\delta(e)}}{\prod_{b = 1}^{\delta(e)} \frac{b\lambda_0}{\delta(e)} \prod_{b = 1}^{\delta(e)} \frac{b\lambda_\infty}{\delta(e)}}.
\end{equation*}
In the case where there is no basepoint, so $\beta(e) =0$, only the
denominator of the above appears.
The numerator, on the other hand, can be rewritten as
\begin{equation*}
  \prod_{b = \delta(e) - \beta(e) + 1}^{\delta(e)} \frac{b\lambda_0}{\delta(e)} = \prod_{b = 0}^{\beta(e) - 1} \left(\lambda_0 - \frac{b\lambda_0}{\delta(e)}\right).
\end{equation*}
Using
\begin{multline*}
  e^{\C^*}(-R\pi_*(\cP_1^\vee \otimes \C_{(\lambda)}))
  = e^{\C^*}((R^0\pi_*(\omega_{\pi} \otimes \cP_1))^\vee \otimes \C_{(\lambda)}) \\
  = e^{\C^*}\bigg(\bigg(R^0\pi_* \O\big((\beta(e) - 1)[0] - [\infty]\big)\bigg)^\vee \otimes \O_{\widetilde Z_e}(-1) \otimes \C_{(\lambda)}\bigg)
  = \lambda_0^{-1} \prod_{b = 0}^{\beta(e) - 1} \left(\lambda_0 - \frac{b\lambda_0}{\delta(e)}\right),
\end{multline*}
we can write
\begin{equation}
\label{eq:zcont}
  \frac{1}{e^{\C^*}(R^0\pi_*(f^* T_{\P_e/\cC_e}))} = \frac{\lambda_0 e^{\C^*}(-R\pi_*(\cP_1^\vee \otimes \C_{(\lambda)}))}{\prod_{b = 1}^{\delta(e)} \frac{b\lambda_0}{\delta(e)} \prod_{b = 1}^{\delta(e)} \frac{b\lambda_\infty}{\delta(e)}},
\end{equation}
which is the contribution to the inverse Euler class of the virtual normal bundle coming from the sections $\vec{z}$.

Together, \eqref{eq:xcont} and \eqref{eq:zcont} yield a contribution of
\begin{equation*}
  \frac{\delta(e)}{\lambda_0 d_{m(h)}} \bigg[J^{\epsilon, \tw}(q,z)\big|_{z = \lambda_0/\delta(e)}\bigg]_{q^{\beta(e)}}
  \cdot \frac 1{\prod_{b = 1}^{\delta(e)} \frac{b\lambda_0}{\delta(e)} \prod_{b = 1}^{\delta(e)} \frac{b\lambda_\infty}{\delta(e)}}
\end{equation*}
to the localization.
  
Finally, there are additional contributions from automorphisms of
$C_e$.
These have nontrivial torus weight only if there is an unmarked
ramification point, so the edge $e$ must be incident to an unstable
vertex $v$ of valence $1$.
The contribution from such automorphisms to the inverse Euler class of
the virtual normal bundle is the tangent line class at the point
corresponding to $v$, which is
\begin{equation*}
  \frac{\lambda_{j(v)}}{\delta(e)}.
\end{equation*}

\subsubsection{Node contributions}
\label{subsubsec:nodes}

In the normalization exact sequence for
$R\pi_* f^* T_{\P((\cP_1^\vee \otimes \C_{(\lambda)}) \oplus
  \O)/\cC}$, a node corresponding to a vertex $v$ contributes
\begin{equation*}
  e^{\C^*}(N_{j(v)}) = \lambda_{j(v)}
\end{equation*}
to the inverse Euler class of the virtual normal bundle.
The deformations in $D_{g,n,\beta}$ smoothing the
node also give a contribution of
\begin{equation*}
  \frac{1}{e^{\C^*}(T_1) + e^{\C^*}(T_2)} = \frac{d_{m(h)}}{e^{\C^*}(\overline T_1) + e^{\C^*}(\overline T_2)},
\end{equation*}
where $T_1$ and $T_2$ denote the tangent lines of the orbifold curve
at the two branches of the node, $\overline T_1$ and $\overline T_2$
denote the corresponding tangent lines pulled back from the underlying
coarse curve, and where $h$ is any of the half-edges of $e$.
This inverse makes sense since at least one of the $T_i$ (say $T_1$)
corresponds to an edge $e$, so we have
\begin{equation*}
  e^{\C^*}(T_1) = \frac{\lambda_{j(v)}}{\delta(e)}.
\end{equation*}

\subsubsection{Virtual cycle}
\label{subsubsec:virtual}

We now look at the fixed part of the perfect obstruction theory in
more detail, decomposing it along certain nodes.

We first introduce some notation and terminology.  Let $\Lambda$ be a localization graph, and let $\Gamma$ be the prestable dual graph of a generic element of the associated fixed locus in $\P X^{\epsilon}_{g,n,\beta,\delta}$.  We call a vertex $v$ of $\Gamma$ \emph{very stable} if the
corresponding irreducible component is not contracted under the
forgetful map
$\P \Xe_{g,n,\beta,\delta} \rightarrow \Xe_{g,n,\beta}$.
Let $\Gamma'$ be obtained from $\Gamma$ according to the forgetful map
$\P \Xe_{g,n,\beta,\delta} \rightarrow \Xe_{g,n,\beta}$, with an
extra leg on each very stable vertex wherever a rational tail is
contracted.

Denote by $V_{E(\Lambda)}(\Gamma) \subset V(\Gamma)$ the vertices whose corresponding irreducible component is an edge component of $\Lambda$, and let
\[T = \{v \in V(\Gamma) \; | \; v \text{ not very stable, } v \notin V_{E(\Lambda)}(\Gamma), \; \beta(v) > 0\}.\]
For each $v \in T$, define
$\Xecto_{0, n(v), \beta(v)}$ as the locus inside
$\Xe_{0, n(v), \beta(v)}$ where $x_1(q_k) = \cdots = x_M(q_k) = 0$ for
all $k \in \{1, \dotsc, n(v)\}$ except for the unique leg that is
closest to one of the very stable vertices.\footnote{In fact, $X^{\epsilon,\ct-1}_{0,n(v),\beta(v)} = \Ze_{0,n(v),\beta(v)}$ in the Landau--Ginzburg phase, but we use different notation to clarify the parallel argument in the geometric phase.}  Define $FZ'_{\Lambda}$ analogously to $FZ_{\Lambda}$, but replacing the factor of $\Ze_{0,n(v),\beta(v)}$ by $\Xecto_{0,n(v),\beta(v)}$ for each $v \in T$.  Let
\[\iota: FZ_{\Lambda} \hookrightarrow FZ'_{\Lambda}\]
be the inclusion.

\begin{figure}
\begin{subfigure}{\textwidth}
\centering
\xymatrix@C-=0.35cm@R-=0.35cm{
\Gamma & &*++[o][F-]{\;}\ar@{-}[d]&&&*++[o][F-]{\;}\ar@{-}[dl]\ar@{-}[dr]&&&*++[o][F-]{\;}\ar@{-}[d]&&*++[o][F-]{\;}\ar@{-}[d]&\\
 & \beta=0\hspace{-0.25cm}& *++[o][F=]{\;}\ar@{-}[d]\ar@{-}[r] & *++[o][F-]{\;}\ar@{-}[r] & *++[o][F=]{1}\ar@{-}[dr] && *++[o][F=]{1}\ar@{-}[dl]\ar@{-}[r] & *++[o][F-]{\;}\ar@{-}[r] & *++[o][F=]{\;}\ar@{-}[r] & *++[o][F-]{\;}\ar@{-}[r] & *++[o][F=]{\;}\ar@{-}[r] & *++[o][F-]{\;}\\
& & *++[o][F-]{\;} &&& *++[o][F-]{\;} &&&\beta>0 && \beta>0 &
}

\end{subfigure}

\vspace{1cm}

\begin{subfigure}[r]{.4\textwidth}
\begin{xy}
\xymatrix@C-=0.3cm@R-=0.3cm{
\Gamma' &&*++[o][F=]{1}\ar@{-}[l]\ar@/_1pc/@{-}[rr]\ar@/^1pc/@{-}[rr]&& *++[o][F=]{1}\ar@{-}[r] &
}
\end{xy}
\end{subfigure}
\hspace{2cm}
\begin{subfigure}[r]{.4\textwidth}
\begin{xy}
\xymatrix@C-=0.3cm@R-=0.3cm{
&&&&&\\
&*++[o][F=]{1}\ar@{-}[l]\ar@/_1pc/@{-}[rr]\ar@/^1pc/@{-}[rr]&& *++[o][F=]{1}\ar@{-}[r] &*++[o][F=]{\;}\ar@{-}[r]\ar@{-}[u]   & *++[o][F=]{\;}\ar@{-}[r]\ar@{-}[u]   & \\
&&&&&
}
\end{xy}
\end{subfigure}
\caption{On the first line, the dual graph $\Gamma$ associated to a localization graph $\Lambda$, where the number within a vertex indicates its genus and empty vertices are understood as genus zero.   Vertices with a double circle correspond to vertex components of $\Lambda$, vertices with a single circle correspond to edge components of $\Lambda$, and the genus-zero vertex components are labeled according to whether their $\beta$-degree equals zero.  The second line shows the graph $\Gamma'$ and the geometry represented by the fiber product \eqref{eq:fiberproduct}.}
\label{fig:1}
\end{figure}
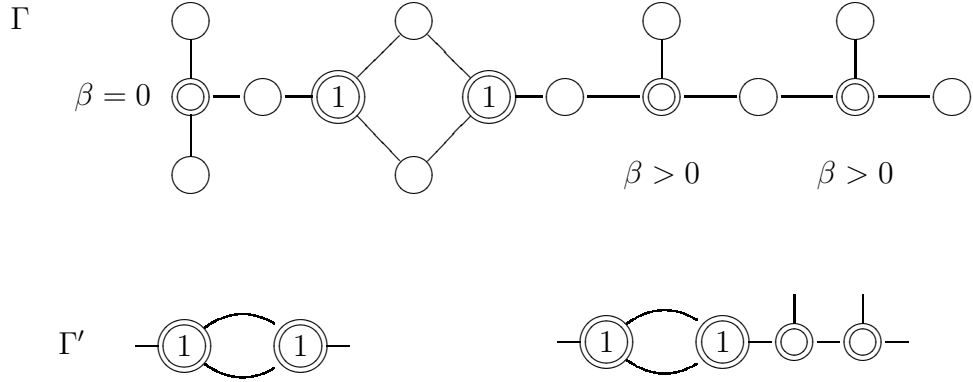

If
\begin{equation}
\label{eq:fiberproduct}
  \widetilde p\colon FZ'_\Lambda \to \Xect_{\Gamma'} \times_{(\ri X^\ct)^{|T|}} \prod_{v \in T} \Xecto_{0, n(v), \beta(v)}
\end{equation}
is defined by contracting vertices that are neither very stable nor in
$T$ (see Figure~\ref{fig:1}), and 
\begin{equation*}
  \widetilde\Delta\colon (\ri X^\ct)^{|T|} \to (\ri X^\ct)^{|T|} \times (\ri X^\ct)^{|T|}
\end{equation*}
is the diagonal map, then the aim of this section is to prove the following formula:
\begin{lemma}
\label{lem:vclass}
  \begin{equation*}
    \iota_* [FZ_\Lambda]^\vir
    = \widetilde p^* \widetilde\Delta^! \left([Z^{\epsilon}_{\Gamma'}]^\vir \times \prod_{v \in T} [\Xecto_{0, n(v), \beta(v)}]^\vir\right)
  \end{equation*}
\end{lemma}

The first step toward proving Lemma~\ref{lem:vclass} is to re-write $[FZ_{\Lambda}]^{\vir}$ via a more convenient perfect obstruction theory.  Thus far, $[FZ_{\Lambda}]^{\vir}$ has been defined using the fixed part of (the absolute perfect obstruction theory induced by) $\widetilde{\mathbb{E}}^{\bullet}$, which is a perfect obstruction theory for $FX^{\ct}_{\Lambda}$ relative to $D_{g,n,\beta}$; we denote this fixed part by $\mathbb E_{\mathrm{old}}^\bullet$.  On  the other hand, if\begin{equation*}
  \mathcal E := \bigoplus_{i = 1}^M (\L^{\otimes w_i} \otimes \O(-{\textstyle\sum_{k = 1}^n} \Delta_k)) \oplus \bigoplus_{j = 1}^N (\L^{\otimes -d}\otimes \omega_{\pi,\log}),
\end{equation*}
in which $\cC$ is the universal curve on $FX^{\ct}_{\Lambda}$ and $\mathcal{L}$ the universal line bundle, then the complex $(R\pi_*\mathcal{E})^{\vee}$ is a perfect obstruction theory for $FX^{\ct}_{\Lambda}$ relative to the smooth Artin stack $D_{\Lambda}$ of curves $C$ in 
\begin{equation*}
  \prod_{\substack{v \text{ stable}\\ j(v) = 0}} \fM^\tw_{g(v), n(v)} \times \prod_{\text{edges }e} \{\mathrm{point}\} \times \prod_{\substack{v \text{ stable}\\ j(v) = \infty}} \fM^\tw_{g(v), n(v)}
\end{equation*}
together with a line bundle $L$.
Here, $\fM^\tw_{g, n}$ denotes the stack of prestable twisted curves.
Let $\mathbb E_{\mathrm{new}}^\bullet$ be the fixed part of the
corresponding absolute perfect obstruction
theory.\footnote{$R\pi_* \mathcal E$ has moving factors corresponding
  to edges with basepoints.}
\begin{lemma}
  $\mathbb E_{\mathrm{old}}^\bullet$ and
  $\mathbb E_{\mathrm{new}}^\bullet$ are quasi-isomorphic.
\end{lemma}
\begin{proof}
  This follows from the analysis of the virtual normal bundle in the
  previous sections.
  As we have seen, the curve deformation factors corresponding to smoothing the
  nodes are moving.
  Therefore, for $\mathbb E_{\mathrm{old}}^\bullet$, we can work
  relative to the Artin stack $D_\Gamma$ instead of $D_{g, n, \beta}$.
  We have also seen that, if $\P:= \P(\mathcal{P}_1^{\vee} \oplus \O)$, then the factor $R\pi_* f^* \mathbb{L}_{\P/\cC}$ in $\widetilde{\mathbb{E}}^{\bullet}$ is moving
  except for one trivial factor for each edge in $\Lambda$.
  This factor is identified with the corresponding edge automorphism
  factor in the cotangent complex of $D_\Gamma$.
  Finally, the remaining edge automorphisms are moving.
  Therefore, working relative to $D_\Lambda$ and removing the factor
  $R\pi_* f^* \mathbb{L}_{\P/\cC}$ from $\widetilde{\mathbb{E}}^{\bullet}$
  gives $\mathbb E_{\mathrm{old}}^\bullet$, which is the same as the description of
  $\mathbb E_{\mathrm{new}}^\bullet$.
\end{proof}

To continue, we separate the dual graph $\Gamma$ corresponding to
$\Lambda$ into two parts $\Gamma_g$ and $\Gamma_0$, where $\Gamma_g$
is formed by all very stable vertices together with all chains of
rational components connecting them, and $\Gamma_0$ consists of the
remaining vertices of $\Gamma$, which form a disjoint union of trees
of rational components.
Correspondingly, we can decompose the universal curve
\begin{equation*}
  \cC = \cC_g \cup \cC_0.
\end{equation*}
Let $\Delta = \mathcal{C}_g \cap \mathcal{C}_0$, and define the following subsheaf:
\begin{equation*}
  \mathcal E' := \bigoplus_{i = 1}^M (\L^{\otimes w_i} \otimes \O(-\Delta - {\textstyle\sum_{k = 1}^n} \Delta_k))|_{\cC_g} \oplus \bigoplus_{j = 1}^N ( \L^{\otimes -d} \otimes \omega_{\pi,\log}) \subset \mathcal E.
\end{equation*}
The cokernel of the inclusion $\mathcal{E}' \hookrightarrow \mathcal{E}$ is given by
\begin{equation*}
  \mathcal E'' := \bigoplus_{i = 1}^M (\L^{\otimes w_i} \otimes \O(-{\textstyle\sum_{k = 1}^n} \Delta_k))|_{\cC_0}.
\end{equation*}
Note that by Remark~\ref{rem:genuszero}, we have
$\pi_* \mathcal E'' = 0$, so that
$\pi_* \mathcal E' \cong \pi_* \mathcal E$.
Therefore, $(R\pi_* \mathcal E')^\vee$ provides an alternative perfect
obstruction theory for $FX^\ct_\Lambda$ relative to $D_\Lambda$; moreover, the cosection $\sigma$ restricts to a cosection $\sigma'$ on this alternative perfect obstruction theory.  The degeneracy locus of $\sigma'$ is precisely the locus $FZ_{\Lambda}'$ defined above, so we denote by $[FZ_{\Lambda}']^{\vir}$ the induced cosection-localized virtual class.  Cosection-localized pullback then implies the following:
\begin{lemma}
  \begin{equation*}
    \iota_* [FZ_\Lambda]^\vir = e(R^1 \pi_* \mathcal E'')^{\mathrm{fix}} \cap [FZ'_\Lambda]^{\vir}
  \end{equation*}
\end{lemma}

The vector bundle $R^1 \pi_* \mathcal E''$ can be easily decomposed.
First, $R^1 \pi_* \mathcal E''$ has a summand for each connected
component of $\cC_0$.
Then, using short exact sequences similar to
$0 \to \mathcal E' \to \mathcal E \to \mathcal E'' \to 0$, we can
inductively split $e(R^1 \pi_* \mathcal E'')^{\mathrm{fix}}$ into a
factor for each irreducible component of $\cC_0$.
Each time we split at a node, we make the choice to twist down by
$\Delta$ on the side closest to the trunk of the tree.
Now, we can check that the only fixed factors of
$e(R^1 \pi_* \mathcal E'')$ are on the vertices of $\Gamma$ in the set
$T$.
Therefore, we can write
\begin{equation*}
  e(R^1 \pi_* \mathcal E'')^{\mathrm{fix}} = \prod_{v \in T} e(R^1\pi_* \mathcal E''_v),
\end{equation*}
where $\mathcal E''_v$ denotes the restriction of $\mathcal E''$ to the
component of $\cC_0$ corresponding to $v$ and twisted down at all but
the node closest to $\cC_g$.

The factor of $[FZ'_{\Lambda}]^{\vir}$, on the other hand, is pulled back under $\widetilde{p}$.  More specifically, let $\Lambda'$ denote the lower-right graph in Figure~\ref{fig:1}, and let $FX^{\ct}_{\Lambda'}$ be the target of $\widetilde{p}$.  Then we have a cartesian diagram
\begin{equation*}
  \xymatrix{
    FX^{\ct}_{\Lambda} \ar[r]^{\widetilde p} \ar[d] & FX^{\ct}_{\Lambda'} \ar[d] \\
    D_\Lambda \ar[r] & D_{\Lambda'},
    },
\end{equation*}
where $D_{\Lambda'}$ denotes the stack of curves $C$ in
\begin{equation*}
  \prod_{v \text{ very stable}} \fM^\tw_{g(v), \m(v)} \times \prod_{v \in T} \fM^\tw_{0, n(v)}
\end{equation*}
together with a line bundle $L$.  Note that the lower horizontal map is smooth.

Since the map $p\colon \cC_g \to \cC'_g$ contracting chains of
rational components of $\cC_g$ is log \'etale (in the sense that
$p^* \omega_{\cC'_g, \log} = \omega_{\cC_g, \log}$ and therefore
also $p^* \L = \L$), the perfect obstruction theory
$(R\pi_* \mathcal E')^\vee$ of $FX^{\ct}_\Lambda$ relative to
$D_\Lambda$ is pulled back from the analogously-defined perfect
obstruction theory of $FX^{\ct}_{\Lambda'}$ relative to the stack
$D_{\Lambda'}$.
The cosection $\sigma'$ is also pulled back from the analogously-defined cosection $\sigma'$ of $FX^\ct_{\Lambda'}$.
Therefore,
\begin{equation}
  \label{eq:contract}
  \iota_* [FZ_\Lambda]^{\vir} = \prod_{v \in T} e(R^1\pi_* \mathcal E''_v) \cap \widetilde p^* [FX^{\ct}_{\Lambda'}]^{\vir}_{\sigma'}.
\end{equation}

Finally, we split $[FX^{\ct}_{\Lambda'}]^{\vir}_{\sigma'}$ along nodes joining the $T$-components to the rest of the curve.  For this, first note that the universal curve $\cC'$ over $FX^{\ct}_{\Lambda'}$ is of the form
\begin{equation*}
  \cC' = \cC'_g \cup \bigcup_{v \in T} \cC_v,
\end{equation*}
where $\cC_v$ is the rational component corresponding to $v$.
The curves intersect in a set $\Delta'$ of nodes, with $|\Delta'| = |T|$.
We form the cartesian diagram
\begin{equation*}
  \xymatrix{
    FX^{\ct}_{\Lambda'} \ar[r] \ar[d] & FX^{\ct'}_{\Lambda'} \ar[d] \\
    (\ri X^\ct)^{|T|} \ar[r]^-{\widetilde\Delta} & (\ri X^\ct)^{|T|} \times (\ri X^\ct)^{|T|},
    }
\end{equation*}
where
\begin{equation*}
  FX^{\ct'}_{\Lambda'} := \Xect_{\Gamma'} \times \prod_{v \in T} \Xecto_{0, n(v), \beta(v)}.
\end{equation*}
Note that the perfect obstruction theory of $FX^{\ct}_{\Lambda'}$
relative to $D_{\Lambda'}$ is defined by $(R\pi_* \mathcal E')^\vee$,
for
\begin{equation*}
  \mathcal E' := \mathcal E'_x \oplus \mathcal E'_p := \bigoplus_{i = 1}^M (\L^{\otimes w_i} \otimes \O(-\Delta - {\textstyle\sum_{k = 1}^n} \Delta_k))|_{\cC'_g} \oplus \bigoplus_{j = 1}^N (\L^{\otimes -d} \otimes \omega_{\pi,\log}),
\end{equation*}
where $\Delta \subset \Delta'$ are the nodes at $\cC'_g$.
By the normalization exact sequence, we may write $R\pi_* \mathcal E'$
as the cone of
\begin{equation*}
  R\pi_* \mathcal E'_x|_{\cC'_g} \oplus R\pi_* \mathcal E'_p|_{\cC'_g} \oplus \bigoplus_{v \in T} R\pi_* \mathcal E'_p|_{\cC'_v} \to \mathcal E'_p|_{\Delta'}.
\end{equation*}
Define $D_{\Lambda''}$ in the same way as $D_{\Lambda'}$, but
instead of one line bundle $L$ on all of $C$, take one line bundle
$L_g$ on $\cC'_g$ and one line bundle $L_v$ on each $\cC_v$ for
$v \in T$.
We can then use the cone of
\begin{equation*}
  R\pi_* \mathcal E'_x|_{\cC'_g} \oplus R\pi_* \mathcal E'_p|_{\cC'_g} \oplus \bigoplus_{v \in T} R\pi_* \mathcal E'_p|_{\cC'_v} \to \bigoplus_{v \in T} T_{\ri X^\ct},
\end{equation*}
which is induced by the map $\mathcal E'_p|_q \to T_{\ri X^\ct}$ coming from the Euler sequence at each node $q \in \Delta'$, to
define a perfect obstruction theory of $FX^{\ct}_{\Lambda'}$ relative to
$D_{\Lambda''}$.
Note that the two relative perfect obstruction theories of
$FX^{\ct}_{\Lambda'}$ yield the same absolute perfect obstruction
theory.

For $FX^{\ct'}_{\Lambda'}$, we can define a perfect obstruction theory
relative to $D_{\Lambda''}$ via
\begin{equation*}
  \bigoplus_{i = 1}^M (\L_g^{\otimes w_i} \otimes \O(-\Delta - {\textstyle\sum_{k = 1}^n} \Delta_k)) \oplus \bigoplus_{j = 1}^N \left((\L_g^{\otimes -d} \otimes \omega_{\pi,\log}) \oplus \bigoplus_{v \in T}  (\L_v^{\otimes -d} \otimes \omega_{\pi,\log})\right),
\end{equation*}
where $\L_g$ and $\L_v$ are the universal line bundles on $\cC'_g$ and
$\cC_v$.
Then, the relative perfect obstruction theories of $FX^{\ct}_{\Lambda'}$
and $FX^{\ct'}_{\Lambda'}$ relative to $D_{\Lambda''}$ are compatible
with respect to the morphism $\widetilde\Delta$.
On the direct product $FX^{\ct'}_{\Lambda'}$, we may define a
cosection $\sigma''$ as in Section~\ref{subsec:vcycle}, but working
only on the component $\cC_g$.
This cosection pulls back to the cosection $\sigma'$ of
$FX^{\ct}_{\Lambda'}$.
Therefore, by cosection-localized pullback, we have
\begin{equation}
  \label{eq:split}
  [FX^{\ct}_{\Lambda'}]^\vir_{\sigma'} = \widetilde\Delta^! \left([FX^{\ct'}_{\Lambda'}]^\vir_{\sigma''}\right).
\end{equation}

Combining \eqref{eq:contract} and \eqref{eq:split}, we have shown that
\[\iota_*[FZ_{\Lambda}]^{\vir} = \widetilde{p}^*\widetilde\Delta^! \left([FX^{\ct}_{\Lambda'}]^{\vir}_{\sigma''}\right) \cap \prod_{v \in T} e(R^1\pi_* \mathcal E''_v).\]
Since $FX^{\ct}_{\Lambda'}$ is defined as a product, $[FX^{\ct'}_{\Lambda'}]^\vir_{\sigma''}$ clearly splits as a product of virtual cycles for each factor.  The factor of $X^{\epsilon,\ct}_{\Gamma'}$ yields a virtual cycle $[Z^{\epsilon}_{\Gamma'}]^{\vir}$, while the factor corresponding to each $v \in T$, after combining with $e(R^1 \pi_* \mathcal E''_v)$, yields
$[\Xecto_{0, n(v), \beta(v)}]^\vir$.
Thus, we recover the statement of Lemma~\ref{lem:vclass}.

\subsubsection{Comparison}
\label{subsubsec:comparison}

We are finally ready to express the total contribution of each localization graph $\Lambda$ in the form needed for what follows.  Let the dual graph $\Gamma'$ associated to $\Lambda$ be defined as in Subsection~\ref{subsubsec:virtual}, and let
\begin{equation*}
  p\colon FX^\ct_\Lambda \to \Xect_{\Gamma'}
\end{equation*}
be the map that forgets $L_2$ and $\vec{z}$ and contracts unstable
components (while preserving markings at the stable vertices to keep
track of where trees of rational components were attached).
Then, by Lemma~\ref{lem:vclass} and the virtual normal bundle
computations of the previous subsections, we can write
\begin{equation*}
  p_* \left(\frac{[FZ_\Lambda]^\vir}{e(N_\Lambda^\vir)}\right)
  = \deg(p) \cdot [Z^{\epsilon}_{\Gamma'}]_\tw^\vir \cdot \prod_{h \in H(\Gamma')} \ev_h^*(F_h(\psi_h)),
\end{equation*}
where $H(\Gamma')$ is the set of half-edges of $\Gamma'$.
The factors $F_h(z)$ are valued in $\H(\lambda)[z]$ and, in the case
that $h$ is a leg, are obtained by collecting localization factors and
performing integrals over moduli spaces
$Z^\epsilon_{0, n(v), \beta(v)}$ corresponding to $h$.
If $h$ belongs to an edge $e = \{h, h'\}$, we choose any way of
distributing localization factors corresponding to $e$ among $F_h$ and
$F_{h'}$.
In the case where $\beta(v) > 0$, the integrals over
$Z^\epsilon_{0, n(v), \beta(v)}$ are against the virtual cycle
$[\Xecto_{0, n(v), \beta(v)}]_\sigma^\vir$, where
$\Xecto_{0, n(v), \beta(v)}$ is defined as in the previous subsection. 
The degree $\deg(p)$ consists of root stack automorphism factors
$(\delta(e))^{-1}$ for each edge of $\Lambda$, as well as up to two
factors of $d_{m(h)}^{-1}$ for each edge of $\Lambda$ due to the fact
that the edge moduli space uses the ordinary inertia stack while the
fiber products are relative to the rigidified inertia stack.
More precisely, all edges in $\Lambda$ that have a basepoint give only
one factor of $d_{m(h)}^{-1}$ while any other edge gives
$d_{m(h)}^{-2}$.

Note that only a few factors appearing in the localization formula
depend on multiplicities:
\begin{itemize}
\item the prefactors of \eqref{eq:virloc}, which give a factor of
  $d_{m(h)}$ for each edge in the dual graph $\Gamma$
\item the factors $d_{m(h)}$ from Section~\ref{subsubsec:nodes} for
  each edge of $\Gamma$
\item the factors $d_{m(h)}^{-1}$ from Section~\ref{subsubsec:edges}
  for each edge of $\Lambda$ with a base point
\end{itemize}
The second and third factors are included in the $F_h$.
Let us define a modified series $F'_h(\psi_h)$ is obtained from
$F_h(\psi_h)$ by removing the second and third factors, while adding
in the factors $(\delta(e))^{-1}$ from $\deg(p)$.
Then, we can write
\begin{equation*}
 \prod_{h \in \widetilde H(\Lambda)} d_{m(h)} \cdot p_* \left(\frac{[FZ_\Lambda]^\vir}{e(N_\Lambda^\vir)}\right)
  = [Z^{\epsilon}_{\Gamma'}]_\tw^\vir \prod_{e \in E(\Gamma')} d_{m(e)}^2 \prod_{h \in H(\Gamma')} \ev_h^*(F'_h(\psi_h)).
\end{equation*}

The key point in what follows is that the insertion $F'_h$
corresponding to a tree of unstable rational components with
$\beta$-degree zero does not depend on the GIT quotient $X$ if one
regards the hyperplane class $H$ as a formal variable.
Therefore, $F'_h$ equals the corresponding factor in the case where
the GIT quotient $X$ is a single point, under the substitution
$\lambda \rightarrow \lambda_0$.
We consider this setting further in the next subsection.

\begin{remark}
\label{rem:geomloc}
Our calculations of the vertex, edge, and node contributions to the virtual normal bundle carry over directly to the geometric chamber and can be found in \cite{CJR}.  The decomposition of the virtual cycle carried out in Lemma \ref{lem:vclass} also holds in the geometric chamber; this can be obtained by the same proof as above, or, if one assumes the equivalence of the cosection-localized virtual cycle with Ciocan-Fontanine--Kim--Maulik's virtual cycle on $\M_{g,n}^{\epsilon}(Y,\beta)$, it follows from the known decomposition along nodes.
Thus, as we saw in \cite{CJR}, it is also the case that the insertion corresponding to a tree of rational components with $\beta$-degree zero in the geometric chamber equals the corresponding insertion for $X = \{\text{point}\}$ under the substitution $\lambda \rightarrow \lambda_0$, which is the crucial fact needed for what follows.
\end{remark}

\subsection{Equivariant orbifold projective line}
\label{subsec:eqP1}

If the GIT quotient $X$ is replaced by a single point (so that there
are no sections $\vec{x}, \vec{p}$), then the twisted graph space
reduces to the usual moduli space of stable maps to $\P^1$. 
In this section, we summarize explicit computations from
\cite[Section~4]{CJR} of related generating series that play a role in
the twisted graph space localization to follow.

As in Section~\ref{fixedloci}, the $\C^*$-fixed loci in
$\M_{g,n}(\P^1, \delta)$ can be indexed by $n$-legged graphs $\Gamma$,
where each vertex $v$ is decorated by an index
$j(v) \in \{0, \infty\}$ and a genus $g(v)$, and each edge $e$ is
decorated by a degree $\delta(e) \in \mathbb{N}$.
Each vertex $v$ corresponds to a maximal sub-curve of genus $g(v)$
contracted to the single point $j(v) \in \P^1$, or, in the unstable
case where the vertex has genus zero and valence one or two, to a
single point in the source curve.
Each edge $e$ corresponds to a noncontracted component, which is
necessarily of genus zero, and on which the map to $\P^1$ is of the
form $[x:y] \mapsto [x^{\delta(e)}: y^{\delta(e)}]$ in coordinates.

Fix insertions $\alpha_1, \ldots, \alpha_n \in H^*_{\C^*}(\P^1)$, and
let $p\colon \M_{g,n}(\P^1, \delta) \to \M_{g, n}$ be the forgetful map.
Then the localization formula expresses the class
\begin{equation*}
  p_*\left(\prod_{i=1}^n \ev_i^*(\alpha_i) \cap [\M_{g,n}(\P^1,\delta)]^{\vir}\right)
\end{equation*}
as a sum over contributions from each fixed-point graph $\Lambda$.
These expressions can be stated more efficiently by considering the generating series
\begin{equation}
  \label{int}
  \sum_{\delta = 0}^\infty y^{\delta} p_*\left(\prod_{i=1}^n \ev_i^*(\alpha_i) \cap [\M_{g,n}(\P^1,\delta)]^{\vir}\right)
\end{equation}
for a Novikov variable $y$.

Let $\Phi$ denote the sum of all contributions to \eqref{int} from
graphs $\Lambda$ on which there is a vertex $v$ with $g(v) = g$ and
$j(v) = 0$, and such that, after stabilization, the generic curve in
the moduli space corresponding to $\Lambda$ is smooth; the second
condition means that there is no tree emanating from $v$ that
contains more than one marking.
Therefore, emanating from the vertex $v$ on such a graph, there are
$n$ (possibly empty) trees on which at least one marking lies and $l$
trees with no marking, for some integer $l$.
It follows that
\begin{equation}
  \label{Phim}
  \Phi = \sum_{l=0}^{\infty} \frac 1{l!} \pi_{l*}\left(\sum_{i = 0}^g c_i(\mathbb E) \lambda^{g - 1 - i} \prod_{k = 1}^n S(\alpha_k, \psi_k) \prod_{k = n+1}^{n+l} \epsilon(\psi_k)\right),
\end{equation}
in which $\pi_l\colon \M_{g, n + l} \to \M_{g, n}$ denotes the forgetful
map, and $\mathbb E$ is the Hodge bundle.
The series $S(\alpha, z)$ in \eqref{Phim} is the universal generating
series of localization contributions of a tree emanating from a
vertex $v$ with $j(v) = 0$ that contains exactly one of the markings
and has an insertion of $\alpha \in H_{\C^*}^*(\P^1)$.  The series
$\epsilon(z)$, similarly, is the generating series of localization
contributions of a tree containing none of the markings.
Clearly, such a tree needs to have degree at least one, so
$\epsilon(z)$ is a multiple of $y$.

Let $\overline\psi_k$ be the pullback under $\pi_l$ of the class
$\psi_k$ on $\M_{g, n}$.
It is well-known that $\overline\psi_k$ differs from $\psi_k$ exactly
on the boundary divisors of $\M_{g, n}$ where the $k$th marking and
some of the last $l$ markings lie on a rational tail.
By rewriting the classes $\psi_1, \dotsc, \psi_n$ in terms of
$\overline\psi_1, \dotsc, \overline\psi_n$ and boundary divisors, and
for each summand integrating along the fibers of the map forgetting
all markings of the involved boundary divisors, we can rewrite
$\Phi$ in the form
\begin{equation*}
  \Phi = \sum_{l=0}^{\infty} \frac 1{l!} \pi_{l*}\left(\sum_{i = 0}^g c_i(\mathbb E) \lambda^{g - 1 - i} \prod_{k = 1}^n \widetilde S(\alpha_k, \overline\psi_k) \prod_{k = n+1}^{n+l} \epsilon(\psi_k)\right)
\end{equation*}
for modified universal series $\widetilde S(\alpha, z)$.
Surprisingly, the series $\widetilde S(\alpha, z)$ is easier to
compute than $S(\alpha, z)$.
In fact, it is closely related to the $R$-matrix for the equivariant
Gromov--Witten theory of $\P^1$.

For the proof of Theorem~\ref{thm:main}, what is relevant from this
formula is $\widetilde S(\alpha, 0)$, and we use the following result from \cite{CJR}.
(Note that $\widetilde S$ is called $\widetilde S_0$ in \cite{CJR}.)
\begin{lemma}
  \label{lem:P1}
  We have the identities
  \begin{align*}
    \widetilde S(\one_{\P^1}, 0) &= \phi^{-1/4}, \\
    \widetilde S(H, 0) &= \phi^{-1/4}\left(\frac\lambda 2 + \frac\lambda 2 \sqrt\phi\right),
  \end{align*}
  in which
  \begin{equation*}
    \phi := 1 + \frac{4y}{\lambda^2}.
  \end{equation*}
\end{lemma}

\section{Proof of Theorem~\ref{thm:main}}
\label{sec:proof}

We now turn to the proof of Theorem~\ref{thm:graphWC}, which, as
explained above, implies Theorem~\ref{thm:main}.

First, we introduce a bracket notation in order to write
Theorem~\ref{thm:graphWC} in a more economical way:
\begin{multline*}
  \left(\prod_{h \in H} \ev_h^*(\phi_h \psi^{a_h}) \right)^\WC_{\Gamma} :=
  \prod_{h \in H} \ev_h^*(\phi_h \psi^{a_h}) \cap [Z^\epsilon_{\Gamma}]^\vir_\tw \\
  - \sum_{K = \{k_v\}_{v \in V(\Gamma)}} \hspace{-0.1cm} \frac{1}{\prod_v k_v!}
 \hspace{-0.1cm} \sum_{\substack{\vec{\beta} = \{\beta_i^{v}\}_{v \in V(\Gamma), i \in [k_v]}\\ \beta_1^v + \cdots + \beta_{k_v}^v = \beta(v)}}\hspace{-1cm} b_{\vec{\beta} *} c_{*} \left( \prod_{h \in H} \ev_h^*(\phi_h \psi^{a_h}) \prod_{v \in V(\Gamma)} \prod_{j=1}^{k_v}  \ev_{v,j}^*(\mu_{\beta_j^v}^{\epsilon, \tw}(-\psi)) \cap [Z^{\infty}_{\Gamma_K + l}]^{\vir}_\tw\right).
\end{multline*}
Here, $H$ is any subset of the half-edges of $\Gamma_{\n}$, and
otherwise the notation is as in Theorem~\ref{thm:graphWC}.
One can extend the definition of $(\cdots)^{\WC}_{\Gamma}$ by
linearity to allow any insertions in $\H[[\psi]]$.

With this notation, the statement of Theorem~\ref{thm:graphWC} can be
written as
\begin{equation}
  \label{eq:graphWC}
  \sum_{l = 0}^\infty \frac 1{l!} \pi_{l*} \left(\prod_{h \in L(\Gamma)} \ev_h^*(\alpha_h) \prod_{h' \in L'} \ev_{h'}^*(\epsilon^{\lambda_0}(\psi)\one) \right)^\WC_{\Gamma + l} = 0,
\end{equation}
where $L(\Gamma)$ is the leg set of $\Gamma$ and $L'$ stands for the set of the $l$ additional legs.  It suffices to prove the theorem without any insertions $\alpha_h$, with the understanding that the only state-space elements that can be pulled back under our evaluation maps and paired with the virtual class are those of compact type.

\begin{remark}
The proof is essentially formal once one assumes wall-crossing results in genus zero and the localization calculations of Sections \ref{subsec:contributions} and \ref{subsec:eqP1}.  In particular, by \cite{CFKZero} and Remark~\ref{rem:geomloc}, the following proof applies without modification to the geometric chamber, yielding the other half of Theorem~\ref{thm:1}.
\end{remark}

\begin{proof}[Proof of Theorem~\ref{thm:graphWC}]
First, apply induction on the degree $\beta$.  When $\beta =0$, the result is trivially true, since $\epsilon$-stable and $\infty$-stable quasimaps agree in degree zero.

Assume, then, that $\beta >0$ and that the theorem holds for all
$\beta'<\beta$, and apply induction on the decorated graph
$(\Gamma, v_{\bullet}, \n)$.
The partial order on the set of such data used to organize the
induction is given by setting
\begin{equation*}
  (\Gamma', v_{\bullet}', \n') \leq (\Gamma, v_{\bullet}, \n)
\end{equation*}
if $\Gamma$ is obtained from $\Gamma'$ by replacing a subgraph
$\Gamma'' \subset \Gamma'$ containing $v_{\bullet}'$ by the single
vertex $v_{\bullet}$, and $\n'(v) = \n(v)$ for all
$v \in \Gamma \setminus \{v_{\bullet}\}$.
This partial order is illustrated in Figure \ref{fig:2}.
In order to be able to apply induction, it is important to note that
there is no infinite chain with respect to this partial order.

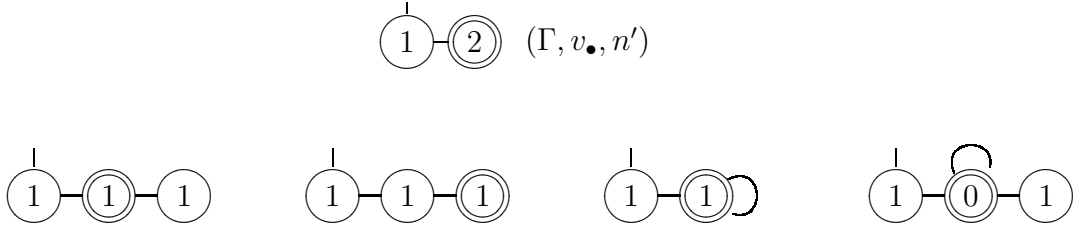
\begin{figure}
\begin{subfigure}{\textwidth}
\hspace{6cm}
\xymatrix@C-=0.2cm@R-=0.2cm{
&\\
*++[o][F-]{1}\ar@{-}[u]\ar@{-}[r] & *++[o][F=]{2}&{(\Gamma, v_{\bullet}, \n)}
}
\end{subfigure}

\vspace{0.75cm}

\begin{subfigure}{0.25\textwidth}
\xymatrix@C-=0.3cm@R-=0.3cm{
&&\\
*++[o][F-]{1}\ar@{-}[u]\ar@{-}[r] & *++[o][F=]{1}\ar@{-}[r] &*++[o][F-]{1}
}
\end{subfigure}
\hspace{1cm}
\begin{subfigure}{0.25\textwidth}
\xymatrix@C-=0.3cm@R-=0.3cm{
&&\\
*++[o][F-]{1}\ar@{-}[u]\ar@{-}[r] & *++[o][F-]{1}\ar@{-}[r] &*++[o][F=]{1}
}
\end{subfigure}
\hspace{1cm}
\begin{subfigure}{0.25\textwidth}
\xymatrix@C-=0.3cm@R-=0.3cm{
&&\\
*++[o][F-]{1}\ar@{-}[u]\ar@{-}[r] & *++[o][F=]{1}\ar@{-}@(ur,dr)& 
}
\end{subfigure}
\hspace{1cm}
\begin{subfigure}{0.25\textwidth}
\xymatrix@C-=0.3cm@R-=0.3cm{
&&\\
*++[o][F-]{1}\ar@{-}[u]\ar@{-}[r] & *++[o][F=]{0}\ar@{-}[r]\ar@{-}@(ul,ur) & *++[o][F-]{1}
}
\end{subfigure}
\caption{An example of a graph $(\Gamma, v_{\bullet},\n)$ and some of the graphs below it in the partial order.  The vertex $v_{\bullet}$ (or, in the second line, $v_{\bullet}'$) is marked with a double circle, and the numbers in each circle indicate the genus.  The data of $\n$ (or, in the second line, $\n'$) is not shown, but can be added arbitrarily as long as the left-most vertex agrees between the first and second lines.}
\label{fig:2}
\end{figure}

Now, suppose that the theorem holds for all
$(\Gamma', v_{\bullet}', \n') < (\Gamma, v_{\bullet}, \n)$.  For the following discussion, we make a choice $N_\infty^\bullet$ of
the set of half-edges at $v_\bullet$, but we usually suppress this
choice from the notation.
Let $N_0^\bullet$ be the other half-edges at $v_\bullet$.
For $i \in \{0, \infty\}$, let $N_i'$ be the set of all half-edges of
$\Gamma$ at vertices in $V_i \setminus \{v_\bullet\}$, and set $N_i = N_i^\bullet \cup N_i'$.

For each integer $\delta \geq 0$, define a new graph $\Gamma_{\delta}$ by, for
each edge $e$ that consists of a half-edge in $N_0$ and a half-edge
in $N_\infty$, subdividing $e$ with a single new vertex $v$ and
setting $g(v) = \beta(v) = 0$.
(This creates new half-edges, but the multiplicities on these new
half-edges are determined.) 
Extend the coloring on vertices of $\Gamma$ to vertices of $\Gamma_{\delta}$
by setting
\begin{equation*}
  V(\Gamma_{\delta}) = V_0 \cup V_{\infty} \cup V_{\text{new}},
\end{equation*}
where $V_{\text{new}}$ denotes the vertices of $\Gamma_{\delta}$ not coming
from vertices of $\Gamma$.
If we define
\begin{equation*}
  D\colon V(\Gamma_{\delta}) \rightarrow \mathbb{N}
\end{equation*}
by
\begin{equation*}
  D(v) =
  \begin{cases}
    \delta & \text{ if } v = v_{\bullet}\\
    1 & \text{ if } v \in V_{\text{new}}\\
    0 & \text{ otherwise},
  \end{cases}
\end{equation*}
then $\Gamma_{\delta}$ defines a moduli space
$\P \Ze_{\Gamma_{\delta}}$ with a finite map to a substack of
$\P\Ze_{g,n+n',\beta,\delta}$.

The proof of the theorem proceeds by computing, via localization, the
difference between the expressions
\begin{equation}
  \label{eq:epsilon}
  \sum_{\delta=0}^{\infty} y^{\delta} p_*\left(\prod_{n \in N_0} \ev_n^*(\one \otimes \varphi_0) \prod_{n \in N_{\infty}} \ev_n^*(\one \otimes \varphi_\infty) \cap [\P \Ze_{\Gamma_{\delta}}]^{\vir} \right)
\end{equation}
and
\begin{multline}
  \label{eq:infinity}
  \sum_{\delta=0}^{\infty} \sum_{K = \{k_v\}_{v \in V}} \frac{1}{\prod_v k_v!}\sum_{\substack{ \vec{\beta}_v = \{\beta_j^v\}_{j \in [k_v]}\\ \beta_1^v + \cdots + \beta_{k_v}^v = \beta(v)}} \sum_{\substack{\vec\delta_v = \{\delta_j^v\}_{j \in [k_v]} \\ \delta_j^v \le \beta_j^v \\ \delta_j^v = 0 \text{ if }v \neq v_\bullet}} y^{\delta + \delta_1^{v_\bullet} + \dotsb + \delta_{k_{v_\bullet}}^{v_\bullet}}\\
  p_* \widetilde b_{\vec{\beta}_v*} \widetilde c_{v*}\left(\prod_{n \in N_0} \ev_n^*(\one \otimes \varphi_0) \prod_{n \in N_{\infty}} \ev_n^*(\one \otimes \varphi_\infty) \prod_{v \in V} \prod_{j=1}^{k_v} \ev_{v,j}^*(\widetilde{\mu}_{\beta_j^v, \delta_j^v}^{\epsilon}(-\psi))  \cap [\P Z^{\infty}_{\Gamma_{\delta,K}}]^{\vir} \right).
\end{multline}
Here,
\begin{equation*}
  p\colon \P\Ze_{\Gamma_{\delta}} \rightarrow \Ze_{\Gamma}
\end{equation*}
is the morphism forgetting $L_2$, $z_1$, and $z_2$ (and stabilizing if
necessary), and
\begin{equation*}
  \varphi_0 = \frac{[0]}\lambda, \qquad \varphi_\infty = -\frac{[\infty]}\lambda.
\end{equation*}
The morphisms $\tilde b_{\vec{\beta}_v}$ and $\tilde{c_v}$ are lifts
of $b_{\vec{\beta}_v}$ and $c_v$ to the twisted graph space, and
$\Gamma_{\delta,K}$ is obtained from $\Gamma_{\delta}$ by setting the
$\beta$-degree of all vertices to zero and adding $k_v$ additional
legs of multiplicity $\frac 1d$ to each vertex $v$.
The mirror transformation $\widetilde\mu^\epsilon$ of the twisted graph space is defined by
\begin{equation*}
  \sum_{\delta = 0}^\infty y^\delta \widetilde\mu_{\beta, \delta}^\epsilon(z)
  = \mu^{\epsilon, \tw}_{\beta}(z) \otimes \varphi_0 - \sum_{\delta = 1}^\beta y^\delta \frac{\left[J^\epsilon_\tw\left(q, \frac{\lambda_0}\delta\right)\right]_{q^\beta} \lambda_\infty}{\left(\frac{\lambda_\infty}\delta + z\right) \prod_{i = 1}^\delta \prod_{j \in \{0, \infty\}} \frac{i\lambda_j}\delta} \otimes \varphi_\infty.
\end{equation*}

The definition of $\Gamma_{\delta}$ is specifically chosen such that only very
special localization graphs $\Lambda'$ contribute to
\eqref{eq:epsilon} and \eqref{eq:infinity}.
In particular, each contributing localization graph $\Lambda'$ for \eqref{eq:epsilon} is
obtained from $\Gamma_{\delta}$ by
\begin{itemize}
\item replacing each $v \in V_{\text{new}}$ by an edge $e$ with $\delta(e)=1$;
\item keeping each vertex $v \in (V_0 \cup V_{\infty}) \setminus \{v_{\bullet}\}$ the same;
\item replacing $v_{\bullet}$ by a localization graph
  $\Lambda_{\delta}^\bullet$ of genus $g(v_{\bullet})$, total $\beta$-degree
  $\beta(v_{\bullet})$, and total edge degree $\delta$.
\end{itemize}
The description of possible localization graphs for
\eqref{eq:infinity} is similar, except that the localization graph
$\Lambda_{\delta}^\bullet$ has $\beta$-degree zero and additional
markings depending on the choice of $K$.
Note that localization graphs $\Gamma'$ for \eqref{eq:epsilon} need to
satisfy $\delta(e) > \beta(v)$ for each genus-zero, valence-one vertex
$v$ with unique incident edge $e$ (see Remark~\ref{rem:deltabeta}).
If a dual graph does not necessarily satisfy this degree condition,
but it satisfies all other conditions of a dual graph contributing to
\eqref{eq:epsilon}, then we refer to it as a ``fake localization
graph''.

For each tuple $K$, each choice of a contributing localization
graph $\Lambda'_K$ to \eqref{eq:infinity}, and each choice of
$\vec\beta_{v_{\bullet}}$ and $\vec\delta_{v_{\bullet}}$, we define
a fake localization graph $\Lambda'$ for \eqref{eq:epsilon} by applying the following operations to $\Lambda'_K$:
\begin{itemize}
\item Remove all extra legs at vertices $v$ with $j(v) = 0$ and add
  their $\beta$-degree to the incident vertex.
\item Replace each extra leg at a vertex $v$ by $j(v) = \infty$
  with a new edge connected to a new genus-zero, valence-one vertex.
  The degree of the new edge is prescribed by
  $\vec\delta_{v_\bullet}$, and the $\beta$-degree of the former extra leg is put onto the new
  vertex.
\end{itemize}
We refer to the sum of all contributions of localization graphs
$\Lambda'_K$ corresponding to the fake localization graph $\Lambda'$
as the ``contribution of $\Lambda'$ to \eqref{eq:infinity}''.

There are several summands in the contribution of $\Lambda'$ to
\eqref{eq:infinity}.
For each vertex $v$ of genus zero and valence one with unique incident
edge $e$ such that $\delta(e) > \beta(v)$, the local contribution to
the localization formula is given by
\begin{multline*}
  \delta_{\beta(v), 0} \frac{\lambda_0}{\delta(e)} + \mu_{\beta(v)}^{\epsilon, \tw}\left(\frac{\lambda_0}{\delta(e))}\right) \\
  + \sum_{\substack{k \ge 1, \beta_0 + \dotsb + \beta_k = \beta(v) \\ k + \beta_0 > 1}} \ev_{1, *} \left(\frac{\prod\limits_{i = 1}^k \ev_{i + 1}^* \left(\mu_{\beta_i}^{\epsilon, \tw}(-\psi_{i + 1})\right)}{\left(\frac{\lambda_0}{\delta(e)} - \psi_1\right) k!} \cap [Z^{\infty}_{0,1+k,\beta_0}]^{\vir}_\tw\right) \\
  = \left[J^{\infty,\tw}\left(q, \sum_\beta q^\beta \mu_\beta^{\epsilon, \tw}\left(-\frac{\lambda_0}{\delta(e)}\right), \frac{\lambda_0}{\delta(e)}\right)\right]_{q^{\beta(v)}},
\end{multline*}
where the first two summands correspond to the cases where $v$ is
unstable and the last summand corresponds to the case where $v$ is
stable.
By the $J$-function wall-crossing \eqref{eq:JWC}, this contribution
is equal to
\begin{equation*}
  \left[J^{\epsilon,\tw}\left(q, \frac{\lambda_0}{\delta(e)}\right)\right]_{q^{\beta(v)}}.
\end{equation*}
The definition of $\widetilde\mu^\epsilon$ is exactly chosen such that
the local contribution to the localization formula of a vertex $v$ of
genus zero and valence one with unique incident edge vanishes if $\delta(e) \leq \beta(v)$.  This means that the contribution to \eqref{eq:infinity} of any fake
localization graph $\Lambda'$ that is not an actual localization
graph contributing to \eqref{eq:epsilon} vanishes.
Because of that, from now on we only consider the contribution of
actual localization graphs to the difference of \eqref{eq:epsilon}
and \eqref{eq:infinity}.

Fix a localization graph $\Lambda'$ contributing to
\eqref{eq:epsilon}, and let $\Gamma'$ denote the graph from
Section~\ref{subsubsec:comparison}.
If $\Lambda'$ has a genus-zero component of positive $\beta$-degree
that is be contracted under the forgetful map, then the total
$\beta$-degree of $\Gamma'$ is strictly less than $\beta$.
It follows that the contribution of $\Lambda'$ to \eqref{eq:epsilon}
equals its contribution to \eqref{eq:infinity}, by the inductive
hypothesis together with an application of the genus-zero
wall-crossing theorem \cite{ClaRo2} and the above application of the $J$-function
wall-crossing.
Thus, for the study of the difference of \eqref{eq:epsilon} and
\eqref{eq:infinity}, we can from now on assume that the localization
graphs have no genus-zero, valence-one vertex of positive
$\beta$-degree.

In order to better package the localization formula for the remaining
graphs, we form a generating series depending on a variable $y$ that
keeps track of the degree $\delta$.
Furthermore, we rearrange the localization graphs $\Lambda'$ according
to the decorated dual graph $(\Gamma', v_{\bullet}', \n')$ associated
to $\Lambda'$, in which $\Gamma'$ is again the stabilization of
$\Lambda'$, $v'_{\bullet}$ is an arbitrary new vertex\footnote{Because
  $\beta > 0$ and using the above assumption on $\Gamma'$, such a
  vertex exists.
  We make a uniform choice for all $\Lambda'$ corresponding to the
  same $\Gamma'$.} with $j(v_\bullet) = 0$, and $\n'(v)$ is the number
of unmarked trees of rational curves in $\Lambda'$ emanating from the
vertex $v \neq v_\bullet$. 
The coloring $V(\Gamma') = V_0' \cup V_{\infty}'$ of the vertices of
$\Gamma'$ is dictated by the labels $j(v)$ on the stable vertices of
$\Lambda'$.

It is important to notice that $\Gamma'$ is obtained from $\Gamma$ by
(possibly) degenerating the vertex $v_{\bullet}$ but leaving all other
vertices the same, and $\Lambda'$ is obtained from $\Gamma'$ by
attaching trees of the following types:
\begin{enumerate}
\item those that connect two vertices of $\Gamma'$ that come from the
  degeneration of $v_{\bullet}$;
\item those that are attached to a single vertex of $\Gamma'$ and contain exactly
  one leg of $v_\bullet$;
\item those that are attached to a single vertex of $\Gamma'$ and
  contain no leg of $v_\bullet$.
\end{enumerate}

First, note that the local contribution of any such tree to the
localization formulas for \eqref{eq:epsilon} or \eqref{eq:infinity} are
evidently identical.
Then, let $S^{\lambda_0}(\varphi_j, z)$ be the generating series
keeping track of the contribution to the localization formula for
\eqref{eq:epsilon} (or \eqref{eq:infinity}) of all trees of the second
type attached to $v'_\bullet$ and with insertion $\varphi_j$ at the
leg, and let $\epsilon^{\lambda_0}(z)$ be the similar generating
series corresponding to trees of the third type attached to
$v'_\bullet$.
We note that $S^{\lambda_0}(\varphi_j, z)$ agrees with
$S(\varphi_j, z)$ from Section~\ref{subsec:eqP1} except for the
substitution $\lambda \to \lambda_0$, and the relation between
$\epsilon$ and $\epsilon^{\lambda_0}$ is the same.
It is not necessary to define generating series for the trees
that are not attached to $v'_\bullet$.

Consider a set of localization graphs $\Lambda'$ that correspond to a
fixed decorated dual graph $(\Gamma', v_{\bullet}', \n')$ and that
agree with each other except for trees of the third type at
$v'_\bullet$.
The contribution of such a set of localization graphs to difference of
\eqref{eq:epsilon} and \eqref{eq:infinity} is essentially
\begin{equation*}
  \sum_{l = 0}^\infty \frac 1{l!} \pi_{l*} \left(\prod_{h \in H} \ev_h^*(U_h(\psi)) \prod_{h \in L'} \ev_h^*(\epsilon^{\lambda_0}(\psi)) \right)^\WC_{\Gamma'+ l},
\end{equation*}
where $H$ and $L'$ to denote the set of half-edges at $v'_\bullet$ or
the $l$ additional half-edges, respectively, and $U_h(\psi)$ stands
for the localization contribution of the (possibly empty) tree at $h$.
The above formula differs from the correct formula by an invertible
factor from the localization contribution of trees that are not
attached to $v'_\bullet$.

We claim that Theorem~\ref{thm:graphWC} would imply that this
contribution is zero. 
Indeed, rewriting the insertions in terms of the $\psi$-classes
$\overline\psi$ pulled back under $\pi_l$ and using a slightly
generalized bracket notation, it becomes
\begin{equation}
  \label{eq:locpull}
  \sum_{l = 0}^\infty \frac 1{l!} \pi_{l*} \left(\prod_{h \in H} \ev_h^*(\widetilde U_h(\overline\psi)) \prod_{h \in L'} \ev_h^*(\epsilon^{\lambda_0}(\psi)) \right)^\WC_{\Gamma' + l},
\end{equation}
where the relation between $U_h$ and $\widetilde U_h$ is the same as
the one between $S$ and $\widetilde S$ discussed in
Section~\ref{subsec:eqP1}.
Adding insertions has the effect of multiplying \eqref{eq:graphWC} by
a cycle, so if Theorem~\ref{thm:graphWC} holds for
$(\Gamma', v_{\bullet}', \n')$, then the contribution of
$(\Gamma', v_{\bullet}', \n')$ to the localization must vanish.
In particular, if $v_{\bullet}$ degenerates at all in the passage from
$\Gamma$ to $\Gamma'$, then the contribution of
$(\Gamma', v_{\bullet}', \n')$ to the difference of \eqref{eq:epsilon}
and \eqref{eq:infinity} vanishes by the induction hypothesis.

Thus, all that remains is the contribution to the localization formula
from $\Lambda'$ for which the decorated stabilized graph is the same
as $(\Gamma, v_{\bullet}, \n)$ except for possible additional legs at
$v_\bullet$.
The contribution from these remaining graphs to the difference of
\eqref{eq:epsilon} and \eqref{eq:infinity} can be written as
\begin{equation*}
  \sum_{l = 0}^\infty \frac 1{l!} \pi_{l*} \left(\prod_{h \in N_0^\bullet} \ev_h^*(\widetilde S^{\lambda_0}(\varphi_0, \overline\psi)) \prod_{h \in N_\infty^\bullet} \ev_h^*(\widetilde S^{\lambda_0}(\varphi_\infty, \overline\psi)) \prod_{h \in L'} \ev_h^*(\epsilon^{\lambda_0}(\psi)) \right)^\WC_{\Gamma + l}.
\end{equation*}
We filter these remaining localization contributions (as well as the
statement of Theorem~\ref{thm:graphWC}) according to cohomological
degree.
This can be made explicit by separating the ``insertion''
$\exp(\epsilon^{\lambda_0})$ into its components of degree $i$:
\begin{equation*}
  \sum_{l = 0}^\infty \frac 1{l!} \epsilon^{\lambda_0}(z_1) \otimes \dotsb \otimes \epsilon^{\lambda_0}(z_l) = \sum_{i=-\infty}^{\infty} \epsilon_i,
\end{equation*}
where the degree of a homogeneous element
\begin{equation*}
  \phi_1 \otimes \cdots \otimes \phi_l \; \lambda^i y^j z_1^{a_1} \cdots z_l^{a_l} \in (\H^{\otimes l} \otimes_{\C} R)[[z_1, \ldots, z_l]]
\end{equation*}
is defined by
\begin{equation*}
  \text{codim}(\phi_1) + \cdots + \text{codim}(\phi_l) + 0 \cdot (i + j) + (a_1 + \cdots + a_l) - l.
\end{equation*}
We can further decompose
\begin{equation*}
  \epsilon_i = \sum_{l = 0}^\infty \epsilon_i^l,
\end{equation*}
such that
$\epsilon_i^l \in (\H^{\otimes l} \otimes_{\C} R)[[z_1,
\ldots, z_l]]$.

We now apply a further induction on $i$ to prove that the
$\epsilon_i$-part of the localization contribution, namely that
\begin{equation}
  \label{eq:locei}
  \sum_{l = 0}^\infty \frac 1{l!} \pi_{l*} \left(\prod_{h \in N_0^\bullet} \ev_h^*(\widetilde S^{\lambda_0}(\varphi_0, \overline\psi)) \prod_{h \in N_0^\bullet} \ev_h^*(\widetilde S^{\lambda_0}(\varphi_\infty, \overline\psi)) \cdot \ev_{H_l}^*(\epsilon_i^l(\psi_{H_l})) \right)^\WC_{\Gamma + l},
\end{equation}
and the $\epsilon_i$-part of Theorem~\ref{thm:graphWC}, that is
\begin{equation}
  \label{eq:thmei}
  \sum_{l = 0}^\infty \frac 1{l!} \pi_{l*} \left(\ev_{H_l}^*(\epsilon_i^l(\psi_{H_l})) \right)^\WC_{\Gamma + l},
\end{equation}
both vanish.
Here, abusing notation somewhat, we denote by $\ev_{H_l}$ the evaluation map at all of the $l$ additional markings and by $\psi_{H_l}$ the corresponding collection of $\psi$-classes.

To start the induction, notice that for $i\ll 0$, the
$\epsilon_i$-part of the localization contribution of $\Gamma$ has
dimension greater than the dimension of each component of
$Z^\epsilon_\Gamma$, so this contribution vanishes.
Suppose, then, that \eqref{eq:locei} vanishes whenever $i < i_0$.
Consider the part of the contribution of $\Gamma$ in dimension
\begin{equation}
  \label{eq:mdim}
  \vdim(Z_\Gamma^\epsilon) - i_0.
\end{equation}
The part in this dimension only involves $\epsilon_i$ for $i \le i_0$.
Moreover, all of its terms $\epsilon_i$ for $i < i_0$ are
obtained from \eqref{eq:thmei} by multiplication by a cohomology
class, and hence they vanish by the inductive hypothesis.
Thus, the part of the localization contribution of dimension
\eqref{eq:mdim} equals
\begin{equation*}
  (\widetilde S(\varphi_0, 0))^{|N_0^\bullet|} (\widetilde S(\varphi_\infty, 0))^{|N_\infty^\bullet|} \sum_{l = 0}^\infty \frac 1{l!} \pi_{l*} \left(\ev_{H_l}^*(\epsilon_i^l(\psi_{H_l})) \right)^\WC_{\Gamma + l}.
\end{equation*}

The key observation, then, is that by Lemma~\ref{lem:P1}, this contribution changes by a
factor of
\begin{equation}
  \label{eq:irr}
  \frac{1 - \sqrt\phi}{1 + \sqrt\phi}
\end{equation}
if we increase the order of $N_\infty^\bullet$ by one (thereby
decreasing the order of $N_0^\bullet$ by one).
This is possible because either $\Gamma$ is the trivial graph and
$v_\bullet$ has all of the legs (of which there are at least one), or
$\Gamma$ is a non-trivial graph and there is therefore at least one
edge at $v_\bullet$.

Since \eqref{eq:irr} has infinitely many negative powers when expanded
as a Laurent series in $\lambda$, this is only possible if the
contribution at hand is identically zero.  This completes the proof that the $\epsilon_{i_0}$-part of the
contribution of $\Gamma$ to the difference between \eqref{eq:epsilon}
and \eqref{eq:infinity} vanishes, and dividing by
\begin{equation*}
  (\widetilde S(\varphi_0, 0))^{|N_0^\bullet|} (\widetilde S(\varphi_\infty, 0))^{|N_\infty^\bullet|}
\end{equation*}
proves the $\epsilon_{i_0}$-part of Theorem~\ref{thm:graphWC}.

\end{proof}

\appendix

\renewcommand\thesection{A}

\section{Wall-crossing when $n=0$ (by Yang Zhou)}



In this appendix, we extend the wall-crossing in the Landau--Ginzburg phase to the case $n=0$, thus completing the proof of Theorem~\ref{thm:1}.  Thanks to Theorem~\ref{thm:main}, in order to express $[Z^{\epsilon}_{g,0,\beta} ]^{\mathrm{vir}}$ in terms of $\infty$-theory virtual cycles, it suffices to express it in terms of $[Z^{\epsilon}_{g,(1/d),\beta} ]^{\mathrm{vir}}$.  This can be done via an analogue of the dilaton equation when $g>1$, and via an analogue of the divisor equation when $N>1$.

The direct generalizations of the dilaton and divisor equations fail in this setting, due to the incompatibility of the virtual cycles under the forgetful map
\[
  Z^{\epsilon}_{g,(1/d),\beta} \to Z^{\epsilon}_{g,0,\beta}.
\]
However, this incompatibility can be remedied by replacing the marking with a
``light" marking, which can indeed be forgotten without affecting the virtual
cycle. Applying the techniques developed in the previous work \cite{Z}, the
virtual cycle with a light marking can then be related to
$[Z^{\epsilon}_{g,(1/d),\beta}]^{\mathrm{vir}}$ via another wall-crossing
formula.

\subsection{The light marking}
The first task is to define $\epsilon$-stable Landau--Ginzburg quasimaps with a light
marking of multiplicity $1/d$.  We assume that $\epsilon$ is not on a wall---i.e., $k\epsilon\neq 1$ for any
integer $k$---and that $2g-2+n+\epsilon\beta>0$.  We then pick a positive rational number $\delta$ small enough so that
\begin{equation}
  \label{eqn_small_delta}
  k\epsilon -1 \quad \text{and} \quad k\epsilon-1+\delta \text{~~have the same
    sign for any integer $k$.} \tag{*}
\end{equation}

\begin{definition}
\label{def:quasimap} An {\it $\epsilon$-stable Landau--Ginzburg quasimap} to $Z$
with a light marking of multiplicity $1/d$ and $n$ heavy markings consists of an $(n+1)$-pointed
prestable Landau--Ginzburg quasimap of genus $g$ to $Z$ 
\[
  (C; q_1,\ldots, q_{n+1};L;\vec p)
\]
satisfying the following conditions:
\begin{itemize}
\item {\it Light marking}: The isotropy group at $q_1$ is cyclic of order $d$, and
  the line bundle $L$ has multiplicity $1/d$ at $q_1$.
\item {\it Nondegeneracy}: The zero set of $\vec{p}$ is disjoint from
the markings $q_2, \ldots, q_{n+1}$ and the nodes of $C$, and for each zero $q$ of $\vec{p}$, the
order of the zero satisfies
  \begin{equation}
  \label{nondegen2} \text{ord}_q(\vec{p}) \leq \frac{1}{\epsilon}.
\end{equation}
Here, when $q=q_1$, we define the vanishing order of $p_i$ to be the multiplicity of the divisor $(p_i)$ on the coarse curve, or in other words, the integer $k$ such that
\[(p_i) = kd[q_1] + [\text{other points}].\]
In particular, by \eqref{eqn_small_delta}, condition \eqref{nondegen2} for $q=q_1$ is equivalent to $\delta + \epsilon\text{ord}_{q_1}(\vec{p}) < 1$.
\item {\it Stability}:  The $\Q$-line bundle
\begin{equation*}
\label{stab2}
(L^{\otimes -d} \otimes \omega_{\log})^{\otimes \epsilon} \otimes \omega_{\log}(d(\delta -1)\cdot [q_1])
\end{equation*}
is ample.  (Note, here, that the line bundle $\omega_{\log}(d(\delta -1)\cdot [q_1])$ is the pullback of 
$\omega_{|C|} (\delta \overline {q_1}+\overline {q_2} +\cdots+\overline {q_{n+1}})$, where $\overline{ q_i}$ is the image of $q_i$ in the coarse curve $|C|$.)
\end{itemize}
\end{definition}

The stability condition
above is equivalent to $q_1$ being a weight-$\delta$ (``light")
marking and $q_2 ,\ldots, q_{n+1}$ being weight-$1$ (``heavy") markings, in the sense of
\cite{CFKBigI}.  It is clear that the definition is independent of the choice of $\delta$. 

The theory of quasimaps with a light marking is parallel to the case with only heavy markings.
In particular, there is a proper Deligne--Mumford stack
$Z^{\epsilon,\delta}_{g,(1/d)+n,\beta}$ parameterizing genus-$g$,
$\epsilon$-stable Landau--Ginzburg quasimaps of degree $\beta$ to $Z$, with a
light marking of multiplicity $1/d$ and $n$ heavy markings, up to isomorphism.  The compact-type virtual cycle
$[Z^{\epsilon,\delta}_{g,(1/d)+n,\beta}]^{\mathrm{vir}}$ and evaluation
map $\ev_1$ are defined in the same way as before, with the only difference being that $\ev_1$ maps to
the rigidified inertia stack $\ri [\mathbb C^N/\mathbb
C^* ]$ for the action of $\C^*$ on $\C^N$ with weights $(d ,\ldots, d)$. Indeed, it lands in
the component indexed by $\bar 1\in \mathbb Z_d$, which is isomorphic to the quotient
$[\mathbb C^N/\mathbb C^*]$ by the weight-$1$ action.

There is a forgetful morphism
\[
  \tau:  Z^{\epsilon,\delta}_{g,(1/d)+n,\beta} \to Z^{\epsilon}_{g,n,\beta}
\]
forgetting the light marking, which only contracts components on which the degree of $L$ is zero.  Thus, $Z^{\epsilon,\delta}_{g,(1/d)+n,\beta}$ is identified with the universal curve 
of $Z^{\epsilon}_{g,n,\beta}$, and in particular, $\tau$ is flat. Moreover, the universal line
bundles of the two moduli spaces are compatible under $\tau$, and hence, so are the relative perfect obstruction theories.  Thus, we obtain:

\begin{lemma}
  \label{lem:forget-light}
   $\tau^*[Z^{\epsilon}_{g,n,\beta}]^{\mathrm{vir}} = [Z^{\epsilon,\delta}_{g,(1/d)+n,\beta}]^{\mathrm{vir}}$.
\end{lemma}

\subsection{The master space}

Our next goal is to compare the virtual cycle for the moduli space
$Z^{\epsilon,\delta}_{g,(1/d),\beta}$ with one light marking (and no heavy
markings) to the virtual cycle for the moduli space $Z^{\epsilon}_{g,
(1/d),\beta}$ with one heavy marking (and no light markings). The technique for
doing so is localization on a larger moduli space that we refer to as the
``master space''.

Assume that $2g-2+n+1+\epsilon\beta>0$ and $2g-2+n\geq -1$, and let $S$ be any scheme.  Then the $S$-points of the master space are as follows:

\begin{definition}
\label{def:quasimap} An $S$-family of  $\epsilon$-stable Landau--Ginzburg quasimaps to $Z$
{\it with a mixed marking} consists of
\[
  (\pi:C\to S; q_1,q_2 ,\ldots, q_{n+1};L,N; \vec p,v_1,v_2 ),
\]
where
\begin{enumerate}
\item 
    $(\pi:C\to S; \vec q;L;\vec p)$
  is an $S$-family of $(n+1)$-pointed prestable Landau--Ginzburg quasimaps of genus
  $g$ to $Z$;
\item $N$ is a line bundle on $S$;
\item $v_1\in H^0(S,T_{q_1}\otimes N)$ and $v_2\in H^0(S,N)$ are  sections
  without common zeros, where $T_{q_1}$ is the line bundle on $S$ formed by the
  relative tangent spaces to the {\it coarse} curves at $q_1$.
\end{enumerate}
We require that this data satisfy the following conditions on each geometric fiber:
\begin{itemize}
\item {\it Mixed marking:}
  The isotropy group at $q_1$ is cyclic of order $d$, and the line bundle $L$ has multiplicity
  $1/d$ at $q_1$.
\item {\it Nondegeneracy}:
  The zero set of $\vec{p}$ is disjoint from
  markings $q_2 ,\ldots, q_{n+1}$ and the nodes of $C$,
  and for each zero $q$ of $\vec{p}$, the
  order of the zero satisfies
\begin{equation*}
\text{ord}_q(\vec{p}) \leq \frac{1}{\epsilon}.
\end{equation*}
\item {\it Generic Stability}:  The $\Q$-line bundle
\begin{equation*}
(L^{\otimes -d} \otimes \omega_{\log})^{\otimes \epsilon} \otimes \omega_{\log}
\end{equation*}
is ample.
\item
  When $v_1=0$, $\vec p$ does not vanish at $q_1$. 
\item
  When $v_2=0$, the $\Q$-line bundle
\begin{equation*}
(L^{\otimes -d} \otimes \omega_{\log})^{\otimes \epsilon} \otimes \omega_{\log}(d(\delta-1)\cdot [q_1])
\end{equation*}
is ample.
\end{itemize}
\end{definition}

Thus, a family over a point consists of a prestable Landau--Ginzburg quasimap to
$Z$ together with $v_1/v_2 \in T_{q_1}C\cup \{\infty\}$, where $q_1$ behaves as a
heavy marking when $v_1=0$ and as a light marking
when $v_2=0$. We call $q_1$ a \textit{mixed marking} in what follows, whereas
$q_2 ,\ldots, q_{n+1}$ are heavy markings.
\begin{theorem}
There is a proper Deligne--Mumford stack $\widetilde
Z^{\epsilon,\delta}_{g,(1/d)+n,\beta}$ (the ``master space")
parameterizing genus-$g$, $\epsilon$-stable Landau--Ginzburg
quasimaps of degree $\beta$ to $Z$, with a mixed marking of multiplicity $1/d$
and $n$ heavy markings, up to isomorphism.
\end{theorem}

\begin{proof}
  The same argument as in the proof of \cite[Theorem 4]{Z} shows that the moduli problem is
  represented by a Deligne--Mumford stack of finite type over $\mathbb C$. We now
  use the valuative criterion to prove its properness. The proof
  is similar to that of \cite[Theorem 5]{Z}.

  Let $R$ be a Henselian discrete valuation ring with residue field $\mathbb C$. 
  Let $B=\operatorname{Spec} R$, with closed point $b\in B$ and generic point $B^\circ = B\backslash \{b\}$.  For a stable family
  \[
    \xi^\circ = (\pi^\circ:C^\circ\to B^\circ;\vec{q^{\circ}};L^\circ,N^\circ;\vec {p^\circ},v^\circ_1,v^\circ_2)
  \]
over $B^\circ$, we must show that (possibly after a finite
  base change) we can extend $\xi^\circ$ to a $B$-family
  \[
    \xi = (\pi:C\to B;\vec q;L,N;\vec p,v_1,v_2),
  \]
  and the extension is unique up to unique isomorphism.

  Let $|C^\circ|$ be the coarse moduli of $C^\circ$. The line bundle 
  \[
     (L^\circ)^{\otimes -d}\otimes \omega_{\pi^\circ,\mathrm{log}}
  \]
  descends to a line bundle $M^\circ$ on $|C^\circ|$, and the section $\vec {p^\circ}$
  descend to a section $\vec{{f}^\circ}$ of $(M^\circ)^{\oplus
  N}$. Let $\overline{q^\circ_i}$
  be the image of $q_i^\circ$ on $|C^\circ|$,
   and we set
  \[
    \eta^\circ := (|\pi^\circ|:|C^\circ|\to B^\circ; \overline{q_1^\circ} ,\ldots, \overline{q_{n+1}^\circ} ; M^\circ;\vec{{f}^\circ}),
  \]
which is a  family of prestable quasimaps to $\mathbb
  P^{N-1}$.
 
  We claim that, for any $B$-family of prestable quasimaps
  \begin{equation}
    \label{eqn_extension}
    \eta = (|\pi|:|C|\to B; \overline{q_1} ,\ldots, \overline{q_{n+1}}; M;\vec{{f}})
  \end{equation}
  extending $\eta^\circ$, there is unique prestable extension $\xi$ of $\xi^\circ$ whose
  underlying family of quasimaps is $\eta$.  Indeed, given $\eta$, by a standard argument (see, for example, \cite[Theorem
  1.5.1]{AJ} or \cite[Thoerem 4.1.7]{Clader}), there is a unique way (possibly after a finite base change) to add stack structure
  at nodes and markings and extend $L^\circ$ to a representable $d$-th root $L$  with
  \[
    L^{\otimes d} \cong \omega_{\pi,\mathrm{log}}\otimes M^{\vee}.
  \]
  Then, as in \cite[Theorem 5]{Z}, there is a unique
  extension of $(N^\circ,v_1^\circ,v_2^\circ)$ to $(N,v_1,v_2)$ such that $v_1$
  and $v_2$ have no common zeros. 

  The stability condition for $\xi$ can be reformulated in terms of $\eta, v_1$, and $v_2$
  as follows: $\xi$ is stable if and only if $(|\pi|:|C|\to B; \overline{q_1}
  ,\ldots, \overline{q_{n+1}})$ is a
  family of prestable $(n+1)$-pointed curves and 
  for each geometric fiber,
  \begin{itemize}
  \item $\vec{f}$ does not vanish at the nodes;
  \item for each zero $q$ of $\vec{{f}}$, we have
    \[
      \mathrm{ord}_q(\vec f) < \frac{1}{\epsilon};
    \]
  \item the $\mathbb Q$-line bundle
    \[
      M^{\otimes \epsilon}\otimes \omega_{|C|/B}([\overline{q_1}]+\cdots [\overline{q_{n+1}}])
    \]
    is ample;
  \item when $v_1=0$, $\vec f$ does not vanish at $q_1$;
  \item when $v_2=0$, the $\mathbb Q$-line bundle
    \[
      M^{\otimes \epsilon}\otimes
      \omega_{|C|/B}(\delta[\overline{q_1}]+[\overline{q_2}]+\cdots [\overline{q_{n+1}}])
    \]
    is ample.
  \end{itemize}
  We will show that, up isomorphism, there is a unique $\eta$ satisfying the above
  conditions.

  As in the proof of \cite{Z}, 
  we may assume that  $\pi^\circ$ is smooth and both $v^\circ_1$ and
  $v_2^\circ$ are non-zero.
  We first consider the case
  \begin{equation}
  \label{eq:case1}
     2g-2+n+\epsilon \deg M + \delta > 0.
  \end{equation}
  When viewing $\overline{q_1}$ as a light marking, 
  the generic fiber is $\epsilon$-stable. Possibly after finite base change, we extend the
  quasimap $\eta^\circ$ to an $\epsilon$-stable $B$-family with light marking
  $\overline{q_1}$, and we get a unique extension $(N,v_1,v_2)$ such that $v_1$ and $v_2$
  have no common zeros. The only situation that violates the master-space
  stability condition is if
  \begin{equation}
    \label{eqn_bad_case_1}
    v_1(b) = 0 \quad \text{and} \quad \overline{q_1} \text{ is a zero of $\vec f$ in 
      the special fiber}.
  \end{equation}
  If this happens, we blow up the total space of the family at the marking $\overline{q_1}$ of the
  special fiber. It is easy to see that $\vec {f^\circ}$ uniquely extends to a $B$-family of
  prestable quasimaps from the new family of curves to $\mathbb P^{N-1}$. 
  
  We repeat this
  procedure until (\ref{eqn_bad_case_1}) does not hold.  As in 
 \cite{Z} , the vanishing order of $v_1$ at $b$ drops
  by one after each step, so the procedure eventually terminates.
  Finally, we contract the exceptional divisors on which $M$ is trivial. 
  Note that $M$ has nontrivial degree on the exceptional divisor of the last
  blowup, which contains $\overline{q_1}$ and is not contracted. Thus, it is easy to see that
  after contractions, (\ref{eqn_bad_case_1}) still does not hold, so we do get a stable family in the case of \eqref{eq:case1}.

Suppose, now, that
    \begin{equation}
    \label{eq:case2}
       2g-2+n +\epsilon \deg M +  \delta \leq 0.
    \end{equation}
    Since we have assumed that $2g-2+n \geq -1$ and $\delta$ is sufficiently
    small, \eqref{eq:case2} only holds when $g=0$, $n=1$, and $\epsilon \deg M<1$.
    In this case, we can find a $B^\circ$-isomorphism between $|C^\circ|$ and
    $\mathbb P^1\times B^\circ$, identifying $\overline{q^\circ_1}$ with $\{0\}\times B^\circ$, $\overline{q^\circ_2}$
    with $\{\infty\}\times B^\circ$, and
    $v^\circ_1/v^\circ_2$ with the standard tangent vector $\partial/\partial
    z$, where $z$ is
    the coordinate on $\mathbb P^1$. We first take the constant family as the
    extension of the curves, and the markings and the prestable quasimaps extend
    uniquely.
     The only
    situation that violates the master-space
    stability condition is when 
  \begin{equation}
    \label{eqn_generic_unstable_1}
      \text{ $\overline{q_2}$ is a zero of $\vec f$ in the special fiber}.
  \end{equation}
    If this happens, we modify the underlying curves by repeated blowups at $q_2$
    of the special fiber until (\ref{eqn_generic_unstable_1})
    no longer holds. Finally, we blow down the unstable components. This produces
    a stable reduction. Note that the irreducible component of the special fiber
    containing $\overline{q_1}$ needs to be blown
    down precisely when its $M$-degree is zero, which takes
    $v_1$ to $0$ and $\overline{q_1}$ to some point where $\vec f$ does not
    vanish. This completes the construction in the case of \eqref{eq:case2}.
  
    The uniqueness part is standard, so we omit the details. The key is that,
    if $C^\prime$ is a smooth rational subcurve containing $\overline{q_1}$ and only one
    other special point, then $C^\prime$ cannot be contracted unless $\deg
    M|_{C^\prime}=0$, since contracting $C^\prime$ takes $v_1$ to $0$.
\end{proof}

\subsection{Virtual cycle for the master space}
Let $\widetilde X^{\epsilon,\delta}_{g,(1/d), \beta}$ be the master-space analogue of
$X^{\epsilon}_{g,(1/d), \beta}$.  More precisely, an $S$-point of $\widetilde
X^{\epsilon,\delta}_{g,(1/d), \beta}$ is
\[
  (\pi:C\to S;q_1;L,N; \vec p,v_1,v_2 ) \in  \widetilde Z^{\epsilon,\delta}_{g,(1/d),\beta}(S)
\]
together with a section 
\[
  \vec x = (x_1,\cdots,x_M) \in \Gamma(C,{\textstyle\bigoplus_{i=1}^M} L^{\otimes w_i}).
\]
As before, $\widetilde X^{\epsilon,\delta}_{g,
(1/d), \beta}$ admits a perfect obstruction
theory relative to 
the smooth Artin stack $\widetilde D_{g,(1/d),\beta}$ parametrizing only
$(C;q_1;L,N;v_1,v_2)$. The relative perfect obstruction theory is
defined by the same formula as (\ref{eq:pot}), and the marking is narrow. The same formula as before defines a cosection whose
degeneracy locus is $\widetilde Z^{\epsilon,\delta}_{g,(1/d),\beta}$, and we
get a cosection-localized virtual
class
\[
  [  \widetilde Z^{\epsilon,\delta}_{g,(1/d),\beta} ]^{\mathrm{vir}} \in A_*(  \widetilde Z^{\epsilon,\delta}_{g,(1/d),\beta} ).
\]

\subsection{Localization on the master space}
From now on, we assume that $g\geq 1$.

Define a $\mathbb C^*$-action on $\widetilde X^{\epsilon,\delta}_{g,(1/d), \beta}$ by\footnote{We note that this action is opposite to the one in \cite{Z}.}  
\begin{align*}
  t\cdot (\pi:C\to S; q_1;L,N; \vec p,\vec x,v_1,v_2 )  =
  (\pi:C\to S;q_1;L,N; \vec p,\vec x,t^{-1}v_1,v_2 ),\quad t\in \mathbb C^*.
\end{align*}
The perfect obstruction theory is equivariant and the cosection is invariant, and the restriction of the perfect obstruction theory to each
fixed locus has a global resolution. 
The degeneracy locus $\widetilde Z^{\epsilon,\delta}_{g,(1/d),\beta}$ of the cosection has
three types of fixed loci:
\begin{enumerate}
\item $F_0\widetilde Z$ is the vanishing
  locus of $v_1$.
\item $F_\infty\widetilde Z$ is the vanishing locus
  of $v_2$.
\item For each $0<\beta^\prime< 1/\epsilon$, $F_{\beta^\prime}\widetilde
  Z$ is the locus where
  \begin{itemize}
  \item $C=C_g \cup C_0$, where $C_0$ is a smooth rational subcurve and $\deg((L^{\otimes -d} \otimes \omega_{\log})|_{C_0}) = \beta^\prime$;
  \item $q_1\in C_0$ and $C_g\cap C_0$ are the only two special points of $C_0$;
    \item neither $v_1$ nor $v_2$ is zero;
  \item $\vec p$ has vanishing order $\beta^\prime$ at $q_1$.
  \end{itemize} 
\end{enumerate}
Set-theoretically, these are the only fixed loci.
Indeed, if $v_1,v_2
\neq 0$, then the only way that an automorphism can re-scale $t^{-1}v_1$ back to
$v_1$ is by acting nontrivially on the component $C_0$ of $C$ on which $q_1$
lies. This changes the moduli point unless $C_0$ has just one other special point
and all of $\deg(L^{\otimes -d}\otimes {\omega_{C,\log}}|_{C_0})$ is concentrated in a single basepoint at $q_1$,
which forces $C_0$ to have genus zero, as claimed.

We now describe the stack structure of the fixed loci and compute their contributions to the
localization formula. The first two cases are the same as in \cite[Lemma
7.15]{Z}, so we only state the results.   
The key to the proof is that $F_0\widetilde Z$ and $F_\infty \widetilde Z$ are effective
Cartier divisors defined by the vanishing of $v_1$ and $v_2$, respectively. 
\begin{lemma}
\label{lem:fixedloci1}
  The substack $F_0\widetilde Z$ is
  isomorphic to $Z^{\epsilon}_{g,(1/d),\beta}$, where $q_1$ is
  viewed as a heavy marking. Its localization contribution is
  \[
   \frac{[F_0\widetilde Z]^{\mathrm{vir}}} {e_{\mathbb
        C^*}(N^{\mathrm{vir}}_{F_0\widetilde Z})}
    =
    \frac{ [ Z^{\epsilon}_{g,(1/d),\beta} ]^{\mathrm{vir}}}{- z
      -\psi_{1}}.
  \]
  The substack $F_\infty\widetilde Z$ is
  isomorphic to $Z^{\epsilon,\delta}_{g, (1/d),\beta}$, where $q_1$ is
  viewed as a light marking. Its localization contribution is
  \[
    \frac{[F_\infty\widetilde Z]^{\mathrm{vir}}} {e_{\mathbb
        C^*}(N^{\mathrm{vir}}_{F_\infty\widetilde
          Z})}
    =
    \frac{ [ Z^{\epsilon,\delta}_{g,(1/d),\beta} ]^{\mathrm{vir}}}{z
      +\psi_{1}}.
  \]
  Here, $\psi_1$ is the $\psi$-class of the coarse curves at $q_1$.
\end{lemma}

We now come to $F_{\beta^\prime}\widetilde Z$. Recall that  the multiplicity of $L$ at 
$q_1$ is $1/d$, so the compatibility condition (\ref{eq:compatibility}) applied
to $C_g$ implies that the multiplicity of $L|_{C_g}$ at $C_g\cap C_0$ is
\[m_{\beta^\prime}:=\left\langle {\frac{\beta^\prime+1}{d}} \right\rangle.\]
Thus, the
automorphism group of the node $C_g\cap C_0$
is cyclic or order $d_m: = d/\gcd(d\cdot m_{\beta^\prime},d)$. 
Since the node  is balanced, the multiplicity of $L|_{C_0}$ at $C_g\cap C_0$ is
$m_\bullet := \left\langle {-m_{\beta^\prime}} \right\rangle$.

We claim, in fact, that $F_{\beta'}\widetilde Z$ is isomorphic to a fiber product of the
moduli space $Z^{\epsilon}_{g,(m_{\beta^\prime}),\beta-\beta'}$ with the locus
$F^{\epsilon}_{\beta'}$ in the graph space
$\mathcal{G}Z^{\epsilon}_{0,1,\beta'}$. Recall from the definition of the
$J$-function that $F^{\epsilon}_{\beta^\prime}\subset \mathcal
GZ^{\epsilon}_{0,1,\beta^\prime}$ is defined as the fixed locus of the graph space where
the only marking lies at $\infty$ and all of the degree lies over $0$. 
The multiplicity at the marking must be $m_\bullet$.
In the definition of the $J$-function, the virtual normal bundle $N^{\vir}_{F_{\beta^\prime}^{\epsilon}/\G
\Ze_{0,1,\beta^\prime}}$ has two parts:
\begin{enumerate}
\item
  a rank-$1$ part coming from deforming the marking away from $\infty\in \mathbb
  P^1$, whose Euler class is $( -z )$; 
\item the moving part of the relative obstruction theory (\ref{eq:E}), which we
  denote by $N^{\vir,\mathrm{rel}}_{F_{\beta^\prime}^{\epsilon}/\G
\Ze_{0,1,\beta^\prime}}$.
\end{enumerate}

We now form a morphism
\begin{equation}
  \label{map:iota-beta-prime}
  \iota_{\beta^\prime} : Z^{\epsilon}_{g,(m_{\beta^\prime}),\beta-\beta^\prime}
  \times_{\overline{\mathcal  I}Z}  F^{\epsilon}_{\beta^\prime} \to
  F_{\beta^\prime}\widetilde Z.
\end{equation}
Here, let us
denote the marking, $\psi$-class, and evaluation map of
$Z^{\epsilon}_{g,(m_{\beta^\prime}),\beta-\beta'}$ by
\[
q_1^\prime,\psi_1^\prime,\; \text{and }  \ev_1^\prime,
\]
and denote by
\[
  \widehat{\mathrm{ev}}_{\bullet} : F^{\epsilon}_{\beta^\prime} \to \overline{\mathcal  I}Z
\]
the evaluation map of $F^{\epsilon}_{\beta'}$ at the unique marking $\infty$
with the inverted banding.
The fiber product in \eqref{map:iota-beta-prime} is defined via $\ev_1'$ and
$\widehat{\ev}_{\bullet}$.

To define
$\iota_{\beta^\prime}$, take any $S$-point of $Z^{\epsilon}_{g,(m_{\beta^\prime}),\beta-\beta^\prime}
\times_{\overline{\mathcal  I}Z}  F^{\epsilon}_{\beta^\prime}$, which  consists of 
\begin{equation}
  \label{eq:gluing-data}
  \begin{aligned}
    &\eta_g = (\pi_{g}:C_g \to S, q^\prime_1;L_{g}, \vec p_{g})\in
      Z^{\epsilon}_{g,(m_{\beta^\prime}),\beta-\beta^\prime}(S), \\
      &\eta_0 = (\pi_{0}:C_0\to S,q_\bullet;L_0; \vec p_{0})\in F^{\epsilon}_{\beta^\prime}(S) ,\text{ and }\\
   &   \text{a }2\text{-morphism     } 
  \theta: \mathrm{ev}^\prime_{1} \overset{\simeq}{\to} \widehat{\mathrm{ev}}_{\bullet}.
  \end{aligned}
\end{equation}
By the definition of $F^{\epsilon}_{\beta^\prime}$, $C_0$
is $\mathbb P^1\times S$ with additional stack structure along the marking
$q_\bullet = \{\infty\}\times S$. The idea is to glue $q_1^\prime$ with
$q_\bullet$ and place the light marking $q_1$ at $0\in \mathbb P^1$.
Let
\begin{equation}
  \label{eq:tail-stackify}
  \rho: C_0^\prime \to C_0
\end{equation}
be the $d$-th root stack of the divisor $\{0\}\times S$ with universal root
$q_1\subset C_0^\prime$.
We also view $q_\bullet$  as a marking on $C_0^\prime$. Note that,
as relative Cartier divisors, $d[q_1]=\rho^*(\{0\}\times S)$.

Let
$L^\prime_0  = \rho^*L_0\otimes \mathcal  O_E([q_1])$.
Then
\begin{equation}
  \label{eq:P-compare}
  \begin{aligned}
    \rho_*\left( L^\prime_0 \right) = L_0 \quad \text{and} \quad 
    \rho_*\left({L^\prime_0}^{\otimes -d}\otimes \omega_{C^\prime_0/S,\mathrm{log}}\right )  = 
    {L_0}^{\otimes -d}\otimes \omega_{C_0/S,\mathrm{log}}.
  \end{aligned}
\end{equation}
Via the second isomorphism, $\vec p_0$ induces a section $\vec{p^\prime_0} \in \Gamma(({L^\prime_0}^{\otimes -d} \otimes
\omega_{C_0^\prime/S,\log})^{\oplus N})$. Thus, we get a family of genus-$0$ Landau--Ginzburg quasimaps
\[
  \eta^\prime_0 = (\pi_0\circ\rho: C_0^\prime \to S; q_1, q_\bullet;L^\prime_0;\vec {p^\prime_0}).
\]
Note that the evaluation map of $\eta^\prime_0$ at $q_\bullet$ is also $\mathrm{ev}_\bullet$.
The map $\iota_{\beta^\prime}$ is defined by gluing the Landau--Ginzburg quasimaps
$\eta_g$ and $\eta^\prime_0$ along the markings $q^\prime_{1}$ and $q_\bullet$ and setting
$N=\mathcal  O_{S}$, $v_2\equiv 1$, and $v_1 \equiv
\partial/\partial z$, where $z$ is the coordinate on $\mathbb P^1$. 
This defines a morphism to $\widetilde Z^{\epsilon,\delta}_{g,(1/d),\beta}$ that maps
each closed point to $F_{\beta^\prime}\widetilde Z$.  Moreover, the resulting family 
 is $\mathbb C^*$-invariant, so it factors through
$F_{\beta^\prime}\widetilde Z$.

\begin{lemma}
  \label{lem:contr1}
  The morphism $\iota_{\beta^\prime}$ is an isomorphism onto $F_{\beta^\prime}\widetilde
  Z$.   
\end{lemma}
\begin{proof}
  We construct a quasi-inverse. Given any family
  \begin{equation}
    \label{eq:glued-family}
    (\pi: C\to S; q_1;L, N; \vec p,v_1,v_2)\in \widetilde
    Z^{\epsilon,\delta}_{g,(1/d),\beta}(S)
  \end{equation} in
  $F_{\beta^\prime}\widetilde Z$, we must recover the gluing
  data (\ref{eq:gluing-data}). By the definition of $F_{\beta^\prime}\widetilde Z$, over any closed
  point $s\in S$, (\ref{eq:glued-family}) comes from gluing the data
  (\ref{eq:gluing-data}), and the key is to show that we can split the node in families over any (possibly non-reduced)
  base $S$.
  
  Since the family is $\mathbb C^*$-fixed, $\mathbb C^*$ acts on $C$ and the
  morphism $\pi:C\to S$ is invariant. On the  fiber $C_s = C_{g,s}\cup C_{0,s}$
  over each closed point $s\in S$,
  $\mathbb C^*$ acts on $C_{0,s}$ fixing $q_1$ and the node
  $q_{n}:= C_{g,s}\cap C_{0,s}$, and $t\in \mathbb C^*$ maps the
  tangent vector $t^{-1}v_1/v_2$ back to $v_1/v_2$.
  Hence, $\mathbb C^*$ acts nontrivially on the tangent space to $C_{0,s}$ at
  $q_{\mathrm{n}}$. It is obvious that $\mathbb C^*$ acts trivially on $C_{g,s}$.
  Hence it acts nontrivially on the first order deformation smoothing the node
  $q_{\mathrm{n}}$.
  Since the family is fixed by $\mathbb C^*$, the $S$-family of curves $C$ can be
  decomposed as $C_g$ and $C_0$ glued at a pair of markings. 

  Now it is clear how to recover (\ref{eq:gluing-data}): restricting the Landau--Ginzburg quasimaps
  to $C_g$ and $C_0$ separately, we recover $\eta^\prime_0$ and $\eta_g$.  It is easy to
  see that $\eta^\prime_0$ is equivalent to $\eta_0$, and this defines a quasi-inverse to $\iota_{\beta^\prime}$.
\end{proof}

Let $\mathrm{pr}_1: F_{\beta^\prime}\widetilde Z \to Z^{\epsilon}_{g,
  (m_{\beta^\prime}),\beta-\beta^\prime}$ and $\mathrm{pr}_2: F_{\beta^\prime}\widetilde Z \to
F^{\epsilon}_{\beta^\prime}$ be the projections. 

\begin{lemma}
  \label{lem:contr2}
  We have
  \[
    [ F_{\beta^\prime}\widetilde{Z} ]^{\mathrm{vir}} = \mathrm{pr}_1^*[ Z^{\epsilon}_{g,(m_{\beta^\prime}),\beta-\beta^\prime}
    ]^{\mathrm{vir}},
  \]
  and the inverse of the equivariant Euler class of the virtual normal bundle
  restricted to $F_{\beta^\prime}\widetilde Z$
  is
  \[
    \frac{1}{e_{\mathbb C^*}( N^{\mathrm{vir}}_{F_{\beta^\prime}\widetilde
          Z})} =
    \mathrm{pr}_1^*\Big(  \frac{d_m}{-\psi^\prime_{1} - z}\Big) \cdot \mathrm{pr}_2^* \Big(   \frac{1}{e_{\mathbb C^*}
        { (N^{\vir,\mathrm{rel}}_{F_{\beta^\prime}^{\epsilon}/\G
            \Ze_{0,1,\beta^\prime}})
        } }\Big).
  \]

\end{lemma}
\begin{proof}
  The proof is similar to that of Lemma~\ref{lem:vclass}. 

  Let $\mathfrak Z\subset \widetilde D_{g,(1/d),\beta}$ be the reduced,
  locally-closed substack where $q_1$ is on a rational tail of degree $\beta^\prime$ and
  $v_1$ and $v_2$ are both nonzero. The normal bundle of $\mathfrak Z$ is
  moving.
  As in the proof of
  Lemma~\ref{lem:vclass}, $[ F_{\beta^\prime}\widetilde Z ]^{\mathrm{vir}}$ can
  be defined by the fixed part of the absolute obstruction theory induced by
  \begin{equation}
    \label{eq:obs-thy-appd}
       \left(R\pi_*\left({\textstyle\bigoplus_{i=1}^M}
        \left(\L^{\otimes w_i}\left(-\textstyle \sum_{k=1}^n\Delta_k\right)\right)
        \oplus {\textstyle \bigoplus_{j=1}^N}\cP\right) \right)^\vee,
  \end{equation}
  which is a relative perfect obstruction theory relative to $\mathfrak Z$.

  The universal curve decomposes as $\mathcal C_g\cup
  \mathcal C_0$, where $\mathcal C_0$ is the rational tail containing $q_1$. Let
  $\Delta_0\subset \mathcal C_0$ and $\Delta_g\subset \mathcal C_g$ be the
  node $\mathcal C_g\cap \mathcal C_0$. The contributions from the subsheaf
  $\bigoplus_{j=1}^N\cP|_{\mathcal C_0}(-\Delta_0)$
  and the quotient sheaf
  $ \bigoplus_{i=1}^M \left(\L^{\otimes w_i}\left(-\textstyle
      \sum_{k=1}^n\Delta_k\right)\right)|_{\mathcal C_0}$ to \eqref{eq:obs-thy-appd} are both moving, and this moving part is identified with
  $\mathrm{pr}_2^*(N^{\vir,\mathrm{rel}}_{F_{\beta^\prime}^{\epsilon}/\G
    \Ze_{0,1,\beta^\prime}})$ thanks to the relation (\ref{eq:P-compare}).
  The remaining part is the fixed part, given by 
  \[
       \left(R\pi_*\left(\textstyle{\bigoplus_{i=1}^M}
        \left(\L^{\otimes w_i}\left(-\Delta_g- {\textstyle
            \sum_{k=1}^n}\Delta_k\right)|_{\mathcal C_g}\right)
        \oplus \textstyle{\bigoplus_{j=1}^N}\cP|_{\mathcal C_g}\right) \right)^\vee.
  \]
  This is exactly the  pullback of the compact-type perfect obstruction theory of
  $Z^{\epsilon}_{g,(m_{\beta^\prime}),\beta-\beta^\prime}$ relative to $D_{g,(m_{\beta^\prime}),\beta-\beta^\prime}$.
  Moreover, the cosections are compatible and the forgetful morphism $
  \mathfrak Z\to D_{g,(m_{\beta^\prime}),\beta-\beta^\prime}$ is \'etale. 
  Thus, the identity of virtual cycles follows from cosection-localized pullback. The moving part of the tangent complex of
  $\widetilde D_{g,(1/d),\beta}$ along $\mathfrak Z$ is the normal bundle of $\mathfrak Z$, whose
  Euler class is $\frac{1}{d_m}(-\psi^\prime_{1}-z)$.
\end{proof}

Having established the contributions from each fixed locus, we can collect all of them to get the localization formula on the master space.  In particular, there is a morphism
\[\rho: \widetilde Z^{\epsilon,\delta}_{g,(1/d),\beta} \rightarrow Z^{\epsilon,\delta}_{g,(1/d),\beta}\]
that forgets $N,v_1$, and $v_2$ and stabilizes as necessary, and we calculate
the coefficient of $z^{-2}$ in the class $\rho_*(\psi_1\cap [ \widetilde
Z^{\epsilon,\delta}_{g,(1/d),\beta}
  ]^{\mathrm{vir}})$ by $\C^*$-localization. This coefficient is zero, but
expressing it in terms of contributions from each fixed locus yields a
nontrivial relation. To calculate these contributions, we first note that, under
the isomorphism $F_0\widetilde{Z} \cong Z^\epsilon_{g,(1/d),\beta}$ from
Lemma~\ref{lem:fixedloci1}, the restriction of $\rho$ to $F_0\widetilde{Z}$ is the
morphism
\[
  c: Z^\epsilon_{g,(1/d),\beta} \to  Z^{\epsilon,\delta}_{g,(1/d),\beta}
\]
that replaces the heavy marking with a light marking and then stabilizes. Under
the isomorphism $F_{\infty}\widetilde{Z} \cong Z^{\epsilon,
\delta}_{g,(1/d),\beta}$, the restriction of $\rho$ to $F_{\infty}\widetilde Z$ is the
identity. Finally, under the isomorphism $F_{\beta'}\widetilde{Z}$ from
Lemma~\ref{lem:contr1}, the restriction of $\rho$ to $F_{\beta'}\widetilde{Z}$ is
the composition of $\mathrm{pr_1}$ with 
the morphism
\[
  b: Z^\epsilon_{g,(m_{\beta^\prime}),\beta-\beta^\prime} \to Z^{\epsilon,\delta}_{g,(1/d),\beta}
\]
that replaces the last marking with a light marking that is also a basepoint of order
$\beta^\prime$.

Using Lemma~\ref{lem:fixedloci1}, Lemma~\ref{lem:contr1}, and Lemma~\ref{lem:contr2} together with the projection formula, the relation that results from localization is as follows:
\begin{equation} 
  \label{eq:dilaton-wc}
 \psi_1\cap  [Z^{\epsilon,\delta}_{g,(1/d),\beta} 
  ]^{\mathrm{vir}}=    c_*\left( \psi_{1} \cap [ Z^{\epsilon}_{g,(1/d),\beta}
      ]^{\mathrm{vir}}\right) 
    - \sum_{\beta^\prime}b_*\left(
       ( \mathrm{ev}^\prime_{1})^*\left( \mu^\epsilon_{\beta^\prime}(-\psi^\prime_{1}) \right) \cap
      [ Z^\epsilon_{g,(m_{\beta^\prime}),\beta-\beta^\prime} ]^{\mathrm{vir}}\right).
\end{equation}

Another relation is obtained by calculating $\rho_*(\ev_{1}^*(H) \cap [ \widetilde
    Z^{\epsilon,\delta}_{g,(1/d),\beta} ]^{\mathrm{vir}})$, where $H$ is the ``hyperplane class'' on $[\mathbb C^{N}/\mathbb C^*]\subset \ri
Z$.  Concretely, $\ev_{1}^*(H)$ is the Euler class of the line bundle
$L^{\otimes -d}\otimes \omega_{C,\log}|_{q_{1}}$, and the restriction of $\ev_{1}^*(H)$ to
$F_{\beta^\prime}\widetilde Z$ is equal to ${(\ev_{1}^{\prime})}^*(H) + \beta^\prime z$.  The coefficient of $z^{-1}$ in this pushforward is zero, but expressing it in terms of contributions from each fixed locus and applying the same reasoning as above yields:
\begin{equation} 
  \label{eq:divisor-wc}
  \begin{aligned}
  \ev^*_{1}(H)\cap [Z^{\epsilon,\delta}_{g,(1/d),\beta} 
  ]^{\mathrm{vir}}     = &   c_*\left( \ev^*_{1}(H)\cap [ Z^{\epsilon}_{g, (1/d),\beta}]^{\mathrm{vir}}\right) \\
    +& \beta^\prime \sum_{\beta^\prime}b_*\left(
      (\mathrm{ev}^\prime_{1})^*\left(\mu^\epsilon_{\beta^\prime}(-\psi^\prime_1) \right)  \cap
      [ Z^\epsilon_{g,(m_{\beta^\prime}),\beta-\beta^\prime} ]^{\mathrm{vir}}\right)\\
    - &\sum_{\beta^\prime}b_*\left(
      (\mathrm{ev}^\prime_{1})^*\left(\frac{H}{\psi^\prime_1} \mu^\epsilon_{\beta^\prime}(-\psi^\prime_1)\right)  \cap
      [ Z^\epsilon_{g,(m_{\beta^\prime}),\beta-\beta^\prime} ]^{\mathrm{vir}}\right).
  \end{aligned}
\end{equation}
Here, in the expression $\frac{H}{\psi^\prime_1} \mu^\epsilon_{\beta^\prime}(-\psi^\prime_1)$, we simply
discard the term of $\mu^\epsilon_{\beta^\prime}(-\psi^\prime_1)$ that does not involve $\psi^\prime_1$.

\subsection{Wall-crossing without a heavy marking}
Equations \eqref{eq:dilaton-wc} and \eqref{eq:divisor-wc} can be viewed as wall-crossing formulas for converting a light marking to a heavy marking.  Equipped with them, we are prepared to prove the $n=0$ case of the main wall-crossing theorem.

From now on, we compute on $Z^{\epsilon}_{g,0,\beta}$ and suppress all the
pushforward notations.  In particular, the class that was denoted $\psi'_1$ in the previous subsection will now be denoted simply $\psi_1$, and similarly for $q_1'$ and $\ev_1'$.  The meaning of $\ev_i^*(H)$ and $\psi_i$ are determined by the virtual cycle that they are capped with.

First, suppose $N>1$.  We then have $\beta[ Z^{\epsilon}_{g,0,\beta} ]^{\mathrm{vir}}  =
   \ev_1^*(H)\cap [ Z^{\epsilon,\delta}_{g,(1/d),\beta}
]^{\mathrm{vir}}$ by Lemma~\ref{lem:forget-light}. Combining this with
Theorem~\ref{thm:main} and 
(\ref{eq:divisor-wc}), we get 
\begin{equation}
  \label{eq:using-divisor}
\begin{aligned}
  \beta[ Z^{\epsilon}_{g,0,\beta} ]^{\mathrm{vir}}   &=  \sum_{\beta_0+\cdots+\beta_k = \beta}
     \frac{1}{k!} \ev_1^*(H)
     \prod_{i=1}^k\ev_{1+i}^*(\mu_{\beta_i}^\epsilon(-\psi_{1+i})) \cap [
     Z^\infty_{g,(1/d)+k,\beta_0} ]^{\mathrm{vir}} \\
   +& \sum_{\beta_0+\cdots+\beta_k + \beta^\prime = \beta}
     \frac{\beta^\prime}{k!}
    \ev_1^*\left(\mu^\epsilon_{\beta^\prime}(-\psi_1)\right)
    \prod_{i=1}^k \ev_{1+i}^*(\mu_{\beta_i}^\epsilon(-\psi_{1+i})) \cap [
    Z^\infty_{g,(m_{\beta^\prime})+k,\beta_0} ]^{\mathrm{vir}} \\
   -&
    \sum_{\beta_0+\cdots+\beta_k + \beta^\prime = \beta}
    \frac{1}{k!}
    \mathrm{ev}_{1}^*\left(\frac{H}{\psi_1} \mu^\epsilon_{\beta^\prime}(-\psi_1)\right)
    \prod_{i=1}^k \ev_{1+i}^*(\mu_{\beta_i}^\epsilon(-\psi_{1+i}))\cap [
    Z^\infty_{g,(m_{\beta^\prime})+k,\beta_0} ]^{\mathrm{vir}}.
\end{aligned}
\end{equation}
Applying the divisor equation
  \begin{equation}
    \begin{aligned}
      \ev_1^*(H)\cdot \prod_{i=1}^k \psi^{a_i}_{1+i} \cap
      [Z^\infty_{g,(1/d)+k,\beta_0} ]^{\mathrm{vir}} 
      = & \beta_0\prod_{i=1}^k \psi^{a_i}_{i} \cap
      [Z^\infty_{g,k,\beta_0} ]^{\mathrm{vir}} \\
      &+ \sum_{j=1}^k\left(  \psi_{j }^{a_j-1} 
      \prod_{i\neq j}\psi_{i}^{a_i}\right)\cap
      [Z^\infty_{g,k,\beta_0} ]^{\mathrm{vir}}
    \end{aligned}
  \end{equation}
for $\infty$-theory, we then obtain
\begin{equation}
  \label{eq:divisor}
\begin{aligned}
 & \ev_1^*{(H)}\cdot \prod_{i=1}^k \ev_{1+i}^*(\mu_{\beta_i}^\epsilon(-\psi_{1+i}))
  \cap [Z^\infty_{g,(1/d)+k,\beta_0} ]^{\mathrm{vir}} 
  =  \beta_0\prod_{i=1}^k \ev_i^*(\mu_{\beta_i}^\epsilon(-\psi_{i}))
  \cap    [Z^\infty_{g,k,\beta_0} ]^{\mathrm{vir}} \\
   &+\sum_{j=1}^k\left( \ev_j^*(\frac{H}{\psi_{j}}\mu_{\beta_i}^\epsilon(-\psi_{j}))  \prod_{i\neq j} 
      \ev_i^*(\mu_{\beta_i}^\epsilon(-\psi_{i}))\right)\cap    [Z^\infty_{g,k,\beta_0} ]^{\mathrm{vir}}.
\end{aligned}
\end{equation}
Substituting (\ref{eq:divisor}) into (\ref{eq:using-divisor}) and rearranging the
summation using the symmetry of the marked points,  we get
\begin{equation}
\begin{aligned}
\beta[ Z^{\epsilon}_{g,0,\beta} ]^{\mathrm{vir}} =\beta \sum_{\beta_0+\cdots+\beta_k = \beta}
     \frac{1}{k!} 
     \prod_{i=1}^k\ev_{i}^*(\mu_{\beta_i}^\epsilon(-\psi_{i})) \cap [
     Z^\infty_{g,k,\beta_0} ]^{\mathrm{vir}},
\end{aligned}
\end{equation}
which is exactly the statement of the wall-crossing theorem in this case.

Now, suppose that $g>1$.  In this case, we have $(2g-2)[ Z^{\epsilon}_{g,0,\beta} ]^{\mathrm{vir}} =
    \psi_1\cap [ Z^{\epsilon,\delta}_{g,(1/d),\beta} ]^{\mathrm{vir}} $.
    Combining this with Theorem~\ref{thm:main} and
(\ref{eq:dilaton-wc}), we get 
\begin{equation}
  \label{eq:using-dilaton}
\begin{aligned}
  (2g-2)&[ Z^{\epsilon}_{g,0,\beta} ]^{\mathrm{vir}} 
  = \sum_{\beta_0+\cdots+\beta_k = \beta}
     \frac{\psi_1}{k!}
          \prod_{i=1}^k \ev_{1+i}^*(\mu_{\beta_i}^\epsilon(-\psi_{1+i}))
        \cap  [ Z^\infty_{g,(1/d)+k,\beta_0} ]^{\mathrm{vir}} \\
  & - \sum_{\beta_0+\cdots+\beta_k+\beta^\prime=\beta}
     \frac{1}{k!}
      \mathrm{ev}_{1}^*\left( \mu^\epsilon_{\beta^\prime}(-\psi_{1}) \right) 
       \prod_{i=1}^k \ev_{1+i}^*(\mu_{\beta_i}^\epsilon(-\psi_{1+i})) [
     Z^\infty_{g,(1/d)+k,\beta_0} ]^{\mathrm{vir}}.
\end{aligned}
\end{equation}
Using the dilaton equation for $\infty$-theory, we get
\begin{equation}
  \frac{\psi_1}{k!} \prod_{i=1}^k \ev_{1+i}^*(\mu_{\beta_i}^\epsilon(-\psi_{1+i}))
        \cap [ Z^\infty_{g,(1/d)+k,\beta_0} ]^{\mathrm{vir}} =
        \frac{2g-2+k}{k!} \prod_{i=1}^k \ev_{i}^*(\mu_{\beta_i}^\epsilon(-\psi_{i}))
        \cap  [ Z^\infty_{g,k,\beta_0} ]^{\mathrm{vir}} .
\end{equation}
Substituting this into (\ref{eq:using-dilaton}), we again obtain the desired wall-crossing formula.

The only remaining case is when $g=N=1$, but this can be handled exactly analogously to the work of Guo--Ross \cite{GR1} for the quintic threefold.

\bibliographystyle{abbrv}
\bibliography{biblio}

\end{document}